\documentclass{article}
\usepackage{amsmath,amssymb,amsfonts,amsthm}
\usepackage{graphicx} % Required for inserting images
\usepackage{hyperref}
\usepackage[capitalise]{cleveref}
\usepackage{enumitem}

\usepackage{algorithm}
\usepackage{algorithmic}
\usepackage{dsfont}
\usepackage{mathrsfs}  

\usepackage{graphicx}
\usepackage{epstopdf}
\usepackage{diagbox}

\newtheorem{theorem}{Theorem}[section]
\newtheorem{corollary}[theorem]{Corollary}
\newtheorem{lemma}[theorem]{Lemma}
\newtheorem{assumption}[theorem]{Assumption}
\newtheorem{example}[theorem]{Example}
\newtheorem{proposition}[theorem]{Proposition}
\newtheorem{definition}[theorem]{Definition}
\newtheorem{remark}[theorem]{Remark}
\newtheorem{Algorithm}[theorem]{Algorithm}
\counterwithin{theorem}{section}

\newcommand{\closure}{\operatorname{cl}}

\newcommand{\dd}{\mathrm{d}}

\newcommand{\interior}{\operatorname{int}}

\newcommand{\supp}{\operatorname{supp}}

\newcommand{\proj}{\pi}

\newcommand{\WOT}{\overset{\normalfont\textsc{wot}}{\to}}
\newcommand{\Id}{\mathrm{Id}}

\newcommand{\D}{D}

\DeclareMathOperator*{\esssup}{ess\,sup}
\DeclareMathOperator*{\essinf}{ess\,inf}

\newcommand{\red}{\mathrm{red}}

\newcommand{\sym}{\mathrm{sym}}
\newcommand{\Lip}{\mathrm{Lip}}

\newcommand{\DD}{\mathcal{D}}

\newcommand{\FF}{\mathcal{F}}
\newcommand{\GG}{\mathcal{G}}
\newcommand{\HH}{\mathcal{H}}
\newcommand{\II}{\mathcal{I}}
\newcommand{\JJ}{\mathcal{J}}

\newcommand{\MM}{\mathcal{M}}
\newcommand{\NN}{\mathcal{N}}

\newcommand{\QQ}{\mathcal{Q}}
\newcommand{\RR}{\mathcal{R}}

\newcommand{\TT}{\mathcal{T}}

\newcommand{\BBB}{\mathfrak{B}}
\newcommand{\bbb}{\mathfrak{b}}

\newcommand{\BBBB}{\mathscr{B}}
\newcommand{\LLLL}{\mathscr{L}}

\newcommand{\tol}{\texttt{tol}}

\renewcommand{\d}{\,\mathrm{d}}
\newcommand{\Qad}{Q_{\mathrm{ad}}}

\newcommand{\N}{\mathbb{N}}

\newcommand{\R}{\mathbb{R}}
\newcommand{\Z}{\mathbb{Z}}

\newcommand{\Com}[1]{}

\makeatletter
\newcommand{\oset}[3][0ex]{%
  \mathrel{\mathop{#3}\limits^{
    \vbox to#1{\kern-2\ex@
    \hbox{$\scriptstyle#2$}\vss}}}}
\makeatother

\newcommand{\weakly}{\rightharpoonup}
\newcommand{\weaklystar}{\stackrel\star\rightharpoonup}
\newcommand{\compactly}{\oset{\textup{c}}{\hookrightarrow}}

\usepackage{xcolor}
\definecolor{darkgreen}{rgb}{0,0.8,0}
\newcommand{\ccnote}[1]{\marginpar{\color{darkgreen}{CC: #1}}}

\title{
First- and second-order constrained optimisation for inverse homogenisation of unidirectionally varying microstructure
}
\author{Constantin Christof\footnote{Technical
University of Darmstadt, Department of Mathematics, Dolivostr.~15, 64293 Darmstadt, Germany, email: \texttt{christof@mathematik.tu-darmstadt.de}} \and Malte A.~Peter\footnote{University of Augsburg, Institute of Mathematics \& Centre for Advanced Analytics and Predictive Sciences (CAAPS),
Universitätsstr.~12a, 86159 Augsburg, Germany, email: \texttt{malte.peter@math.uni-augsburg.de}}
\and Michael Ulbrich\footnote{
Technical
University of Munich,
Department of Mathematics, School of Computation, Information and Technology, Boltzmannstr.~3,
85748 Garching b.\ München, Germany, email: \texttt{mulbrich@ma.tum.de}}
}
\date{\today}

\begin{document}

\maketitle

\begin{abstract}
    This paper is concerned
    with inverse homogenisation problems that
    aim to identify microstructure features of a  
    body from macroscopic measurements.
    The unknown microstructure is assumed to vary 
    macroscopically in one spatial direction
    and affects the 
    effective tensor in an 
    elliptic partial differential equation 
    that governs the 
    behaviour of the body on the macroscale.
    We establish the well-posedness of the 
    PDE-constrained minimisation 
    problems that model the task of identifying 
    the microstructure,
    derive 
    first-order necessary optimality conditions,
    and prove the convergence of a gradient-projection 
    algorithm for 
    the studied problem class in an infinite-dimensional 
    setting. 
    We moreover demonstrate that recently obtained 
    results on generalised differentiability 
    properties of 
    solution maps of bilateral obstacle problems 
    make it possible to set up and analyse a novel
    semismooth Newton method in the considered 
    situation.
    The functional analytic setting that we 
    use for 
    our investigation is chosen such that it covers 
    various elliptic PDEs simultaneously 
    (classical second-order equations, bi-Laplace problems, 
    linear elasticity, etc.) and explicitly exploits 
    the problem structure to obtain global convergence 
    guarantees for the studied
    solution algorithms in terms of the 
    Lipschitz constants of the involved operators. 
    Our analysis is accompanied by 
    numerical experiments which show that 
    the new semismooth Newton approach 
    requires far less expensive projection evaluations 
    than the gradient-descent method 
    and, thus, offers a significant speedup.
\end{abstract}


\bigskip
\textbf{MSC2020-Codes:} 35B27, 49M05, 49M15, 49M41, 49J52, 49K20.

\section{Introduction}
\label{sec:intro}
This paper is concerned with inverse 
homogenisation problems 
that are 
governed by elliptic partial differential equations (PDEs) 
and 
aim to
identify 
microstructure features, which are
periodic on the microscale but
unidirectionally varying on the macroscale,
by means of macroscopic measurements. 
A prototypical example of this type of inverse problem
is the following one 
for a two-dimensional 
elliptic PDE on a 
rectangular domain 
$\Omega = (0,L) \times (0,1)$, $L>0$,
that is microscopically 
perforated with a perforation radius varying in $x_1$-direction: 
\begin{equation}
\label{eq:model_problem_elliptic}
\left.
	\begin{aligned}
		\text{Minimise} 
		\quad &  
        \frac12 \|u - u_\mathrm{obs}\|_{L^2(\Omega)}^2 + 
        \frac{\nu}{2} \|q\|_{H^1(0,L)}^2 \\
        \text{w.r.t.}
        \quad &u \in H_D^1(\Omega), \quad  q \in H^1(0,L),\\
		\text{s.t.}
		\quad & 
        -\nabla\cdot (B(q)\nabla u)=f \text{ in } \Omega,
\quad
u=0 \text{ on } \Gamma_D,
\quad 
\vartheta \cdot B(q)\nabla u = 0 \text{ on } \Gamma_N,
  \\
		\text{and} \quad & a \leq q \leq b \text{ a.e.\ in } (0,L).
\end{aligned}
\quad \right \}
\end{equation}
Here, 
the boundary $\partial \Omega$ of $\Omega$
is assumed to be decomposed disjointly into a nonempty 
relatively open Neumann part
$\Gamma_N \subset \partial \Omega$ and a
closed Dirichlet part $\Gamma_D \subset \partial \Omega$ 
of positive surface measure;  
$\vartheta$ denotes the outer unit normal vector field of $\Omega$;
$u_\mathrm{obs} \in L^2(\Omega)$ is a given function representing 
data measured, e.g., in an experiment; 
$\nu > 0$ is a Tikhonov regularisation parameter;
$u$ denotes the PDE solution; 
$q \in H^1(0,L) \hookrightarrow L^\infty((0,L) \times (0,1))$
models the varying radius of the microperforation 
(i.e.\ $q(x_1)$ is the radius of the microperforation on the 
slice $\{x_1\} \times (0,1)$ of $\Omega$); 
$0 < a < b < 0.5$ are constant bounds on the perforation radius;
$f \in L^2(\Omega)$ is given; 
$B\colon [a,b] \to \R^{2\times 2}$ is a matrix field such that 
$B(r) \in \R^{2\times 2}$ is equal to the effective tensor
that is obtained when the cell-problems in a classical periodic 
homogenisation approach are solved for a periodicity 
cell $Y := (0,1)^2$ that has a hole of radius $r \in [a,b]$ centred at $(1/2, 1/2)$;
the spaces $L^2(\Omega)$, $ H_D^1(\Omega) := \{\phi \in H^1(\Omega) \mid 
\phi = 0 \text{ on } \Gamma_D\}$, and
$H^1(0,L)$ 
are defined as usual; and the governing PDE is understood weakly. 
Note that, from the viewpoint of physics,  
\eqref{eq:model_problem_elliptic} 
can be understood as a model 
for the problem of identifying
how the microperforation of a microscopically 
perforated elastic membrane (occupying $\Omega$,
clamped on $\Gamma_D$, and free on $\Gamma_N$) 
varies in $x_1$-direction by pushing 
onto the membrane with a force $f$ and observing 
the resulting out-of-plane displacement $u_\mathrm{obs}$. 
The $\nu$-term serves as a regularisation in this context that 
is supposed to compensate noise in the measured deformation data;
see \cite{Isak:1998} for an introduction to the related field of inverse problems. 
The main purpose of this work is to develop a theoretical 
and algorithmic framework for the analysis and numerical 
solution of inverse homogenisation problems like 
\eqref{eq:model_problem_elliptic}. 
More precisely, we discuss 
the existence of minimisers,
first-order necessary optimality conditions, 
the convergence of an infinite-dimensional
gradient projection method, 
and a novel semismooth Newton 
approach for inverse homogenisation problems that aim to identify
macroscopically unidirectionally varying 
microscopic features from macroscale measurements. 
A main feature 
of our analysis is that it is conducted 
in a general functional analytic setting that
covers not only the model application \eqref{eq:model_problem_elliptic}
but also numerous other inverse homogenisation problems
involving various ansatz spaces for $q$, objective functions, 
and partial differential equations; see
\cref{sec:abstractsetting}. 
For the main results of our work, 
we refer the reader to \cref{th:solvabilityP,th:first_order_necessary_condition,th:convgradproj,th:convergence_Newton}.
Before we discuss our findings in more detail, 
let us put our work into perspective: 

Materials specifically designed for certain purposes often owe their properties to microstructural features like 
fibre-reinforcements or some form of porosity. 
Particularly in fabricated materials, the microstructure is often periodic (up to production quality), to simplify the production process
or to leverage effects caused by periodic microstructures specifically. 
When it comes to designing paths of wave propagation in materials,
for example, periodic microstructures are particularly relevant;
see \cite{Davies_Roadmap:2025} for an overview of the state-of-the-art and current challenges in this field.

For the durability assessment and/or quality inspection of structures whose properties arise from microscopic features, it is often of vital importance to be able to detect and identify defects and macroscopic variations of microstructural properties (fibre-densities, perforation patterns, etc.).
As microscopic defects can rarely be measured directly, 
these types of variations often have to be identified 
based on the macroscopic behaviour that they induce, for instance, 
altered load-bearing capabilities or a changed response to mechanical excitations. Periodic microstructures that
are slowly modulated macroscopically in only one spatial direction
are of particular relevance in this context.
They arise, for example, in  additive manufacturing processes that rely on a layer-wise printing procedure, 
or when laminates are considered. In metamaterial applications, unidirectionally slowly modulated  microstructures (typically called \emph{graded} in this context) are also often leveraged to increase the bandwidth of desired effects, e.g.~rainbow trapping or rainbow reflection \cite{Ben:2019,DaPonti:2021}.%\pagebreak

Mathematically,
the task of identifying
microstructure features from 
macroscopic measurements gives rise to
\emph{inverse homogenisation problems} like
\eqref{eq:model_problem_elliptic}, that is, 
optimisation problems governed by partial differential equations which aim to minimise an objective function depending on macroscopic quantities 
with respect to a parameter modelling material properties
on the microscale. 
The effects induced by the 
microstructures are typically 
incorporated via upscaled equations obtained 
from periodic homogenisation approaches in these problems, 
that is, the effective parameters in 
the governing PDE (in the case of the model problem \eqref{eq:model_problem_elliptic}, the field $B(q)$) are obtained with the help of solutions to auxiliary problems formulated on unit (periodicity) cells. In total, this gives rise to a rich and intricate 
structure, 
in which the solution operators of the cell problems 
(in \eqref{eq:model_problem_elliptic} encoded in
the mapping $[a,b] \ni r \mapsto B(r) \in \R^{2\times 2}$) 
act by superposition on the parameter describing the 
microstructure variation (in \eqref{eq:model_problem_elliptic} the function $q \in H^1(0,L)$) to give rise to the effective tensor field 
in the governing equation (in \eqref{eq:model_problem_elliptic}
the field $B(q) \in L^\infty(\Omega;\R^{2\times 2})$); also see \cite{Pet07,Pet07b,GahPetPopWie25} for this structure in the context of evolving microstructures. 

In the literature, inverse homogenisation problems 
have been tackled from various angles in the past, albeit 
often in simplified settings, e.g.\ without the 
assumption that the microstructure is macroscopically 
varying. (In the context of the model problem  \eqref{eq:model_problem_elliptic}, this scenario 
corresponds to the situation that 
the function $q$ describing the radius of the microperforation 
is constant throughout $\Omega$.) 
See, for instance, 
\cite{Loch:2022},
where 
the main axes of ellipsoidal microinclusions 
in the field of linear elasticity 
are identified in a fully periodic setting 
without any macroscopic variations,
and \cite{Loch:2024} where the same problem 
is discussed in an instationary scenario. 
For further related works 
on inverse homogenisation problems, 
we also refer the reader to    
\cite{All:2019,Bendsoe1988,Elena_Cherkaev_2001,Haslinger1995,BriLegRug26,Lipton2006,Sch:2015,Sch:2016}
and the references therein. An aspect that is particularly challenging 
when studying inverse homogenisation problems 
with spatially varying microstructures is their effective numerical 
solution. The reason for this is that the homogenisation 
background as well as the mathematical requirements necessary 
for ensuring well-posedness cause these problems 
to possess a very demanding structure. Most notably, one has to deal with the 
following difficulties when studying problems like \eqref{eq:model_problem_elliptic}:
\begin{enumerate}[label=\roman*)]
\item\label{item:problemn:i} 
As the optimisation variable $q$ enters nonlinearly in the principal part 
of the governing PDE in problems like \eqref{eq:model_problem_elliptic}, one requires the ansatz space for $q$
to possess suitable compactness properties 
to be able to ensure the existence of minimisers 
via the direct method of the calculus of variations.
(Note in particular that an $L^2$-ansatz for $q$ is insufficient here
and that the presence of the nonlinear Nemytskii operator $B$ also makes it impossible
to work with the concept of $H$-convergence in \eqref{eq:model_problem_elliptic}; cf.\ \cite{Deck:2011,Deck:2012}.) Further, the design variable 
$q$ is typically subject to distributed pointwise 
inequality constraints in inverse homogenisation 
problems---in \eqref{eq:model_problem_elliptic}, for example, 
the constraint $ a \leq q \leq b$ a.e.\ in $(0,L)$ ensures that the circular hole of the perforation is entirely 
contained in the periodicity cell $Y$. In total, 
this results in an optimisation problem 
involving a high-regularity optimisation variable $q$
with constraints---a problem type that is known to be challenging. 

\item\label{item:problemn:ii} 
Owing to \ref{item:problemn:i},
any form of projected-gradient descent method for the 
numerical solution of a problem like \eqref{eq:model_problem_elliptic}
has to perform projections in a high-regularity space---typically
an $H^s$-space for some $s > 0$---onto the admissible set 
to be consistent and mesh-independent. 
These projections entail high numerical costs,
in particular given the fact that gradient methods 
typically require dozens of iterations to achieve a reasonable 
accuracy; cf.\ \cite{Blank2017}. 

\item\label{item:problemn:iii} 
In a problem like \eqref{eq:model_problem_elliptic},
the map $B\colon [a,b] \to \R^{2\times 2}$ acts by superposition 
on the distributed design variable $q$. As the matrix 
$B(r)$ is obtained by 
solving the associated cell problems for the perforation radius $r$, 
this means that 
one requires the solutions 
of the cell-problems not only for single cases 
but for a continuous 
range of microstructure parameters.
As solving infinitely many cell problems 
is typically impossible, 
this makes it necessary to work with surrogate models 
that emulate the exact relationship between the 
microstructure parameters and the macroscopic behaviour that 
they induce; see, e.g.\ \cite{All:2019,Bendsoe1988}. 
\end{enumerate}

In the present paper, we are primarily concerned 
with points \ref{item:problemn:i} and \ref{item:problemn:ii} 
above. Our main contributions are as follows: 
First, we develop a 
general functional analytic framework 
for inverse homogenisation problems 
involving macroscopically unidirectionally 
varying microstructure features; see \cref{sec:abstractsetting}.
In this setting, we discuss the 
existence of solutions, first-order necessary optimality conditions, 
and associated differentiability and adjoint calculus 
results; see \cref{subsec:well-posedness_etc}.
Second, we analyse 
the convergence of an infinite-dimensional gradient projection algorithm 
for the numerical solution of the considered class of 
inverse homogenisation problems in \cref{subsec:gradproj}. 
Although gradient projection 
algorithms constitute the state of the art in the field 
of inverse homogenisation---see, e.g., \cite{Loch:2022,Loch:2024,Lipton2006,All:2019}---their convergence 
on the infinite-dimensional level 
seems to have not been discussed so far in this area. 
While the convergence analysis of our gradient projection method
is quite standard from the viewpoint of optimisation theory, 
it also exploits the particular structure 
of the considered inverse 
problems and yields computable bounds for the
admissible step sizes in terms of the Lipschitz constants 
of the involved operators, thus providing additional 
value in comparison with general purpose results on the topic;
see \cref{th:convgradproj,cor:global_convergence_tracking_type}.
Lastly (and most importantly), we develop and
discuss a novel semismooth Newton approach for the numerical 
solution of problems like \eqref{eq:model_problem_elliptic}
in \cref{subsec:semiNewton}; see 
\cref{th:convergence_Newton,cor:invert_Newton_large_Tikhonov}. 
This method has the 
main advantage that it allows to reduce the number of the 
costly projections mentioned in \ref{item:problemn:ii} drastically
in comparison with gradient-based approaches
and, thus, offers a significant speed-up
over existing 
algorithms.
It is based on recent results on the Newton 
differentiability of projections in Sobolev spaces 
onto sets with pointwise bounds, see \cite{ChristofWachsmuthUniObstacle,ChristofWachsmuthBiObstacle,Rauls2020}---an operation that arises very naturally 
when inverse homogenisation problems are considered 
due to the structure described in \ref{item:problemn:i} above. 
We highlight that, 
for the convenience of the reader, we revisit and extend 
the results of \cite{ChristofWachsmuthUniObstacle,ChristofWachsmuthBiObstacle,Rauls2020}
in \cref{sec:appendix} of this paper 
in a form that is tailored 
to the 
application to problems like \eqref{eq:model_problem_elliptic}
involving unidirectionally varying microstructures. 
We expect that the results proven there in a 
(largely) self-contained form along with the 
analysis in \cref{subsec:semiNewton} are also of interest 
for other application areas, e.g., phase-field models;
see \cite{Blank2017}. 
Finally, we would like to mention at this point that 
our focus on microstructures that vary 
only in one spatial direction is not only due to the 
practical relevance of this type of macroscopic
 dependence but also
 due to the fact that our functional analytic 
 setting benefits greatly from the one-dimensionality
 of the domain
 of the design variable and the additional 
 regularity that it provides (in the context of \eqref{eq:model_problem_elliptic}, for example, the embedding 
 $H^1(0,L) \hookrightarrow L^\infty(0,L)$). We expect that, for 
 particular instances of inverse homogenisation problems, 
 e.g.~linear elasticity, the assumption of unidirectionality 
 of the macroscopic variations can be removed by employing Gröger-regularity; see \cite{Herzog2011}. 
 We leave this topic for future research. 

 We conclude this introduction with a short outline 
 of the rest of the paper: 
 \cref{subsec:notation} below introduces basic notation
 and conventions that are needed for our analysis. 
 In \cref{sec:homogen}, we 
 briefly recall basics from periodic 
 homogenisation theory, in particular the role
 of the cell problems and their relationship to the 
 effective tensor in a homogenised equation. 
\Cref{sec:abstractsetting} is then concerned with 
our abstract functional analytic setting. 
Here we analyse the well-posedness 
of the considered class of inverse 
homogenisation problems (\cref{subsec:well-posedness_etc}),
analyse our gradient projection method (\cref{subsec:gradproj}),
and develop our semismooth Newton approach (\cref{subsec:semiNewton}). 
\cref{sec:num_exp} contains numerical experiments which demonstrate 
that our semismooth Newton method is competitive and 
indeed requires far less 
projection operations than a gradient-descent method.
\cref{sec:appendix} finally establishes the Newton differentiability 
results needed for \cref{subsec:semiNewton} by revisiting and 
generalising the analysis of \cite{ChristofWachsmuthUniObstacle,ChristofWachsmuthBiObstacle,Rauls2020};
see \cref{th_main_appendix} 
for the main result on this topic.

\subsection{Notation and basic concepts}
\label{subsec:notation}

In this paper, we use the standard 
symbols 
$\N$, $\N_0$, $\Z$, $\R^+$, and $\R$
to denote the sets of natural numbers, 
nonnegative integers, 
integers, positive real numbers, 
and real numbers, respectively.
For the Euclidean norm on $\R^d$
and the Frobenius norm on $\R^{d \times d}$
for some $d \in \N$, we use the notation $|\cdot|$,
and for the Euclidean scalar product on $\R^d$, the 
dot ``$\cdot$''. The set of symmetric matrices 
$\{B \in \R^{d\times d} \mid B = B^\top\}$ is denoted by 
$\R^{d \times d}_{\sym}$. 
For norms and inner products 
on (real) vector spaces, we use the symbols
$\|\cdot\|_X$ and $(\cdot, \cdot)_X$, respectively, 
with the subscript $X$ referring to the space under consideration. 
Given a subset $D$ of a normed space $(X,\|\cdot\|_X)$,
we write $\closure(D)$, $\interior(D)$, and $\partial D$
for the closure, interior, and boundary of $D$ in $X$,
respectively. A subscript is added in this notation 
to highlight the underlying topology if necessary, i.e.~$\closure_X(D)  := \closure(D)$ etc.
The indicator function 
of a set $D$ is denoted by 
$\mathds{1}_D\colon X \to \{0,1\}$. 
For a nonempty closed and convex set $D$,
we define $\TT(x; D) := \closure_X(\R^+(D - x))$
to be 
the tangent cone of $D$ at $x \in D$.
If $(X_i, \|\cdot\|_{X_i})$, $i=1,2,3$, are normed spaces,
then 
 $\LLLL(X_1, X_2)$
 denotes the space of 
 linear and continuous functions from 
 $X_1$ to $X_2$ and 
 $\BBBB(X_1, X_2;X_3)$
 the space of 
 continuous bilinear functions 
 from $X_1\times X_2$ to $X_3$.
For the special case of the topological dual space
of a normed space $(X, \|\cdot\|_X)$, 
we write $X^* := \LLLL(X, \R)$
and set $\langle x^*, x\rangle_X := x^*(x)$
for all $x^* \in X^*$, $x \in X$.
The identity map and the identity matrix
are denoted by $\Id$.
Given $T \in \LLLL(X_1, X_2)$,
we write $T^* \in \LLLL(X_2^*, X_1^*)$ for the 
adjoint operator of $T$ and 
$\ker(T)$ for the kernel of $T$. 
Strong, weak, and weak-star convergences 
are indicated with the arrows
$\to$, $\weakly$, and $\weaklystar$, respectively,
and for convergence with respect to (w.r.t.) the weak 
operator topology, we use the symbol $\WOT$;
see \cite[Section 1.1]{ChristofWachsmuthUniObstacle}.
If $(X_1, \|\cdot\|_{X_1})$ is continuously
embedded into $(X_2, \|\cdot\|_{X_2})$,
then we write $X_1 \hookrightarrow X_2$, and if
the embedding is compact, 
then we highlight this with the notation 
$X_1 \compactly X_2$. 
For a Hilbert space $X$, the Riesz 
isomorphism is denoted by $\iota_R\colon X^* \to X$,
the orthogonal complement of a closed
subspace $W \subset X$ is denoted by $W^\perp$,
and the projection in $X$ onto a nonempty 
convex and closed set $D \subset X$
is denoted by $\pi_D\colon X \to D$. 
If $F\colon X_1 \to X_2$ is a function between
normed spaces $(X_i, \|\cdot\|_{X_i})$, $i=1,2$,
and $x, h \in X_1$, then 
we write $F'(x;h)$
for the directional derivative 
of $F$ at $x$ in direction $h$
and $F'(x)$ for the Gâteaux/Fréchet derivative 
of $F$ at $x$, provided these derivatives exist;
see \cite[Section 2.2]{BonnansShapiro2000} for details.
Recall that a function $F\colon X_1 \to X_2$
that is Gâteaux differentiable with a continuous 
derivative $X_1 \ni x \mapsto F'(x) \in \LLLL(X_1, X_2)$
is already Fréchet differentiable by the mean value theorem;
see \cite[Proposition 3.2.15]{Drabek2007}.
Such a function $F$ is called $C^1$. 
For second Fréchet derivatives of a function $F\colon X_1 \to X_2$, 
we write
$F''(x) \in \BBBB(X_1, X_1; X_2)$,
and if $F$ has several arguments $x_1, x_2, \ldots$, 
then we use the usual symbols  
$\partial_1$, $\partial_{x_1}$,
$\partial_{x_1 x_2}$ etc.\ to denote partial 
first- and second-order Fréchet derivatives of $F$.
In the special case that $X_1$ is Hilbert and 
$X_2 = \R$, we moreover define 
$\nabla_{X_1} F(x) := \iota_R(F'(x)) = (\iota_R \circ F')(x)$ 
to be the gradient of $F$ at $x$ w.r.t.\ the 
$X_1$-inner product. Here, the circle $\circ$ denotes a 
composition, as usual. Recall that 
a function $F\colon X_1 \to X_2$
endowed with a 
generalised, potentially set-valued 
differential $\D F\colon X_1 \rightrightarrows \LLLL(X_1, X_2)$
is called \emph{Newton differentiable} if, 
for all $x \in X_1$, we have 
\begin{equation}
\label{eq:NewtonDiffDef}
\sup_{G \in \D F(x + h)} 
\frac{\|F(x + h) - F(x) - G h\|_{X_2}}{\|h\|_{X_1}} \to 0
\quad\text{ for } 0 <  \|h\|_{X_1} \to 0.
\end{equation}
In this case, $\D F$ is also called a \emph{Newton derivative} of $F$.
If $\D F(x)$ is a singleton for all $x \in X_1$, 
then we simply 
write $\D F\colon X_1 \to \LLLL(X_1, X_2) $. 
For the level sets of a function 
$F\colon X_1 \to X_2$, we use the 
shorthand notation 
$\{F = c\} := \{x \in X_1 \mid F(x) = c\}$.
Analogous conventions are used
for other preimages, e.g., sub- and superlevel sets 
in the case $X_2 = \R$.

Given a Lebesgue-measurable set 
$\Omega \subset \R^d$, $d \in \N$,
we denote the vector space of 
(equivalence classes of) real-valued Lebesgue measurable 
functions on $\Omega$ by $L^0(\Omega)$
and the classical 
Lebesgue spaces by $L^p(\Omega)$, $p \in [1,\infty]$.
If $(X, \|\cdot\|_X)$ is a separable Banach space, 
then we further write 
$L^p(\Omega;X)$, $p \in [1,\infty]$,
for the classical Lebesgue-Bochner spaces 
on $\Omega$ with values in $X$,
and if $\Omega \subset \R^d$ is open, 
then we use the standard symbols 
$H_0^k(\Omega)$, 
$H^s(\Omega)$,
$W^{k,p}(\Omega)$
with $k \in \N$, $p \in [1,\infty]$, $s \in [0, \infty)$
for the (fractional and ordinary) Sobolev spaces on $\Omega$;
see, e.g., 
\cite{Adams2003,Necas2012}
for details. 
For weak/distributional derivatives
of elements of the Sobolev spaces
and related differential operators, 
we use standard notation, for instance, $\nabla$ for the 
gradient and $\nabla \cdot$ for the divergence. 
A subscript is added to these symbols if derivatives are taken w.r.t.\ a certain variable.
If $\Omega$ is a bounded Lipschitz domain, 
then evaluations of elements 
of the Sobolev spaces on $\partial \Omega$
are always understood in the sense of traces.
The Lebesgue spaces associated with the 
surface measure of $\partial \Omega$ are then 
denoted as usual by 
$L^p(\partial \Omega)$, $L^p(\partial \Omega; \R^d)$
etc.\ for $p \in [1,\infty]$.
The dual space of $H_0^1(\Omega)$ (with pivot space $L^2(\Omega)$)
is denoted by $H^{-1}(\Omega)$.
Given an open set $\Omega \subset \R^d$, $d \in \N$,
with closure $\overline{\Omega} := \closure(\Omega)$,
we further use the standard symbols
$C(\Omega)$, $C(\overline{\Omega})$, $C^{0,1}(\Omega)$,
$C^k(\Omega)$, and $C^k(\overline{\Omega})$
with $k \in \N \cup \{\infty\}$ 
for the vector spaces 
of continuous, Lipschitz continuous, 
$k$-times continuously differentiable, and smooth functions 
on $\Omega$ and $\overline{\Omega}$, respectively 
(with the definition 
$C^k(\overline{\Omega}) := \{\phi \mid \phi = \tilde \phi \text{ on }
\overline{\Omega} \text{ for some } \tilde \phi\in C^k(\R^d)\}$). 
If $\Omega$ is additionally bounded, then we
also define $\MM(\overline{\Omega})$ to be the space 
of signed finite Radon measures on $\overline{\Omega}$,
so that $C(\overline{\Omega})^* \cong \MM(\overline{\Omega})$,
and denote by $\supp(\cdot)$ the measure-theoretic support.
Note that, if $\Omega$ is an interval, 
then 
we typically drop the brackets of the interval in 
our notation for function spaces, i.e.\
we write $H^1(0, 1)$ 
instead of $H^1((0,1))$ etc.

\section{A short recap of periodic homogenisation theory}
\label{sec:homogen}

Before 
we introduce the functional analytic setting that 
we consider in the remainder of this paper, 
we briefly recall basic results 
on the periodic homogenisation of PDEs,
upscaled effective equations, and the field of inverse homogenisation;
see \cite{Cior:2018,Cior:1999}. 
For the sake of simplicity, 
we focus on a classical second-order linear elliptic PDE with homogeneous 
Dirichlet boundary conditions and a 
two-scale coefficient matrix
in what follows, that is 
\begin{equation}\label{eq:model_PDE}
-\nabla \cdot (\BBB_\varepsilon \nabla u_\varepsilon)=f 
\text{ in }\Omega,\qquad 
	u_\varepsilon =0  \text{ on }\partial\Omega.
\end{equation}
Our standing assumptions are as follows:
\begin{assumption}[{Standing assumptions for \cref{sec:homogen}}]
\label{ass:standing:homo}~
\begin{enumerate}[label=\roman*)]
\item $\Omega \subset \R^d$, $d \in \N$, is a bounded Lipschitz domain
with boundary $\partial \Omega$ (the macroscopic body);
\item $f \in H^{-1}(\Omega)$ is given and fixed; 
\item $\{\BBB_\varepsilon\}_{\varepsilon > 0}$
is a bounded family in $L^\infty(\Omega;\R^{{d\times d}})$
and there exists $\alpha > 0$ such that 
\begin{equation}
\label{eq:ellipt_ass_epsilon}
    \alpha |\zeta|^2 \leq \BBB_\varepsilon(x)\zeta \cdot \zeta\qquad 
    \forall \zeta \in \R^d\qquad \text{f.a.a.\ } x\in \Omega
    \qquad \forall \varepsilon > 0.
\end{equation}
\end{enumerate}
\end{assumption}

Note that, under the conditions above, 
the boundary-value problem \eqref{eq:model_PDE} 
has a unique weak solution $u_\varepsilon \in H_0^1(\Omega)$
for all $\varepsilon > 0$
by the lemma of Lax--Milgram; see \cref{lem:well-posedness_state} in \cref{sec:abstractsetting}. 

The situation that one should have in mind in the 
context of periodic homogenisation is that the 
coefficient
matrix $\BBB_\varepsilon$ in 
\eqref{eq:model_PDE} varies along 
two length scales, one much smaller than the other, 
with $\varepsilon$-periodicity on the smaller scale.
That is, one usually has 
$\BBB_\varepsilon (x) := \BBB(x, x/\varepsilon)$
with a function $\BBB\colon \Omega \times \R^d \to \R^{d\times d}$
that is periodic in the second argument 
w.r.t.\ a \emph{periodicity cell} $Y \subset \R^d$,
typically $Y:=(0, 1)^d$.
As a prototypical example,
one can consider here the 
function  
$\BBB\colon \Omega \times \R^d \to \R^{d\times d}$
that is obtained by extending the definition 
\begin{equation}
\label{eq:coeff_perforation_example}
\BBB(x,y) := 
\big [
c_1 + (c_2 - c_1)
\mathds{1}_{\{\zeta \in Y \colon |\zeta - m_Y| \leq q(x_1)\}}(y)
\big ] \Id\qquad \forall (x,y) \in \Omega \times Y
\end{equation}
$Y$-periodically in $y$ to $\Omega \times \R^d$, 
with $\Omega := (0,L) \times D$ for some bounded 
Lipschitz domain $D \subset \R^{d-1}$ and $L>0$, 
$Y:=(0, 1)^d$, 
two constants $0 < c_2 < c_1$,
$m_Y := (\frac12,\ldots, \frac12)^\top$ 
as the midpoint of $Y$, and a function 
$q \in H^1(0,L) \subset C([0,L])$ satisfying
$0 < q < \frac12$ in $[0,L]$. 
For this type of matrix field,
$\BBB_\varepsilon (x) := \BBB(x, x/\varepsilon)$
describes a structure 
with periodic spherical inclusions of length-scale $\varepsilon > 0$
whose radius $\varepsilon q$ varies macroscopically 
with $x_1$---a scenario analogous to that outlined after 
\eqref{eq:model_problem_elliptic}. 

While the boundary-value 
problem \eqref{eq:model_PDE} is well-posed for all $\varepsilon > 0$
owing to the ellipticity condition \eqref{eq:ellipt_ass_epsilon}, 
it is infeasible to compute its solution $u_\varepsilon$ 
numerically if the scale, on which $\BBB_\varepsilon$ varies spatially, is 
equal to $\varepsilon$ and $\varepsilon$ is very small. It is then instrumental to pass to
the limit $\varepsilon \to 0$ in \eqref{eq:model_PDE} and
to work with an effective limit problem whose solution
can be used as an
approximation of $u_\varepsilon$ for $0 < \varepsilon \ll 1$. 
A suitable tool for the derivation of 
an effective limit equation is 
the so-called \emph{periodic unfolding operator}; see, for instance, 
\cite[Chapter 1.1]{Cior:2018}.

\begin{definition}[Periodic unfolding]\label{def:periodic_unfold_op}
Let 
$\varepsilon > 0$  be fixed.
Consider the cell
$Y := (0,1)^d$ and set 
\[
\Lambda_\varepsilon:=
\left \{\xi\in\mathbb{Z}^d: \varepsilon(Y+\xi)\subset\Omega\right \},
\qquad
\Omega_\varepsilon:=\interior \left (\bigcup_{\xi\in\Lambda_\varepsilon}\varepsilon(\closure(Y)+\xi) \right ).
\]
For each $\zeta \in \R^d$,
define
$\left \lfloor\zeta\right \rfloor$ to be the 
unique element of $\Z^d$ that satisfies 
$\zeta - \left \lfloor\zeta\right \rfloor\in [0,1)^d$. Then, 
the periodic unfolding operator 
$\mathcal{T}_\varepsilon$ is the function
$\mathcal{T}_\varepsilon:L^0(\Omega_\varepsilon)\to L^0(\Omega\times Y)$ given by
	\begin{equation*}
		\mathcal{T}_\varepsilon(\phi)(x,y):=\begin{cases}
        \displaystyle
		\phi \left (\varepsilon \left \lfloor \frac{x}{\varepsilon}\right \rfloor+\varepsilon y \right) &\text{ for a.e. }(x,y)\in\Omega_\varepsilon\times Y,\\ \displaystyle
		0 &\text{ for a.e. }(x,y)\in (\Omega\setminus\Omega_\varepsilon)\times Y.
		\end{cases}
	\end{equation*}	
    \end{definition}

We remark that 
$\mathcal{T}_\varepsilon$ maps $L^s(\Omega_\varepsilon)$
linearly and continuously 
into $L^s(\Omega\times Y)$ for all $s \in [1,\infty]$; 
see \cite[Proposition 1.8]{Cior:2018}. With $\mathcal{T}_\varepsilon$ at hand,
one can prove: 

\begin{theorem}[{\cite[Theorem 3.5]{Cior:2018}}]\label{theo:homog_prob_not_scale_sep_gen}
Let $Y$ and $\mathcal{T}_\varepsilon$ be as 
in \cref{def:periodic_unfold_op}. 
Let $u_\varepsilon$ be the unique weak solution to \eqref{eq:model_PDE}
and suppose that there exists $\BBB \in L^\infty(\Omega \times Y)$
such that 
\begin{equation}
\label{eq:unfolding_limit}
\mathcal{T}_\varepsilon(\BBB_\varepsilon)\to \BBB\text{ pointwise a.e.~in }\Omega\times Y
\end{equation}
holds for $\varepsilon \to 0$, where $\mathcal{T}_\varepsilon$ acts componentwisely. 
Define 
\[
H^1_{\mathrm{per},0}(Y)
:=
\left \{
\phi \in H^1(Y)
\,\colon
\phi \text{ is periodic and } 
\int_Y \phi\,\mathrm{d}y = 0
\right \}.
\]
Then, $\BBB$
satisfies $\alpha |\zeta|^2 \leq \BBB(x, y)\zeta \cdot \zeta$
for all $\zeta \in \R^d$ and a.a.\ $(x,y) \in \Omega \times Y$ 
and there exist $u\in H^1_{0}(\Omega)$ and 
$\hat{u}\in L^2(\Omega;H^1_{\mathrm{per},0}(Y))$ such that
\begin{align*}
u_\varepsilon\weakly u \text{ in }H_0^1(\Omega),
\quad
\mathcal{T}_\varepsilon(u_\varepsilon)\to u \text{ in }L^2(\Omega;H^1(Y)),
~~\text{and}~~
\mathcal{T}_{\varepsilon}(\nabla u_\varepsilon)\rightharpoonup\nabla u+\nabla_y\hat{u} \text{ in }L^2(\Omega\times Y)^{d}
\end{align*}
for $\varepsilon \to 0$.
Further, the tuple  
$(u,\hat{u})$ is the unique solution of the variational problem
\begin{equation}
\label{eq:cell_problem_prototypical}
\begin{aligned}
&\text{Find }(u,\hat{u}) \in  H^1_{0}(\Omega) \times L^2(\Omega;H^1_{\mathrm{per},0}(Y))
\text{ such that}
\\
&\int_{\Omega\times Y}\BBB(x,y)(\nabla_x u(x)+\nabla_y\hat{u}(x,y))(\nabla_x \phi(x)+\nabla_y\hat{\phi}(x,y))\,\mathrm{d}x\mathrm{d}y=
\left \langle
f, \phi
\right \rangle_{H_0^1(\Omega)}
\\
&\hspace{6.3cm}
\forall (\phi, \hat \phi) \in H^1_{0}(\Omega) \times L^2(\Omega;H^1_{\mathrm{per},0}(Y)).
\end{aligned}
\end{equation}
\end{theorem}

Note that, for a coefficient matrix field 
$\BBB_\varepsilon$
of the form 
$\BBB_\varepsilon (x) := \BBB(x, x/\varepsilon)$
involving a function $\BBB\colon \Omega \times \R^d \to \R^{d\times d}$
that is $Y$-periodic in its second argument with  $Y=(0, 1)^d$,
it follows from \cref{def:periodic_unfold_op} that 
	\begin{equation}
    \label{eq:randome28bhe8n}
		\mathcal{T}_\varepsilon(\BBB_\varepsilon)(x,y):=
        \begin{cases}
        \displaystyle
        \BBB\left (
        x + \varepsilon\left ( \left \lfloor\frac{x}{\varepsilon}\right\rfloor - \frac{x}{\varepsilon} + y \right ),
          y \right)
        &\text{ for a.e. }(x,y)\in\Omega_\varepsilon\times Y,\\ \displaystyle
		0 &\text{ for a.e. }(x,y)\in (\Omega\setminus\Omega_\varepsilon)\times Y.
		\end{cases}
	\end{equation}	
For such a field, one thus has
$\mathcal{T}_\varepsilon(\BBB_\varepsilon) \to \BBB$ a.e.\ in 
$\Omega \times Y$ for $\varepsilon \to 0$ if 
$\BBB$ is sufficiently 
smooth in its first argument for enough $y \in Y$. 
A condition that ensures the pointwise-a.e.\  
convergence in \eqref{eq:unfolding_limit} 
(and, as a consequence, the applicability of 
\cref{theo:homog_prob_not_scale_sep_gen}) in the situation of 
\eqref{eq:randome28bhe8n} is,
for example, that $\BBB\colon \Omega \times \R^d \to \R^{d\times d}$ 
is continuous on 
$(\Omega \times \R^d) \setminus \mathcal{N}$ for
some set $\NN$ of measure zero. Indeed, 
under this assumption, 
the right-hand side of \eqref{eq:randome28bhe8n} 
tends to $\BBB(x,y)$ for all 
$(x,y) \in (\Omega \times Y) \setminus \mathcal{N}$ as $\varepsilon \to 0$,
and the convergence 
$\mathcal{T}_\varepsilon(\BBB_\varepsilon) \to \BBB$ a.e.\ in 
$\Omega \times Y$ for $\varepsilon \to 0$ is readily obtained.
This criterion is satisfied, for instance, 
for the coefficient field in \eqref{eq:coeff_perforation_example}
with the set 
$
\NN := \left \{ (x,y) \in \Omega \times Y ~ 
\colon  |y - m_Y| = q(x_1)\right \},
$
that, as a surface arising from the 
extrusion and revolution of the graph of the 
function $q$ (which has finite length due to the 
assumed regularity $q \in H^1(0,L)$) has $2d$-dimensional volume zero. 

By employing standard arguments (see \cite[Proposition 3.7]{Cior:2018}) 
the problem 
\eqref{eq:cell_problem_prototypical} 
in \cref{theo:homog_prob_not_scale_sep_gen}
can also be written in a more convenient form: 
The weak limit  $u \in H_0^1(\Omega)$
of the functions $u_\varepsilon$
is the unique weak solution to the \emph{effective equation}
\begin{equation}
\label{eq:effective_PDE}
-\nabla \cdot (\hat \BBB \nabla u ) =f 
\text{ in }\Omega,\qquad 
	u  =0  \text{ on }\partial\Omega,
\end{equation}
involving 
the \emph{effective matrix field} $\hat \BBB \in L^\infty(\Omega; \R^{d\times d})$
defined by 
$\hat \BBB(x)=(\hat \bbb_{ij}(x))_{1\leq  i,j \leq d}$ with  
\begin{equation}\label{def:A_hom_x}
\hat \bbb_{ij}(x)
= \int_Y \bbb_{ij}(x,y) + \BBB(x,y) \nabla_{y}\chi_{j}(x,y) \cdot e_i \,\mathrm{d}y
\quad \text{f.a.a.\ } x \in \Omega
\qquad\forall 1 \leq i,j \leq d.
\end{equation}
Here, $\bbb_{ij}$ are the components of $\BBB$,
$e_i$, $i=1,\ldots, d$, are the unit vectors, 
and $\chi_j$, $j=1,\ldots, d$, are the unique solutions to the \emph{cell problems}
\begin{equation}\label{eq:homog_weak_cell_problem_gen}
\begin{aligned}
\text{Find }
\chi_j \in L^\infty(\Omega;H^1_{\mathrm{per},0}(Y))
\text{ such that } 
\int_Y \BBB(x,y)\left(\nabla_y \chi_j (x,y)+e_{j}\right)\cdot&\nabla_y \phi(y)\,\mathrm{d}y=0
\\
&\forall \phi\in H^1_{\mathrm{per},0}(Y).
\end{aligned}
\end{equation}
Note that, if $\BBB$ is $Y$-periodic in its second argument, then 
\eqref{eq:homog_weak_cell_problem_gen} can also be written 
strongly as
\begin{equation*}
 -\nabla_y\cdot \big ( \BBB(x,\cdot)(\nabla_y \chi_j (x,\cdot)+e_{j}) \big )=0 \text{ in }Y,
    \quad
 \chi_j(x,\cdot)\,\text{$Y$-periodic with } \int_Y \chi_j(\cdot,y) \dd y=0 \text{ in } \Omega.
 \end{equation*}
Let us now again 
consider the example in \eqref{eq:coeff_perforation_example}.
For this type of coefficient field, we obtain from the formulas above that 
$\hat \BBB\in L^\infty(\Omega; \R^{d\times d})$ is given by 
$\hat \BBB(x) = B(q(x_1))$ for a.a.\ $x \in \Omega = (0,L) \times D$, where 
$B\colon (0, 1/2) \to \R^{d\times d}$ is the function 
that assigns to each number $r \in (0,1/2)$ the matrix 
$B(r) := (b_{ij}(r))_{1 \leq i,j \leq d}$ 
characterised by the ($x$-independent) system 
\begin{equation}
\label{eq:r_system}
\begin{gathered}
b_{ij}(r)
= \int_Y \bbb_{r,ij}(y) + \BBB_r(y) \nabla_{y}\hat \chi_{j}(y) \cdot e_i \,\mathrm{d}y
\qquad\forall 1 \leq i,j \leq d,
\\
\hat \chi_j \in  H^1_{\mathrm{per},0}(Y),
\quad
\int_Y \BBB_r(y)\left(\nabla_y \hat \chi_j (y)+e_{j}\right)\cdot\nabla_y \phi(y)\,\mathrm{d}y=0
\quad\forall \phi\in H^1_{\mathrm{per},0}(Y)
\end{gathered}
\end{equation}
with 
\begin{equation}
\label{eq:coeff_again_42}
\BBB_r(y) 
:= (\bbb_{r,ij}(y))_{1 \leq i,j \leq d}
:=
\big [
c_1 + (c_2 - c_1)
\mathds{1}_{\{\zeta \in Y \colon |\zeta - m_Y| \leq r\}}(y)
\big ] \Id\qquad \forall y \in Y. 
\end{equation}
For the case of spherical microinclusions
with a macroscopically unidirectionally varying 
radius $0 < q(x_1) < 1/2$, $q \in H^1(0,L)$, 
it thus follows that the homogenised
coefficient matrix in the effective equation 
\eqref{eq:effective_PDE} is obtained by applying a function
$B\colon (0, 1/2) \to \R^{d\times d}$ in a superposition sense 
to $q$. In particular, 
if one wants to identify a radius distribution $q \in H^1(0,L)$ that 
gives rise to a certain macroscopic state $u_\mathrm{obs}$,
then 
one ends up with a problem that has precisely 
the structure \eqref{eq:model_problem_elliptic},
with a Nemytskii operator entering in the principal part 
of the governing effective equation
and with all the challenges \ref{item:problemn:i}-\ref{item:problemn:iii}
described in \cref{sec:intro}. Note that the derivation above 
shows that the difficulty in point \ref{item:problemn:iii}
in \cref{sec:intro} is indeed quite severe as 
the system
\eqref{eq:r_system} is rather involved 
and can only very rarely be solved analytically. 
(See 
 \cite[Sections 5.3, 5.4]{Cior:1999}
for some exceptions.)
Moreover, 
we see in \cref{theo:homog_prob_not_scale_sep_gen} that 
higher regularity properties for $q$ are not only needed 
to ensure the well-posedness of problems like 
\eqref{eq:model_problem_elliptic} but also to 
justify rigorously the homogenisation limit; 
see the discussion after \eqref{eq:randome28bhe8n}.

We remark that, for other elliptic PDEs and 
microstructures different from that in \eqref{eq:coeff_perforation_example},
one observes the same effects as above. 
See, e.g.,  \cite[Section 10]{Cior:1999} and 
\cite{Loch:2020} for analogue results 
in linear elasticity. The structure 
in \eqref{eq:model_problem_elliptic} is thus indeed 
generic for inverse homogenisation problems. 

\section{Abstract problem formulation}
\label{sec:abstractsetting}
Motivated by the considerations in the last section and \eqref{eq:model_problem_elliptic}, 
we henceforth focus on minimisation problems of the following type as general models for inverse homogenisation problems: 
\begin{equation*}
\label{eq:P}
\tag{P}
\left.
	\begin{aligned}
		\text{Minimise} 
		\quad &  J(v, q)  \\
        \text{w.r.t.}
        \quad &v \in V, \quad  q \in Q,\\
		\text{s.t.}
		\quad & 
        \left (A(q)v, w \right )_{L^2(\Omega)^m} = F(w)\quad \forall w \in V
  \\
		\text{and} \quad & a \leq q \leq b \text{ a.e.\ in } (0,L).
\end{aligned}
\quad \right \}
\end{equation*}
Our standing assumptions on the quantities in 
\eqref{eq:P} are as follows:
\begin{assumption}[{Standing assumptions for the study of problem \eqref{eq:P}}]
\label{ass:standing:abstract}~
\begin{enumerate}[label=\roman*)]
\item\label{ass:standing:abstract:i}
$\Omega := (0,L) \times D$
for a nonempty Lebesgue measurable set $D \subset \R^n$, $n \in \N$,
and some $L>0$;
\item $V$ is a closed subspace of $L^2(\Omega)^m$ for some $m \in \N$
endowed with the norm $\|\cdot\|_V := \|\cdot\|_{L^2(\Omega)^m}$;
\item\label{ass:standing:abstract:iii}
$Q$ is a reflexive Banach space 
satisfying  $Q \compactly L^\infty(0,L) \subset  L^\infty((0,L) \times D) = L^\infty(\Omega)$;
\item\label{ass:standing:abstract:iv} $a,b \in L^0(0,L)$ are given functions
such that $\Qad := \{q \in Q \mid a \leq q \leq b \text{ a.e.\ in } (0,L) \} \neq \emptyset$;
\item\label{ass:standing:abstract:v} $A\colon \R \to \R^{m \times m}$ is 
continuously differentiable with 
derivative $A'(s) \in \R^{m \times m}$ for all $s \in \R$,
acting as a Nemytskii operator
$A\colon Q \to L^\infty(\Omega)^{m \times m}$,
and there exists $\alpha  > 0$ such that 
\begin{equation}
\label{eq:elliptic_abs}
    \alpha \|v\|_{L^2(\Omega)^m}^2 \leq 
\left ( 
A(q)v, v
\right )_{L^2(\Omega)^m}
\qquad \forall v \in V\qquad 
\forall q \in \Qad;
\end{equation}
\item $F \in V^*$ is given and fixed;
\item\label{ass:standing:abstract:vii}
$J\colon V \times Q \to \R$ is an objective function such that:
\begin{enumerate}
\item $J$ is Fréchet differentiable with partial derivatives 
$\partial_v J \colon V\times Q \to V^*$, 
$\partial_q J \colon V\times Q \to Q^*$;
\item $J$ is lower semicontinuous 
in the sense that, for 
all $\{v_k\} \subset V$ and $\{q_k \} \subset \Qad$
satisfying $v_k \to v$ in $V$ 
and $q_k \weakly q$ in $Q$,
we have 
\[
\liminf_{k \to \infty} J(v_k, q_k) \geq J(v,q);
\]
\item $J$ is radially unbounded in the sense that 
there exists a function $\rho\colon [0, \infty) \to \R$
satisfying  
\[
\qquad\rho(s) \to \infty \text{ for } s \to \infty
\qquad
\text{and}
\qquad
J(v,q) \geq \rho\left ( \|q\|_Q \right ) \quad \forall (v,q) \in V \times \Qad.
\]
\end{enumerate}
\end{enumerate}
\end{assumption}

As many Sobolev spaces arising in the study of 
elliptic PDEs are isomorphic to closed subspaces 
of $L^2(\Omega)^m$ for suitable $m$, 
the abstract problem setting \eqref{eq:P}
covers numerous interesting situations.
We demonstrate this by means of two examples:

\begin{example}[Scalar elliptic PDEs]
\label{ex:scalar}
Assume that the set $\Omega \subset \R^d$, $d := n + 1$,
in \cref{ass:standing:abstract}\ref{ass:standing:abstract:i}
is a bounded Lipschitz domain,
that the boundary $\partial \Omega$ 
is split disjointly into a (possibly empty) 
relatively open Neumann part
$\Gamma_N \subset \partial \Omega$ and a
closed
Dirichlet part $\Gamma_D \subset \partial \Omega$ of positive surface measure,
and that the bounds 
$a$ and $b$ in 
\cref{ass:standing:abstract}\ref{ass:standing:abstract:iv}
are real numbers satisfying $ a < b$. 
Denote the outer unit normal vector field of $\Omega$ by
$\vartheta \in L^\infty(\partial \Omega; \R^d)$ and 
suppose that we are given a PDE of the type
\[
-\nabla\cdot (B(q)\nabla u)=f \text{ in } \Omega,
\qquad
u=0 \text{ on } \Gamma_D,
\qquad 
\vartheta \cdot B(q)\nabla u = 0 \text{ on } \Gamma_N,
\]
respectively (in weak formulation), 
\begin{equation}
\label{eq:randomeq-36467}
\text{Find }
u \in H_D^1(\Omega)
\text{ such that }
\left (B(q)\nabla u, \nabla \phi \right )_{L^2(\Omega)^d}
=
\left \langle f, \phi\right \rangle_{ H_D^1(\Omega)}
~~\forall \phi \in H_D^1(\Omega),
\end{equation}
where
$q \in \Qad$;
$B\colon \R \to \R^{d \times d}$ is
a continuously differentiable function 
satisfying 
\begin{equation}
\label{eq:randomeq-5242}
\exists \beta > 0\colon \qquad 
\beta |\zeta|^2 \leq B(s)\zeta \cdot \zeta \quad \forall \zeta \in \R^d
\quad \forall s \in [a,b];
\end{equation}
$H_D^1(\Omega) $ is defined by 
$ H_D^1(\Omega) := \{\phi \in H^1(\Omega) \mid 
\phi = 0 \text{ on } \Gamma_D\}$; and $f \in H_D^1(\Omega)^*$.
Then, it follows from the 
continuity of the trace operator that 
$H_D^1(\Omega)$ is 
a closed subspace of $H^1(\Omega)$ 
and, owing to the definition of the $H^1(\Omega)$-norm, that  
\[
V := \left \{
v \in L^2(\Omega)^{d+1}
\mid 
\exists u \in H_D^1(\Omega)\colon 
u = v_1 \text{ and } \partial_i u = v_{i+1} \forall i=1,\ldots, d
\right \}
\]
is a closed subspace of $L^2(\Omega)^{d+1}$
isometrically isomorphic to $( H_D^1(\Omega), \|\cdot\|_{H^1(\Omega)})$. 
In particular, we may rewrite 
\eqref{eq:randomeq-36467} as 
\[
\text{Find }
v\in V \text{ such that }
\left (A(q)v, w \right )_{L^2(\Omega)^m} = F(w)\quad \forall w \in V
\]
with 
\begin{equation}
\label{eq:randomeq626738}
m:= d+1,\qquad 
A(s) := \begin{pmatrix}
0 & 0 \\ 0 &B(s)
\end{pmatrix}
\quad \forall s \in \R
\qquad
\text{and}
\qquad 
F(w)
:=
\left \langle f, w_1\right \rangle_{ H_D^1(\Omega)}\quad \forall w \in V.
\end{equation}
From Friedrichs' inequality 
\cite[Lemma 2.5]{Troeltzsch2010},
\eqref{eq:randomeq-5242},
and the definition of $\Qad$ in 
\cref{ass:standing:abstract}\ref{ass:standing:abstract:iv},
it follows further that
the operator $A$ in \eqref{eq:randomeq626738} satisfies 
\begin{equation*}
\begin{aligned}
\left ( 
A(q)v, v
\right )_{L^2(\Omega)^m}
&=
\left ( 
B(q)[v_i]_{i=2}^{d+1} , [v_i]_{i=2}^{d+1}
\right )_{L^2(\Omega)^d}
=
\left ( 
B(q)\nabla v_1, \nabla v_1
\right )_{L^2(\Omega)^d}
\\
&\geq
\beta \| \nabla v_1\|_{L^2(\Omega)^d}^2
\geq 
\beta c
\| v_1\|_{H^1(\Omega)}^2
= 
\beta c
\|v \|_{L^2(\Omega)^m}^2
\qquad \forall v \in V
\qquad \forall q \in \Qad
\end{aligned}
\end{equation*}
for some constant
$c = c(\Omega, \Gamma_D) > 0$.
By collecting everything, we now obtain that
the isometric 
isomorphism
$( H_D^1(\Omega), \|\cdot\|_{H^1(\Omega)})
\cong 
(V, \|\cdot\|_{L^2(\Omega)^m})$
allows us to reformulate 
the scalar 
elliptic partial differential  equation 
\eqref{eq:randomeq-36467}
as a variational identity that fits precisely into the framework of
the abstract problem 
\eqref{eq:P}. 
This shows in particular that 
the 
abstract framework
in  \eqref{eq:P} 
indeed covers our model problem 
\eqref{eq:model_problem_elliptic}
and the PDEs considered in \cref{sec:homogen};
see also  \cref{sec:num_exp}.
\end{example}

\begin{example}[Linear elasticity]
\label{ex:elastic}
Let $\Omega$, $d$, $\Gamma_N$, $\Gamma_D$,
$a$, $b$, $\vartheta$, and $H_D^1(\Omega)$
be as in \cref{ex:scalar}
and 
assume that 
$B \colon \R \to \LLLL(\R^{d \times d}, \R^{d \times d}) \cong \R^{d\times d \times d \times d}$, 
$s \mapsto B(s) = (B_{ij}^{kl}(s))_{i,j,k,l =1}^d$,
is a continuously differentiable tensor map satisfying
$B_{ij}^{kl}(s) = B_{ji}^{kl}(s)$
for all $i,j,k,l$ and $s \in \R$
as well as 
\begin{equation}
\label{eq:randomeq-52455-2}
\exists \beta > 0\colon \qquad 
\beta |Z|^2 \leq (B(s) Z) : Z
\quad \forall Z \in \R^{d \times d}_{\sym}
\quad \forall s \in [a,b],
\end{equation}
where $|\cdot|\colon \R^{d \times d} \to [0, \infty)$
is the Frobenius norm
and where 
\[
B(s) Z
:=
\left (\sum_{k,l = 1}^d
B_{ij}^{kl}(s) Z_{kl}\right )_{i,j=1}^d
\qquad 
\text{and}
\qquad
Z : W := \sum_{i,j=1}^d Z_{ij}W_{ij}\qquad 
\forall Z,W \in \R^{d \times d}.
\]
Consider the steady-state linear elasticity problem 
\[
-\nabla\cdot ( B(q)\nabla_\sym u)=f \text{ in } \Omega,
\qquad
u=0 \text{ on } \Gamma_D,
\qquad 
\vartheta \cdot B(q)\nabla_\sym u  = 0 \text{ on } \Gamma_N,
\]
respectively (in weak formulation), 
\begin{equation}
\label{eq:randomeq-36467-2}
\text{Find }
u \in H_D^1(\Omega)^d
\text{ such that }
\int_\Omega 
 ( B(q)\nabla_\sym u ) : \nabla_\sym \phi
\d x
=
\left \langle f, \phi\right \rangle_{ H_D^1(\Omega)^d}
~~ \forall \phi \in H_D^1(\Omega)^d,
\end{equation}
where $q \in \Qad$, $f \in [ H_D^1(\Omega)^d ]^*$, and 
\[
\nabla_\sym \phi := \frac12 \left ( \nabla \phi + (\nabla \phi)^\top \right )
\in L^2(\Omega)^{d \times d}
\qquad \forall \phi \in H_D^1(\Omega)^d.
\]
Then, the space 
\[
V := \left \{
v \in L^2(\Omega)^{d(d+1)}
\mid 
\exists u \in H_D^1(\Omega)^d\colon 
u_j = v_{1 + (j-1)(d+1)},
\partial_i u_j = v_{1 + (j-1)(d+1) + i} 
\forall i, j =1,\ldots, d\right \}
\]
is isometrically isomorphic 
to $H_D^1(\Omega)^d$ by the definition of 
$\|\cdot\|_{H^1(\Omega)^d}$ and,
for every $s \in \R$, there exists a 
unique 
$A(s) \in \R^{d(d+1) \times d(d+1)}$
satisfying 
\begin{equation*}
\begin{aligned}
&[A(s)v ]^\top w
\\
&=
\frac14
 \sum_{i,j,k,l = 1}^d
 B_{ij}^{kl}(s) 
 \left ( v_{1+(k-1)(d+1)+l} + v_{1+(l-1)(d+1)+k} \right )
 \left ( w_{1+(i-1)(d+1)+j} + w_{1+(j-1)(d+1)+i} \right )
\end{aligned}
\end{equation*}
for all $v, w \in \R^{d(d+1)}$.
Owing to the definition of $A(s)$ and the 
properties of $B \colon \R \to  \R^{d\times d \times d \times d}$,
the mapping $A\colon \R \to \R^{d(d+1) \times d(d+1)}$,
$s \mapsto A(s)$,
is continuously differentiable and satisfies 
\begin{equation*}
\left (A(q)v, w \right )_{L^2(\Omega)^{d(d+1)}}
=
\int_\Omega
\left (
B(q)\nabla_\sym [v_{1 + (j-1)(d+1)}]_{j =1}^d
\right )
:
\nabla_\sym [w_{1 + (j-1)(d+1)}]_{j =1}^d
\d x 
\end{equation*}
for all $v,w \in V$ and $q \in \Qad$.
Analogously to \cref{ex:scalar},
the above allows us to restate \eqref{eq:randomeq-36467-2}
as 
\[
\text{Find }
v\in V \text{ such that }
\left (A(q)v, w \right )_{L^2(\Omega)^m} = F(w)~~\forall w \in V,
\]
with 
\begin{equation*}
m:= d(d+1)
\qquad
\text{and}
\qquad 
F(w)
:=
\left \langle f, [w_{1 + (j-1)(d+1)}]_{j =1}^d\right \rangle_{ H_D^1(\Omega)^d}\quad \forall w \in V.
\end{equation*}
Note that, from Korn's inequality \cite[Theorem 6.15-4]{Ciarlet2013},
\eqref{eq:randomeq-52455-2}, the assumptions on $a$ and $b$
and the definition of $\Qad$, it follows that 
\begin{equation*}
\begin{aligned}
\left ( 
A(q)v, v
\right )_{L^2(\Omega)^m}
&=
\int_\Omega
\left (
B(q)\nabla_\sym [v_{1 + (j-1)(d+1)}]_{j =1}^d
\right )
:
\nabla_\sym [v_{1 + (j-1)(d+1)}]_{j =1}^d
\d x 
\\
&\geq
\beta 
\int_\Omega 
\left | \nabla_\sym [v_{1 + (j-1)(d+1)}]_{j =1}^d \right |^2
\d x
\\
&\geq 
\beta c
\| [v_{1 + (j-1)(d+1)}]_{j =1}^d\|_{H^1(\Omega)^d}^2
= 
\beta c
\|v \|_{L^2(\Omega)^m}^2
\qquad \forall v \in V
\qquad \forall q \in \Qad
\end{aligned}
\end{equation*}
holds for some constant $c = c(\Omega, \Gamma_D) > 0$.
In combination with what was shown before,
this demonstrates that the 
PDE 
\eqref{eq:randomeq-36467-2} also falls under the scope 
of the setting of \eqref{eq:P}
after invoking the isometric isomorphism 
$( H_D^1(\Omega)^d, \|\cdot\|_{H^1(\Omega)^d})
\cong 
(V, \|\cdot\|_{L^2(\Omega)^{d(d+1)}})$
and that all results established for \eqref{eq:P}
can also be used when studying the 
identification of a unidirectionally distributed 
parameter $q \in \Qad$
in the field of steady-state linear elasticity.
\end{example}

\begin{remark}~
\begin{itemize}
    \item If $a$ and $b$ are constants (as in \cref{ex:scalar,ex:elastic}),
    then it is, of course, sufficient that the map 
    $A$
    in \cref{ass:standing:abstract}\ref{ass:standing:abstract:v}
    is defined on the interval $[a,b]$
    as no values of $A$ outside of $[a,b]$
    are needed in \eqref{eq:P}. 
    Since every $C^1$-mapping 
    $A\colon [a,b] \to \R^{m \times m}$
    with $-\infty < a < b < \infty$ admits
    a $C^1$-extension to the whole of $\R$,
    the assumption that $A$ is defined and continuously 
    differentiable on the entire real line 
    in \cref{ass:standing:abstract} does not cause 
    any loss of generality. 

    \item 
    A purely one-dimensional situation 
    (i.e.\ a problem living solely on an interval $(0,L)$)
    can be transformed 
    into a setting compatible
    with \cref{ass:standing:abstract} by adding an 
    artificial dimension. If we are given 
    a differential equation of the form \eqref{eq:randomeq-36467}
    on an interval $(0,L)$, i.e.\
    \begin{equation*}    
    \text{Find }
    u \in H_D^1(0,L)
    \text{ such that }
    \left (B(q) u', \phi' \right )_{L^2(0,L)}
    =
    \left \langle f, \phi\right \rangle_{ H_D^1(0,L)}
    ~~\forall \phi \in H_D^1(0,L),
    \end{equation*}
    for example,
    then one can simply consider the 
    space 
    \begin{equation*}
    \begin{aligned}
    V := \Big\{
    v \in L^2( (0,L) \times (0,1) )^{2}
     \mid ~&
    \exists u \in H_D^1(0,L)\colon 
    u(x_1) = v_1(x_1,x_2) 
    \text{ and }
    \\
    & u'(x_1) = v_{2}(x_1, x_2) ~\text{ f.a.a.\ } 
    (x_1, x_2) \in (0,L) \times (0,1)
\Big \}
\end{aligned}
\end{equation*}
in the transformation step described in 
\cref{ex:scalar} to arrive at a 
variational problem in $L^2( (0,L) \times (0,1) )^{2}$
with $D := (0,1)$ that fits into the framework of \eqref{eq:P}.
The assumption that the 
$D$-component exists in \cref{ass:standing:abstract}
is thus without loss of generality (w.l.o.g.).
\item If a partial differential equation on a 
non-cylindrical domain $\Omega$ is given, 
then one can again compensate this by adjusting the definition of 
the space $V$ in the transformation step. If, for example, we consider a PDE 
\[
-\nabla\cdot (B(q)\nabla u)=f \text{ in } \Omega,
\qquad
u=0 \text{ on } \Gamma_D,
\qquad 
\vartheta \cdot B(q)\nabla u = 0 \text{ on } \Gamma_N,
\]
as in \cref{ex:scalar} defined on an arbitrary 
bounded Lipschitz domain $\Omega \subset \R^d$, $2 \leq d \in \N$,
satisfying $\Omega \subset (0,L) \times \R^{d-1}$,
then, by choosing a bounded Lipschitz domain 
$D \subset \R^{d-1}$ satisfying $\Omega \subset (0,L) \times D := \hat \Omega$
and by considering the space 
\begin{equation*}
\begin{aligned}
V := \Big\{
v \in L^2(\hat \Omega)^{d+1}
\mid ~
& 
v = 0 \text{ a.e.\ in } \hat \Omega \setminus \Omega\text{ and }
\\
& \exists u \in H_D^1(\Omega)\colon 
u = v_1 \text{ and } \partial_i u = v_{i+1} \text{ a.e.\ in } 
\Omega~\forall i=1,\ldots, d \Big \}
\end{aligned}
\end{equation*}
in the transformation steps in \cref{ex:scalar},
one ends up with a problem of the type \eqref{eq:P}
on the cylindrical domain $\hat \Omega = (0,L) \times D$
that satisfies all of the conditions in \cref{ass:standing:abstract}.
Note that the assumption 
that $q$ varies in $x_1$-direction is also 
w.l.o.g.\ in this context as we can always rotate the coordinate system 
before starting the reformulation. 
Problems in which the parameter 
$q$ varies in an arbitrary fixed direction
are thus also covered.
\item Analogously to \cref{ex:scalar,ex:elastic},
one shows that the setting 
of \eqref{eq:P} also covers higher-order 
elliptic PDEs (e.g., fourth-order elliptic equations 
weakly formulated in $H_0^2(\Omega)$).
\end{itemize}
\end{remark}
%\pagebreak

Having discussed the scope of our abstract framework and that
it indeed covers the applications in \cref{sec:intro,sec:homogen}, we
now 
turn our attention to the analysis of \eqref{eq:P}. 
We start by proving well-posedness 
and by deriving 
first-order necessary optimality conditions.

\subsection{Well-posedness and first-order necessary optimality conditions}
\label{subsec:well-posedness_etc}
To prove the solvability of \eqref{eq:P},
we analyse the solution map of its governing variational 
identity.

\begin{lemma}[Unique solvability of the governing variational identity]
\label{lem:well-posedness_state}
For every $q \in \Qad$, there exists a unique solution 
$v \in V$ of the variational problem 
\begin{equation}
\label{eq:var_id}
\text{Find }
v\in V \text{ such that }
\left (A(q)v, w \right )_{L^2(\Omega)^m} = F(w)\quad \forall w \in V.
\end{equation}
This solution satisfies 
\begin{equation}
\label{eq:LM_estimate}
\|v\|_{L^2(\Omega)^m} \leq \frac{1}{\alpha} \|F\|_{V^*}. 
\end{equation}
\end{lemma}
\begin{proof}
As $V$ is a closed subspace of $L^2(\Omega)^m$, it is Hilbert,
and as $Q \hookrightarrow L^\infty(\Omega)$
and $A\colon \R \to \R^{m \times m}$ is $C^1$
(and thus bounded on bounded sets), we have 
\[
\left | 
\left ( 
A(q)v, w
\right )_{L^2(\Omega)^m}
\right |
\leq
C(q)
\|v\|_{L^2(\Omega)^m}
\|w\|_{L^2(\Omega)^m}
\qquad 
\forall v,w \in V
\]
with a constant $C(q) > 0$ for all $q \in Q$. 
In view of \eqref{eq:elliptic_abs},  
the assertion of the lemma 
now follows immediately 
from the lemma of Lax--Milgram
\cite[Lemma 2.2]{Troeltzsch2010}.
\end{proof}

\begin{definition}[Solution map $S$]
We define $S\colon \Qad \to V$, $q \mapsto v$, to be the 
solution map of \eqref{eq:var_id}.
\end{definition}

\begin{lemma}[Local Lipschitz continuity of $S$]
\label{lem:SlocLip}
Define 
\[
\Lip(A,[c_1, c_2])
:=
\inf
\left \{
C \geq 0
\mid
\|A(s_1) - A(s_2)\|_{\LLLL(\R^m, \R^m)} \leq C |s_1 - s_2|~
\forall s_1, s_2 \in \R\cap [c_1, c_2]\right \}
\]
for all $-\infty \leq c_1 \leq c_2 \leq \infty$.
Then, it holds 
$\Lip(A,[c_1, c_2]) < \infty$ for all 
$-\infty < c_1 \leq c_2 < \infty$ and,
for all  $q_1, q_2 \in \Qad$
satisfying $c_1 \leq q_i \leq c_2$
a.e.\ in $(0,L)$ for $i=1,2$
and $-\infty \leq c_1 \leq c_2 \leq \infty$,
it holds
\begin{equation*}
\|S(q_1) - S(q_2)\|_{L^2(\Omega)^m}
\leq
\frac{\Lip(A,[c_1, c_2])
\|F\|_{V^*}}{\alpha^2} \,
\|q_1 - q_2\|_{L^\infty(0,L)}.
\end{equation*}
\end{lemma}

\begin{proof}
The continuous differentiability
of $A\colon \R \to \R^{m\times m}$
implies that
$\Lip(A,[c_1, c_2]) < \infty$ holds for all 
$-\infty < c_1 \leq c_2 < \infty$.
Suppose now that
$q_1, q_2 \in \Qad$
are as in the assertion of the lemma. 
Define $v_i := S(q_i)$, $i=1,2$. Then, 
\eqref{eq:var_id} yields 
\[
\left (A(q_1)v_1, w \right )_{L^2(\Omega)^m} 
= F(w)
=
\left (A(q_2)v_2, w \right )_{L^2(\Omega)^m} 
\quad \forall w \in V.
\]
In view of \eqref{eq:elliptic_abs},
the definition of $\Lip(A,[c_1, c_2])$ 
and \eqref{eq:LM_estimate},
this implies 
\begin{equation}
\label{eq:randomeq26367}
\begin{aligned}
\alpha \|v_1 - v_2\|_{L^2(\Omega)^m}^2
&\leq
\left (A(q_1)(v_1 - v_2), v_1 - v_2 \right )_{L^2(\Omega)^m}
\\
&=
\left (A(q_2)v_2 - A(q_1)v_2, v_1 - v_2 \right )_{L^2(\Omega)^m}
\\
&\leq
\Lip(A,[c_1, c_2])
\|q_1 - q_2\|_{L^\infty(0,L)}
\|v_2\|_{L^2(\Omega)^m}
\|v_1 - v_2\|_{L^2(\Omega)^m}
\\
&\leq
\frac{\Lip(A,[c_1, c_2])
\|F\|_{V^*}}{\alpha}
\|q_1 - q_2\|_{L^\infty(0,L)}
\|v_1 - v_2\|_{L^2(\Omega)^m}.
\end{aligned}
\end{equation}
The assertion now follows immediately. 
\end{proof}

\begin{lemma}[Complete continuity of $S$]
\label{lem:complete_continuity}
For every sequence 
$\{q_k \} \subset \Qad$
satisfying $q_k \weakly q$ in $Q$ for some $q \in Q$,
it holds $q \in \Qad$ and $S(q_k) \to S(q)$ in $L^2(\Omega)^m$.
\end{lemma}

\begin{proof}
The inclusion $q \in \Qad$ follows
from $Q \compactly L^\infty(0,L)$
and the definition of $\Qad$. The convergence 
$S(q_k) \to S(q)$ in $L^2(\Omega)^m$
is an immediate consequence of 
\cref{lem:SlocLip}
and $Q \compactly L^\infty(0,L)$.
\end{proof}

With the last lemmas at hand, the solvability 
of \eqref{eq:P} follows along standard lines:

\begin{theorem}[{Solvability of \eqref{eq:P}}]
\label{th:solvabilityP}
Problem \eqref{eq:P} has at least one global solution
$(\bar v, \bar q) \in V \times \Qad$.
\end{theorem}

\begin{proof}
To prove this theorem, we use the 
direct method of the calculus of variations:
As $\Qad \neq \emptyset$, we can find 
$\{q_k\} \subset \Qad$ such that 
\[
J(S(q_1), q_1) \geq J(S(q_k), q_k) \to \inf_{q \in \Qad} J(S(q), q)
\quad \text{ for }k\to \infty.
\]
Using \cref{ass:standing:abstract}\ref{ass:standing:abstract:vii},
the above yields that
\[
J(S(q_1), q_1) \geq J(S(q_k), q_k) \geq 
\rho(\|q_k\|_Q)\qquad \forall k \in \N
\]
and, in view of the properties of $\rho$, that 
$\{q_k\}$ is bounded in $Q$. As $Q$ is assumed to be a 
reflexive Banach space, the latter allows us to pass to a subsequence 
(still denoted by the same symbol) such that $q_k \weakly \bar q$ holds in $Q$
for some $\bar q \in Q$. By \cref{lem:complete_continuity},
we have $\bar q \in \Qad$ and $S(q_k) \to S(\bar q)$ in $V$. 
In combination with the lower semicontinuity property
of $J$ in \cref{ass:standing:abstract}\ref{ass:standing:abstract:vii},
this entails 
\[
J(S(\bar q), \bar q)
\leq
\liminf_{k \to \infty}
J(S(q_k), q_k)
=
\inf_{q \in \Qad} J(S(q), q).
\]
Thus, 
$(S(\bar q), \bar q)$ is a global minimiser of \eqref{eq:P} and the proof is
complete.
\end{proof}

Next, we turn our attention to first-order necessary optimality conditions for \eqref{eq:P}. 
To derive these, we require differentiability results for the 
solution mapping $S$. 

\begin{lemma}[Continuous
%Fréchet
differentiability of the 
operator in the governing equation]
\label{lem:aux_dif_op}
Define 
\begin{equation}
\label{eq:FFdef}
\FF\colon V \times Q \to V^*,
\qquad 
\left \langle \FF(v,q), w \right \rangle_V
:=
\left (A(q)v, w \right )_{L^2(\Omega)^m} - F(w)
\quad \forall w \in V.
\end{equation}
Then, $\FF$ is continuously
%Fréchet
differentiable 
with partial derivatives 
\begin{align*}
&\partial_v \FF\colon 
V \times Q \to \LLLL(V, V^*),
&&
\left \langle  \partial_v \FF(v,q)\delta v  , w\right \rangle_{V}
=
\left (A(q)\delta v, w \right )_{L^2(\Omega)^m}
&&\forall w \in V\quad \forall\delta v \in V,
\\
&\partial_q \FF\colon 
V \times Q \to \LLLL(Q, V^*),
&&
\left \langle  
\partial_q \FF(v,q)\delta q , w\right \rangle_{V}
=
\left (A'(q) \delta q\, v , w \right )_{L^2(\Omega)^m}
&&\forall w \in V\quad \forall \delta q \in Q.
\end{align*}
\end{lemma}
% end \Com
\begin{proof}
The operator $\FF$ in \eqref{eq:FFdef}
can be written as $\FF(v,q)=\FF_1(\FF_2(\FF_3(q)),v)-F$ with
the smooth bilinear map $\FF_1\colon 
L^\infty(\Omega)^{m\times m} \times V
\ni
(B,v) 
\mapsto Bv\in L^2(\Omega)^m \hookrightarrow V^*$,
the Nemytskii operator
$\FF_2\colon
L^\infty(\Omega) \ni q \mapsto A(q(\cdot))\in L^\infty(\Omega)^{m\times m}$,
and the continuous embedding 
$\FF_3\colon Q \ni q \mapsto q\in L^\infty(\Omega)$.
Since $A\colon \R\to\R^{m\times m}$ is $C^1$, we have that
$A'$ is uniformly continuous on bounded sets. 
Thus, condition
(UC) on p.~133 of \cite{GKT92} is satisfied,
and we can
apply \cite[Theorem 7, Remark 4]{GKT92} to obtain that
$\FF_2$ is continuously differentiable with
\begin{equation}\label{eq:F2deriv}
(\FF_2'(q)\delta q)(x)=A'(q(x))\delta q(x)
\qquad
\text{f.a.a.\ }x\in \Omega
\qquad 
\forall q,
\delta q\in L^\infty(\Omega).
\end{equation}
The continuous differentiability of 
$\FF$ and the formulas for the 
partial 
derivatives $\partial_v \FF$ and $\partial_q \FF$ now follow easily 
from the chain rule \cite[Proposition 3.4.4]{Schirotzek2007}.
\end{proof}

Using the implicit function theorem, we can now prove:

\begin{lemma}[Continuous
%Fréchet
differentiability of $S$]
\label{lem:Sfrechet}
There exists 
an 
open neighbourhood $O \subset Q$ of $\Qad$
such that 
the solution mapping 
$S\colon \Qad \to V$, $q \mapsto v$,
of \eqref{eq:var_id}
can be extended 
to a continuously
%Fréchet
differentiable function $S\colon O \to V$,
such that 
$\partial_v \FF(S(q),q) \in \LLLL(V,V^*)$
is an isomorphism for all $q \in O$
for the function $\FF$ in \eqref{eq:FFdef},
and such that, 
for every $q \in O$,
the derivative $\delta v := S'(q)\delta q \in V$
of $S$ at $q$
in direction
$\delta q \in Q$ is uniquely characterised by the 
problem 
\begin{equation}
\label{eq.DirDifVE}
\text{Find }
\delta v\in V
\text{ such that }
\left (A(q)\delta v, w \right )_{L^2(\Omega)^m} = 
-
\left (A'(q)\delta q S(q), w \right )_{L^2(\Omega)^m} 
~~\forall w \in V.
\end{equation}
\end{lemma}
\begin{proof}
Let $\tilde q \in \Qad$ be arbitrary but fixed.
Then, it holds $\FF(S(\tilde q), \tilde q) = 0 \in V^*$
by the definitions of $\FF$ and $S$, and 
$\partial_v \FF(S(\tilde q),\tilde q)\colon V \to V^*$
is an isomorphism by the lemma of Lax--Milgram
(cf.\ the proof of \cref{lem:well-posedness_state}).
Owing to \cref{lem:aux_dif_op} and the implicit function theorem 
\cite[Theorem 3.7.2]{Schirotzek2007},
this implies that there exists an 
open neighbourhood $O_{\tilde q}$ of $\tilde q$
such that 
$S$ admits a continuously 
%Fréchet
differentiable extension to $O_{\tilde q}$,
such that $\partial_v \FF(S(q),q) \in \LLLL(V,V^*)$
is an isomorphism for all $q \in O_{\tilde q}$,
 and such that, for all $q \in O_{\tilde q}$,
it holds 
$S'(q) = - \partial_v \FF(S(q),q)^{-1}\partial_q \FF(S(q),q)$.
Since 
$\tilde q \in \Qad$ was arbitrary 
and owing to the formulas for 
$ \partial_v \FF$
and
$ \partial_q \FF$
in \cref{lem:aux_dif_op},
the assertions of the lemma now follow immediately with
$O=\bigcup_{\tilde q\in \Qad} O_{\tilde q}$. 
\end{proof}

For the objective function of \eqref{eq:P}, 
the last result implies:

\begin{lemma}[Derivative of the reduced objective function]
\label{lem:derredObjective}
Let $O \subset Q$ be the open 
neighbourhood from \cref{lem:Sfrechet}. Define
\[
\hat J\colon O \to \R,
\qquad 
\hat J(q) := J(S(q),q).
\]
Then, 
$\hat J$
is Fréchet differentiable on $O$.
Further, for every $q \in O$
with associated $v := S(q)$
and every $\delta q \in Q$, it holds 
\begin{equation}
\label{eq:redJderivative}
 \langle \hat J'(q), \delta q  \rangle_Q
=
\left \langle \partial_q J(v,q), \delta q \right \rangle_Q
+
\int_0^L
p \,\delta q \d s,
\end{equation}
where $p \in L^1(0,L)$
is defined by 
\begin{equation}
\label{eq:pdef}
p(x_1)
:=
\int_D
\big [ A'(q(x_1)) v(x)
\big ]^\top z(x)
\d (x_2,\ldots,x_{n+1})
\quad
\text{f.a.a.\ } x_1 \in (0,L)
\end{equation}
with $z$ as the unique solution of 
\begin{equation}
\label{eq:adjointequation}
\text{Find }
z \in V \text{ such that }
\left (A(q)^\top z, w \right )_{L^2(\Omega)^m} = 
-
 \left \langle \partial_v J(v,q), w \right \rangle_V
~~\forall w \in V.
\end{equation}

\end{lemma}

\begin{proof}
The Fréchet differentiability of $\hat J$
follows from our assumptions 
on $J$, \cref{lem:Sfrechet} and the 
chain rule \cite[Proposition 3.4.4]{Schirotzek2007}.
Suppose now that
some 
$q \in O$ with associated 
$v :=S(q) \in V$ is given and let 
$\FF$ be as in \eqref{eq:FFdef}. 
Then, we know from the properties of $O$
that $\partial_v \FF( v, q) \in \LLLL(V,V^*)$
is an isomorphism. 
This implies that the adjoint 
$\partial_v \FF(v,q)^* \in \LLLL(V,V^*)
\cong \LLLL(V^{**}, V^*)$
given by 
\[
\left \langle 
\partial_v \FF(v,q)^* w_1,
w_2
\right \rangle_V
:=
\left \langle 
\partial_v \FF(v,q) w_2,
w_1
\right \rangle_V
=
\left (
A(q)^\top w_1, w_2
\right )_{L^2(\Omega)^m}
\qquad
\forall w_1, w_2 \in V
\]
is an isomorphism, too. 
In particular, there is a 
unique $z \in V$ such that 
$\partial_v \FF(v,q)^* z = - 
\partial_v J(v,q)$ holds, 
i.e.\ such that $z$ satisfies \eqref{eq:adjointequation}.
Suppose now that $\delta q \in Q$ is arbitrary.
Then, the chain rule,
\eqref{eq:adjointequation}
with $w:= S'(q)\delta q \in V$,
\eqref{eq.DirDifVE} with $w:=z$
and the Theorem of Fubini
imply that 
\begin{equation*}
\begin{aligned}
\langle \hat J'(q), \delta q \rangle_Q
&=
\left \langle \partial_v J(v,q), S'(q)\delta q \right \rangle_V
+
\left \langle \partial_q J(v,q), \delta q \right \rangle_Q
\\
&=
-
\left (A(q)^\top z, S'(q)\delta q \right )_{L^2(\Omega)^m}
+
\left \langle \partial_q J(v,q), \delta q \right \rangle_Q
\\
&=
-
\left (A(q)S'(q)\delta q, z \right )_{L^2(\Omega)^m}
+
\left \langle \partial_q J(v,q), \delta q \right \rangle_Q
\\
&=
\int_\Omega
\left [ A'(q)\delta q S(q)\right ]^\top z \d x 
+
\left \langle \partial_q J(v,q), \delta q \right \rangle_Q
\\
&=
\int_0^L
\int_D
\big [ A'(q(x_1)) v(x)
\big ]^\top z(x)
\d (x_2,\ldots,x_{n+1})
\,
\delta q(x_1) \d x_1
+
\left \langle \partial_q J(v,q), \delta q \right \rangle_Q.
\end{aligned}
\end{equation*}
This establishes the formula 
for the directional derivative in \eqref{eq:redJderivative}
and completes the proof.
\end{proof}

As a consequence of 
\cref{lem:derredObjective}, we obtain:

\begin{theorem}[First-order necessary optimality conditions for 
\eqref{eq:P}]
\label{th:first_order_necessary_condition}
Suppose that $\bar q$ is locally optimal for \eqref{eq:P}, i.e.\ 
that $\bar q$ is a local minimiser of the reduced problem
\begin{equation*}
\label{eq:Pred}
\tag{P$_{\red}$}
	\begin{aligned}
		\text{Minimise} 
		\quad &  \hat J(q)  \\
        \text{w.r.t.}
        \quad & q \in \Qad
\end{aligned}
\end{equation*}
with $\hat J (q) := J(S(q),q)$ as in \cref{lem:derredObjective}. 
Then, $\bar q$ satisfies 
the Bouligand stationarity condition 
\begin{equation}
\label{eq:Bouligand}
\bar q \in \Qad,\qquad
\langle \hat J'(\bar q), q - \bar q \rangle_Q
\geq 0\quad \forall q \in \Qad
\end{equation}
and this condition is equivalent to the stationarity system
\begin{equation}
\label{eq:KKT}
\begin{aligned}
&
\bar v \in V,\qquad~ \left (A(\bar q)\bar v, w \right )_{L^2(\Omega)^m} = F(w)\qquad\qquad\qquad~~~~ \forall w \in V && \text{(state equation)}
\\
&
\bar z   \in V,\qquad~ 
  (A(\bar q)^\top \bar z, w  )_{L^2(\Omega)^m} = 
-
 \left \langle \partial_v J(\bar v,\bar q), w \right \rangle_V
\quad~\; \forall w \in V
&& \text{(adjoint equation)}
\\
&\left.
\begin{aligned}
&\bar q \in \Qad,
\quad
\left \langle \partial_q J(\bar v,\bar q), q-\bar q \right \rangle_Q
+
\int_0^L
\bar p (q - \bar q) \d s
\geq 0
\quad\, \forall q \in \Qad
\\
&\bar p \in L^1(0,L),
\quad
\bar p =
\int_D
\big [ A'(\bar q ) \bar v
\big ]^\top \bar z
\d (x_2,\ldots,x_{n+1})
\end{aligned}
\right \}
&& \text{(variational inequality)}.
\end{aligned}
\end{equation}
\end{theorem}
\begin{proof}
It follows straightforwardly 
from the Fréchet differentiability established for 
$\hat J$ in \cref{lem:derredObjective}
and the convexity of $\Qad$ that a local solution 
$\bar q$ of \eqref{eq:Pred} satisfies 
\eqref{eq:Bouligand}. The equivalence of
\eqref{eq:Bouligand} 
to the system 
\eqref{eq:KKT} 
is an immediate consequence of 
the formula \eqref{eq:redJderivative}.
\end{proof}

\begin{remark}\label{rem:highreg}
The system 
\eqref{eq:KKT} can be used 
to obtain higher regularity 
for the local minimisers
of \eqref{eq:Pred}. In 
the case $Q := H_0^1(0,L)$, 
$J(v,q) := J_1(v) + \frac\nu2 \|q\|_{H_0^1(0,L)}^2$
for some $J_1\colon V\to \R$
and 
$\nu > 0$, 
and 
$a,b \in \R$ with $ a < 0 <  b $,
for example, the variational inequality 
in \eqref{eq:KKT} becomes 
\[
\bar q \in \Qad,
\quad
\left ( \bar q, q-\bar q \right )_{H_0^1(0,L)}
+
\int_0^L
\frac{\bar p}{\nu}
(q - \bar q) \d s
\geq 0
\quad\, \forall q \in \Qad
\]
with 
$\Qad = \{q \in H_0^1(0,L) \mid a \leq q \leq b \text{ a.e.\ in }(0,L)\}$,
and it follows from 
$\bar p \in L^1(0,L)$ and standard regularity results 
for elliptic obstacle problems
that 
$\bar q \in W^{2,1}(0,L)$;
cf.\ \cite[Theorem 2.3]{ChristofWachsmuth2021SSC} and its proof.
Note that this 
$W^{2,1}(0,L)$-regularity is obtained independently 
of the underlying 
governing equation (cf.\ the PDEs in 
\cref{ex:scalar,ex:elastic}) and independently of 
the precise structure of $J_1$ (as long as $J_1$
is such that the conditions on $J$ in \cref{ass:standing:abstract} are satisfied,
of course).
\end{remark}

\subsection{A gradient-projection algorithm}
\label{subsec:gradproj}

We are now in the position to
construct and analyse a gradient-projection algorithm 
for the numerical solution of the problem 
\eqref{eq:P} (respectively, \eqref{eq:Pred}) in function space.
As already mentioned, the following results are rather standard 
from the viewpoint of optimisation theory; see 
\cite{Blank2017} and the references therein. 
We discuss them in detail,
on the one hand, because
we will use the gradient-projection algorithm 
(which can be seen as the state of the art in the field 
of inverse homogenisation; see \cite{Loch:2022,Loch:2024,Lipton2006,All:2019}) as a point of reference 
in the numerical experiments in \cref{sec:num_exp} and, 
on the other hand, since several of the following 
estimates etc.\ are also of relevance 
for our semismooth Newton approach; see 
\cref{cor:global_convergence_tracking_type,cor:invert_Newton_large_Tikhonov}.

Naturally, to talk about gradients and 
projections, we have to restrict the analysis to the 
Hilbert space case. Throughout this subsection, 
we work with the following standing assumptions:

\begin{assumption}[{Standing assumptions 
for \cref{subsec:gradproj}}]
\label{ass:standing_grad_proj}
In addition to the conditions in \cref{ass:standing:abstract},
we assume that:
\begin{enumerate}[label=\roman*)]
\item $Q$ is a Hilbert space with Riesz isomorphism $\iota_R\colon Q^* \to Q$;
\item There exists a constant $C_1 > 0$ such that
\begin{equation*}
\begin{aligned}
&\left \| \partial_v J(v, q) \right \|_{V^*}
\leq
C_1 
\qquad 
\forall v \in V,
\|v\|_V \leq \frac{1}{\alpha} \|F\|_{V^*},
\qquad 
\forall q \in \Qad;
\end{aligned}
\end{equation*}
\item There exists a constant $C_2 > 0$ such that
\begin{equation*}
\begin{aligned}
&\left \| \partial_v J(v_1, q_1) -  \partial_v J(v_2, q_2) \right \|_{V^*}
\leq
C_2 \left ( \|v_1 - v_2\|_V + \|q_1 - q_2\|_Q\right )
\\
&\hspace{1cm}
\forall v_1, v_2 \in V,
\|v_i\|_V \leq \frac{1}{\alpha} \|F\|_{V^*},
i=1,2,
\qquad 
\forall q_1, q_2 \in \Qad;
\end{aligned}
\end{equation*}
\item There exists a constant $C_3 > 0$ such that
\begin{equation*}
\begin{aligned}
&\left \| \partial_q J(v_1, q_1) -  \partial_q J(v_2, q_2) \right \|_{Q^*}
\leq
C_3 \left ( \|v_1 - v_2\|_V + \|q_1 - q_2\|_{Q}\right )
\\
&\hspace{1cm}
\forall v_1, v_2 \in V,
\|v_i\|_V \leq \frac{1}{\alpha} \|F\|_{V^*},
i=1,2,
\qquad 
\forall q_1, q_2 \in \Qad;
\end{aligned}
\end{equation*}
\item It holds 
\[
C_4 :=
 \Lip(A,[c_a, c_b]) < \infty,
 \qquad
 C_5 :=
\Lip \left (A', [ c_a, c_b] \right )
 < \infty,
\]
and
\[
 C_6 :=
 \sup_{s \in [c_a, c_b]}
 \|A'(s)\|_{\LLLL(\R^m, \R^m)}
 < \infty,
\]
where 
\[
c_a := \essinf_{s \in (0,L)} a(s)
\qquad\text{and}\qquad
c_b := \esssup_{s \in (0,L)} b(s),
\]
and where 
$\Lip(\cdot, [\cdot, \cdot])$
is defined as in \cref{lem:SlocLip}. 
\end{enumerate}
\end{assumption}
We further define 
\begin{equation}
\label{eq:C7_def}
C_7 := \|\iota\|_{\LLLL(Q, L^\infty(0,L))},
\end{equation}
where 
$\iota \in \LLLL( Q , L^\infty(0,L))$
denotes the embedding 
$Q \hookrightarrow L^\infty(0,L)$
from 
\cref{ass:standing:abstract}\ref{ass:standing:abstract:iii}.
As a starting point for our analysis, we prove:

\begin{lemma}[Lipschitz estimate for the derivative 
of the reduced objective function]
\label{lem:LipRedFunction}
The reduced objective function
$\hat J\colon O \to \R $ in \cref{lem:derredObjective} satisfies 
\begin{equation*}
\begin{aligned}
  \|
\hat J'(q_1)
-
\hat J'(q_2)
 \|_{Q^*}
\leq
\hat C \|q_1 - q_2\|_{Q}
\qquad \forall q_1, q_2 \in \Qad
\end{aligned}
\end{equation*}
with $\hat C$ given by 
\begin{equation}
\label{eq:QuellDerFreude}
\begin{aligned}
&\hat C :=
\Bigg (
\frac{C_3 C_4 C_7 \|F\|_{V^*} }{\alpha^2}
+
C_3
 +
 \frac{C_1 C_5 C_7^2 \|F\|_{V^*}}{\alpha^2} 
+ 
\frac{C_2 C_4 C_6 C_7^2\|F\|_{V^*}^2}{\alpha^4} 
\\
&\qquad\qquad\qquad\qquad
+ \frac{C_2 C_6 C_7 \|F\|_{V^*} }{\alpha^2}
+
2\frac{C_1C_4 C_6 C_7^2  \|F\|_{V^*} }{\alpha^3}
\Bigg ).
\end{aligned}
\end{equation}
\end{lemma}
\begin{proof}
Let $q_1, q_2 \in \Qad$ 
with associated 
$v_i := S(q_i)$, $i=1,2$,
be given,
and let $p_1, p_2 \in L^1(0,L)$
and $z_1, z_2 \in V$ be the quantities 
associated with $q_1, q_2$ in 
\cref{lem:derredObjective}.
Then, it follows from 
\eqref{eq:elliptic_abs}
and the same estimate as in 
the proof of 
\cref{lem:well-posedness_state}
that 
\[
\|z_i\|_{L^2(\Omega)^m} \leq \frac{C_1}{\alpha}
\qquad i=1,2,
\]
and we may calculate analogously 
to \eqref{eq:randomeq26367} that 
\begin{equation*}
\begin{aligned}
&\alpha \|z_1 - z_2\|_{L^2(\Omega)^m}^2
\\
&\leq
\left (A(q_1)^\top(z_1 - z_2), z_1 - z_2 \right )_{L^2(\Omega)^m}
\\
&=
\left (A(q_1)^\top z_1 - A(q_2)^\top z_2, z_1 - z_2 \right )_{L^2(\Omega)^m}
+
\left (A(q_2)^\top z_2 - A(q_1)^\top z_2, z_1 - z_2 \right )_{L^2(\Omega)^m}
\\
&=
\left \langle
- \partial_v J(v_1, q_1)
+
\partial_v J(v_2, q_2),
z_1 - z_2 \right \rangle_{V}
+
\left (A(q_2)^\top z_2  - A(q_1)^\top z_2, z_1 - z_2 \right )_{L^2(\Omega)^m}
\\
&\leq
\left (
C_2 \left ( \|v_1 - v_2\|_V + \|q_1 - q_2\|_{Q}\right )
+
\frac{C_1C_4}{\alpha}
\|q_1 - q_2\|_{L^\infty(0,L)}
\right )
\|z_1 - z_2\|_{L^2(\Omega)^m}.
\end{aligned}
\end{equation*}
In combination with 
\cref{lem:SlocLip}, this yields
\begin{equation}
\label{eq:zestimate}
\|z_1 - z_2\|_{L^2(\Omega)^m}
\leq
\left (
\frac{C_2 C_4 C_7\|F\|_{V^*}}{\alpha^3} \,
+ \frac{C_2}{\alpha}
+
\frac{C_1C_4C_7}{\alpha^2}
\right )
\|q_1 - q_2\|_Q.
\end{equation}
For the functions $p_1, p_2$, it now follows that 
\begin{equation*}
\begin{aligned}
\|p_1 - p_2\|_{L^1(0,L)}
&\leq
\int_\Omega 
\left |
\big [ A'(q_1 ) v_1
\big ]^\top z_1
-
\big [ A'(q_2 ) v_2
\big ]^\top z_2
\right |
\d x
\\
&\leq
\int_\Omega 
\left |
\big [ A'(q_1 ) v_1
\big ]^\top z_1
-
\big [ A'(q_2 ) v_1
\big ]^\top z_1
\right |
+
\left |
\big [ A'(q_2 ) v_1
\big ]^\top z_1
-
\big [ A'(q_2 ) v_2
\big ]^\top z_1
\right |
\\
&\qquad
+
\left |
\big [ A'(q_2 ) v_2
\big ]^\top z_1
-
\big [ A'(q_2 ) v_2
\big ]^\top z_2
\right |
\d x,
\end{aligned}
\end{equation*}
where 
\begin{equation*}
\begin{aligned}
&\int_\Omega 
\left |
\big [ A'(q_1 ) v_1
\big ]^\top z_1
-
\big [ A'(q_2 ) v_1
\big ]^\top z_1
\right |
\d x
\\
&\quad\leq
\|A'(q_1 ) - A'(q_2 )\|_{L^\infty(0,L;\LLLL(\R^m, \R^m))}
\|v_1\|_{L^2(\Omega)^m}
\|z_1\|_{L^2(\Omega)^m}
\\
&\quad\leq
\frac{C_1 C_5 C_7 \|F\|_{V^*}}{\alpha^2} 
\|q_1 - q_2\|_Q
\end{aligned}
\end{equation*}
and (due to \cref{lem:SlocLip})
\begin{equation*}
\begin{aligned}
 \int_\Omega 
\left |
\big [ A'(q_2 ) v_1
\big ]^\top z_1
-
\big [ A'(q_2 ) v_2
\big ]^\top z_1
\right |
\d x
 \leq
 \frac{C_1 C_4 C_6 C_7 \|F\|_{V^*}}{\alpha^3}
\|q_1 - q_2\|_{Q}
\end{aligned}
\end{equation*}
and (due to \eqref{eq:zestimate})
\begin{equation*}
\begin{aligned}
&\int_\Omega 
\left |
\big [ A'(q_2 ) v_2
\big ]^\top z_1
-
\big [ A'(q_2 ) v_2
\big ]^\top z_2
\right |
\d x
\\
&\quad \leq
\left (
\frac{C_2 C_4 C_6 C_7\|F\|_{V^*}^2}{\alpha^4} \,
+ \frac{C_2 C_6  \|F\|_{V^*} }{\alpha^2}
+
\frac{C_1C_4 C_6 C_7  \|F\|_{V^*} }{\alpha^3}
\right )
\|q_1 - q_2\|_Q.
\end{aligned}
\end{equation*}
From
\cref{lem:SlocLip,lem:derredObjective}
and
our assumptions on $\partial_q J$,
we now obtain that 
\begin{equation*}
\begin{aligned}
 \|
\hat J'(q_1)
-
\hat J'(q_2)
 \|_{Q^*}
 &\leq
  \|
\partial_q J(v_1,q_1)
-
\partial_q J(v_2,q_2)
 \|_{Q^*}
 +
 C_7 \|p_1 - p_2\|_{L^1(0,L)}
 \\
&\leq C_3 \left ( \|v_1 - v_2\|_V + \|q_1 - q_2\|_{Q}\right )
 +
  C_7 \|p_1 - p_2\|_{L^1(0,L)}
 \\
&\leq  \left ( 
\frac{C_3 C_4 C_7 \|F\|_{V^*} }{\alpha^2}
+
C_3
\right )
\|q_1 - q_2\|_{Q}
 +
  C_7 \|p_1 - p_2\|_{L^1(0,L)}
\\
&\leq  
\Bigg (
\frac{C_3 C_4 C_7 \|F\|_{V^*} }{\alpha^2}
+
C_3
 +
 \frac{C_1 C_5 C_7^2 \|F\|_{V^*}}{\alpha^2} 
+ 
\frac{C_2 C_4 C_6 C_7^2\|F\|_{V^*}^2}{\alpha^4} 
\\
&\qquad\qquad\qquad
+ \frac{C_2 C_6 C_7 \|F\|_{V^*} }{\alpha^2}
+
2\frac{C_1C_4 C_6 C_7^2  \|F\|_{V^*} }{\alpha^3}
\Bigg )
\|q_1 - q_2\|_Q
\\
&= \hat C \|q_1 - q_2\|_Q.
  \end{aligned}
\end{equation*}
This completes the proof. 
\end{proof}

Recall that the gradient $\nabla_Q \hat J(q)$
of the reduced objective function 
$\hat J\colon O \to \R$
in \cref{lem:LipRedFunction}
w.r.t.\ the inner product $(\cdot, \cdot)_Q$
is given by $\nabla_Q \hat J(q) := \iota_R(\hat J'(q))$
for all $q \in O$ and 
that $\proj_{\Qad}\colon Q \to Q$
denotes the projection in $Q$ 
onto the convex, closed and nonempty set $\Qad$.
Using these two quantities, we can reformulate the necessary 
optimality condition in \cref{th:first_order_necessary_condition}.

\begin{lemma}[Reformulation of the necessary optimality condition]
\label{lem:grad_opt_reform}
Let $\tau > 0$ be
arbitrary and let
$\bar q \in Q$ be given. 
Then, $\bar q$ is Bouligand stationary 
for \eqref{eq:P} in the sense of 
\eqref{eq:Bouligand}
if and only if 
\[
\bar q - \proj_{\Qad}
 (
\bar q - \tau \nabla_Q \hat J(\bar q)
  )
= 0.
\]
\end{lemma}
\begin{proof}
We have 
\begin{equation*}
\begin{aligned}
&\bar q \in \Qad,\qquad
\langle \hat J'(\bar q), q - \bar q \rangle_Q
\geq 0\quad \forall q \in \Qad
\\
&\iff\qquad
\bar q \in \Qad,\qquad
(\bar q - \bar q + \tau \nabla_Q \hat J(\bar q), q - \bar q )_Q
\geq 0\quad \forall q \in \Qad
\\
&\iff\qquad
\bar q = \proj_{\Qad}
 (
\bar q - \tau \nabla_Q \hat J(\bar q)
  ).
\end{aligned}
\end{equation*}
This proves the assertion.
\end{proof}

From \cref{lem:LipRedFunction}, we may now deduce the following classical 
decrease lemma:

\begin{lemma}[Decrease for gradient-projection steps]
\label{lem:descent_property}
Let $\hat J\colon O \to \R$ and $\hat C > 0$ be as in 
\cref{lem:LipRedFunction}.
Then, for all $\tau \in (0, 2/\hat C)$ and all 
$q \in \Qad$, we have 
\begin{equation}
\label{eq:descent_prop}
\hat J(q)
-
\hat J(\proj_{\Qad}
 (q - \tau \nabla_Q \hat J(q)))
 \geq
\left (
\frac{1}{\tau} - \frac{\hat C}{2}
\right )
\left \|
q - \proj_{\Qad}
 (
q - \tau \nabla_Q \hat J(q)
  )
\right \|_Q^2
\geq 0.
\end{equation}
\end{lemma}

\begin{proof}
Let $\tau \in (0, 2/\hat C)$ and  
$q \in \Qad$ be given.
Define 
$\hat q := \proj_{\Qad}
 (q - \tau \nabla_Q \hat J(q))$. 
 From the convexity of $\Qad$,
\cref{lem:LipRedFunction}
and
the mean value theorem 
 \cite[Theorem 3.2.6]{Drabek2007},
 we obtain that 
 \begin{equation}
 \label{eq:randomeq273h}
 \begin{aligned}
\hat J(\hat q)
-
\hat J(q)
&=
\int_0^1
\left ( \nabla_Q \hat J(q + s(\hat q - q)), \hat q - q\right )_Q
\d s
\\
&=
\left ( \nabla_Q \hat J(q), \hat q - q\right )_Q
+
\int_0^1
\left \langle  \hat J'(q + s(\hat q - q)) - \hat J'(q), \hat q - q\right \rangle_Q
\d s
\\
&\leq
\left ( \nabla_Q \hat J(q), \hat q - q\right )_Q
+
\int_0^1
\hat C s \|\hat q - q\|_Q^2 \d s
\\
&=
\left ( \nabla_Q \hat J(q), \hat q - q\right )_Q
+
\frac{\hat C}{2} \|\hat q - q\|_Q^2.
\end{aligned}
 \end{equation}
By the definitions of $\hat q$ and $\proj_{\Qad}$,
 it further holds that 
  \begin{equation}
  \label{eq:randomeq72hu3omn}
 - \frac1\tau\left \| \hat q - q \right \|_Q^2
-
 \left ( \nabla_Q \hat J(q) , \hat q - q\right )_Q
 =
 \frac1\tau
 \left ( \hat q - (q - \tau \nabla_Q \hat J(q)) , q - \hat q\right )_Q
\geq 0.
 \end{equation}
 Combining \eqref{eq:randomeq273h} with 
 \eqref{eq:randomeq72hu3omn}
 yields \eqref{eq:descent_prop} as desired. 
\end{proof}

\Cref{lem:descent_property} implies that 
gradient projection steps in $Q$
with a sufficiently small constant step size $\tau > 0$
produce a sequence of iterates with monotonously
decreasing objective-function values. This motivates the following algorithm.

\begin{Algorithm}[Gradient projection method for the solution of \eqref{eq:P}]\label{alg:gradproj}
~\hspace{-10cm}
\begin{algorithmic}[1]
\STATE{Choose an initial guess $q_0 \in \Qad$,
a step size $\tau > 0$ and a tolerance $\tol \geq 0$.}
\FOR{$k=0,1,2,3,\ldots$}
\STATE{
Calculate $v_k \in V$ by solving the problem 
\[
\text{Find }
v_k \in V
\text{ such that }\left (A(q_k)v_k, w \right )_{L^2(\Omega)^m} = F(w)\quad\forall w \in V.
\]%
}
\STATE{Calculate $z_k \in V$ by solving the problem 
\[
\text{Find }
z_k   \in V \text{ such that }
  (A(q_k)^\top z_k, w  )_{L^2(\Omega)^m} = 
-
 \left \langle \partial_v J(v_k,q_k), w \right \rangle_V
\quad \forall w \in V.
\]%
}
\STATE{Calculate $p_k \in L^1(0,L)$ via the formula
\[
p_k =
\int_D
\big [ A'(q_k ) v_k
\big ]^\top z_k
\d (x_2,\ldots,x_{n+1}).
\]
}
\STATE{Calculate $q_{k+1} \in \Qad$ by solving the 
variational inequality
\begin{equation*}
\begin{aligned}
&\text{Find } q_{k+1} \in \Qad \text{ such that }
\\
&\left ( q_{k+1}  - q_k, q - q_{k+1} \right )_Q
+
\left \langle \tau \partial_q J(v_k,q_k), q - q_{k+1} \right \rangle_Q
+
\int_0^L
\tau p_k (q - q_{k+1}) \d s
\geq 0
\quad \forall q \in \Qad.
\end{aligned}
\end{equation*}
}
\IF{$\|q_{k+1} - q_k\|_Q  \leq \tol$\label{alg:line:7}}
\STATE{{\bf \emph{break}} (Return $q_k$.)}
\ENDIF
\ENDFOR
\end{algorithmic}
\end{Algorithm}
\pagebreak

By means of \cref{lem:descent_property}, one obtains:

\begin{theorem}[{Convergence of \cref{alg:gradproj}}]
\label{th:convgradproj}
Let $\hat J\colon O \to \R$ and $\hat C > 0$ be as in 
\cref{lem:LipRedFunction}.
Suppose that $q_0 \in \Qad$ is arbitrary and that 
$\tau \in (0, 2/\hat C)$. Then, the following is true:
\begin{enumerate}[label=\roman*)]
\item\label{th:convgradproj:i}
If $\tol > 0$ holds, then 
the termination criterion in 
step \ref{alg:line:7} of \cref{alg:gradproj}
is triggered after finitely many iterations
and the returned iterate $q_k$ satisfies 
\begin{equation}
\label{eq:approxstationary}
 \|
q_k - \proj_{\Qad}
 (
q_k - \tau \nabla_Q \hat J(q_k)
  )
\|_Q \leq \tol.
\end{equation}
If $J(v,q) \geq 0$ holds for all $v \in V$ and $q \in \Qad$,
then the index $k$ for which the termination criterion is 
triggered can be estimated by 
\begin{equation}
\label{eq:iteratebound}
0 \leq k \leq 
\frac{2 \tau \hat J(q_0)}{(2 - \hat C\tau )\tol^2}.
\end{equation}
\item\label{th:convgradproj:ii} 
If $\tol = 0$ holds and $\dim(Q) < \infty$, 
then \cref{alg:gradproj} either terminates after 
finitely many iterations with an iterate $q_k$ that is Bouligand stationary 
for \eqref{eq:P} in the sense of \eqref{eq:Bouligand}
or it produces a sequence of iterates $\{q_k\} \subset \Qad$
such that $\{q_k\}$ possesses at least one accumulation 
point and such that every accumulation point of 
$\{q_k\}$ is Bouligand stationary for \eqref{eq:P}.
\end{enumerate}
\end{theorem}
\begin{proof}
From \cref{lem:derredObjective} and the definition of $\proj_{\Qad}$, 
we obtain 
that $q_{k+1} =  \proj_{\Qad}
 (
q_k - \tau \nabla_Q \hat J(q_k)
  )$ holds in \cref{alg:gradproj}.
  In combination with \cref{lem:descent_property}
  and $\tau \in (0, 2/\hat C)$,
  this yields
  \begin{equation*}
  \label{eq:suff_decrease_again}
  \hat J(q_k)
 \geq 
\hat J(q_{k+1})
 +
\left (
\frac{1}{\tau} - \frac{\hat C}{2}
\right )
\left \|
q_k -  
q_{k+1} 
\right \|_Q^2
\geq \hat J(q_{k+1}) 
  \end{equation*}
  for all performed iterations $k$.
  Suppose now that $\tol > 0$ holds and that 
  the termination criterion in step
  \ref{alg:line:7} of the algorithm is never triggered. Then, 
  for all $k \in \N_0$, we have 
  \begin{equation}
  \label{eq:randomestimate6e}
  \hat J(q_0)
 \geq 
\hat J(q_{k})
 +
k \left (
\frac{1}{\tau} - \frac{\hat C}{2}
\right )
\tol^2,
  \end{equation}
  which implies
  $\hat J(q_{k}) \to - \infty$,
  in contradiction with \cref{th:solvabilityP}. 
  Thus, \cref{alg:gradproj} indeed has to terminate 
  finitely in \ref{th:convgradproj:i}.
  It
  follows immediately from the 
  identity 
  $q_{k+1} =  \proj_{\Qad}
 (
q_k - \tau \nabla_Q \hat J(q_k)
  )$ and the form of the termination criterion that the returned iterate $q_k$
  satisfies 
  \eqref{eq:approxstationary}. 
  If, finally, $J(v,q) \geq 0$ holds for all $v \in V$ and $q \in \Qad$
  in \ref{th:convgradproj:i},
  then it follows analogously to \eqref{eq:randomestimate6e}
  that 
    \begin{equation*}
  \hat J(q_0)
 \geq 
(k+1) \left (
\frac{1}{\tau} - \frac{\hat C}{2}
\right )
\tol^2
  \end{equation*}
  holds if the termination criterion 
  is not triggered 
  in the iterations $l=0,\ldots, k$.
  This estimate immediately 
  yields the bound \eqref{eq:iteratebound} 
  on the index $k$ that the criterion is triggered for.

  Let us now assume that 
  $\tol = 0$ and $\dim(Q) < \infty$.
  If, in this case, 
  \cref{alg:gradproj} terminates 
  finitely, then the returned iterate 
  $q_k$ satisfies 
  $\|q_{k+1} - q_k\|_Q = \|\proj_{\Qad}
 (
q_k - \tau \nabla_Q \hat J(q_k)
  ) - q_k\|_Q = 0$
  and is Bouligand stationary by 
  \cref{lem:grad_opt_reform}.
  If \cref{alg:gradproj}
  does not terminate finitely in this setting, then 
  \cref{lem:descent_property} and 
  \cref{th:solvabilityP} imply that 
  \begin{equation*}
  \hat J(q_0)
 \geq 
\hat J(q_{k})
 +
 \sum_{l=0}^{k-1}
\left (
\frac{1}{\tau} - \frac{\hat C}{2}
\right )
\|\proj_{\Qad}(q_{l} - \tau \nabla_Q \hat J(q_l)) - q_l\|_Q^2
\geq
\hat J(q_{k})
\geq
\min_{q \in \Qad} \hat J(q) > -\infty
\quad \forall k \in \N_0.
  \end{equation*}
  In view of the radial unboundedness property 
  in \cref{ass:standing:abstract}\ref{ass:standing:abstract:vii},
  this yields that the sequence 
  $\{q_k\}$ is bounded and that 
  \[
   \sum_{l=0}^{\infty} \|\proj_{\Qad}(q_{l} - \tau \nabla_Q \hat J(q_l)) - q_l\|_Q^2 < \infty.
  \]
  In particular, we have
  $\|\proj_{\Qad}(q_{k} - \tau \nabla_Q \hat J(q_k)) - q_k\|_Q^2 \to 0$ for 
  $k \to \infty$. The assertions 
  in \ref{th:convgradproj:ii}
  now follow immediately from the theorem of 
  Bolzano--Weierstrass, 
  \cref{lem:LipRedFunction}
  and the fact that the projection $\proj_{\Qad}$ is 
  one-Lipschitz (and hence continuous). 
\end{proof}

\begin{remark}~
\begin{itemize}
\item 
It is of little use to consider a 
backtracking
line-search (e.g.\ of Armijo type)
for the function 
$(0,\infty) \ni \tau \mapsto \hat J(q)
-
\hat J(\proj_{\Qad}
 (q - \tau \nabla_Q \hat J(q))) \in \R$
in \cref{alg:gradproj}; cf.\
\cite[Section 9.4.2]{Beck2014}.
Indeed,
such a backtracking procedure 
requires numerous evaluations of the 
projection $\proj_{\Qad}$,
which is prohibitively costly 
in the situation of \eqref{eq:P}
owing to the required properties of the space $Q$. 
As we will see in the remainder of the paper, 
the main advantage of our semismooth 
Newton approach is that it allows to reduce 
the number of evaluations of 
$\proj_{\Qad}$ required for the solution of \eqref{eq:P} significantly 
in comparison with a gradient-descent method. 
\item 
\Cref{lem:grad_opt_reform} shows that the left-hand side 
of \eqref{eq:approxstationary} measures 
the degree 
of stationarity of an iterate $q_k$. 
Case \ref{th:convgradproj:i} of \cref{th:convgradproj}
thus provides us with finite convergence 
to approximate 
stationary points in function space.
As the right-hand side of 
\eqref{eq:iteratebound}
is stable w.r.t.\ discretisations, 
the estimate \eqref{eq:iteratebound} 
establishes a mesh-independence of this convergence property.
\end{itemize}
\end{remark}

We conclude this subsection by proving that, 
for a tracking-type objective, 
we can also obtain strong convergence of the 
iterates $\{q_k\}$
in the infinite-dimensional setting,
provided the involved Tikhonov parameter is large enough. 

\begin{corollary}
[{\Cref{alg:gradproj} for tracking-type problems}]
\label{cor:global_convergence_tracking_type}
Suppose (in addition to our standing assumptions) 
that $J$
has the form 
\[
J(v,q) = J_1(v) + \frac{\nu}{2}\|q\|_Q^2
\]
with a function $J_1\colon V \to \R$ 
and a Tikhonov parameter $\nu \in \R$ satisfying 
\[
\nu > 
\bar C :=
\Bigg (
 \frac{C_1 C_5 C_7^2 \|F\|_{V^*}}{\alpha^2} 
+ 
\frac{C_2 C_4 C_6 C_7^2\|F\|_{V^*}^2}{\alpha^4} 
+ \frac{C_2 C_6 C_7 \|F\|_{V^*} }{\alpha^2}
+
2\frac{C_1C_4 C_6 C_7^2  \|F\|_{V^*} }{\alpha^3}
\Bigg ).
\]
Then, \eqref{eq:P}
possesses precisely one solution
$\bar q \in \Qad$, this solution is the only 
 Bouligand stationary point of \eqref{eq:P}
 and, for all $q_0 \in \Qad$ and all $0 < \tau < 2/(\bar C + \nu)$,
\cref{alg:gradproj} with $\tol = 0$ 
either terminates finitely with $\bar q$
or produces iterates
$\{q_k\}$ that converge linearly in $Q$ to $\bar q$.
\end{corollary}

\begin{proof}
Let $q_1, q_2 \in \Qad$ be given,
define  
$v_i := S(q_i)$, $i=1,2$,
and let $p_1, p_2 \in L^1(0,L)$
and $z_1, z_2 \in V$ denote the quantities 
associated with $q_1, q_2$ in 
\cref{lem:derredObjective}.
Then, we obtain from the estimates 
in the proof of \cref{lem:LipRedFunction} that 
\begin{equation*}
\begin{aligned}
&C_7\|p_1 - p_2\|_{L^1(0,L)}
\\
&\leq  
\Bigg (
 \frac{C_1 C_5 C_7^2 \|F\|_{V^*}}{\alpha^2} 
+ 
\frac{C_2 C_4 C_6 C_7^2\|F\|_{V^*}^2}{\alpha^4} 
+ \frac{C_2 C_6 C_7 \|F\|_{V^*} }{\alpha^2}
+
2\frac{C_1C_4 C_6 C_7^2  \|F\|_{V^*} }{\alpha^3}
\Bigg )
\|q_1 - q_2\|_Q
\\
&= \bar C \|q_1 - q_2\|_Q. 
  \end{aligned}
\end{equation*}
If we define 
$T_\tau\colon \Qad \to \Qad$, 
$T_\tau(q) := \proj_{\Qad}(q - \tau \nabla_Q \hat J(q))$,
for $\tau > 0$, then the above, \cref{lem:derredObjective} and the 
fact that $\proj_{\Qad}$ is one-Lipschitz imply that 
\begin{equation*}
\begin{aligned}
\|T_\tau(q_1) - T_\tau(q_2)\|_{Q}
&=
\|\proj_{\Qad}(q_1 - \tau \nabla_Q \hat J(q_1)) - 
\proj_{\Qad}(q_2 - \tau \nabla_Q \hat J(q_2))\|_{Q}
\\
&\leq
\| q_1 - q_2 - \tau \nabla_Q \hat J(q_1) + \tau \nabla_Q \hat J(q_2) \|_{Q}
\\
&=
\| q_1 - q_2 - \tau \nu (q_1 - q_2) - \tau \iota_R(p_1 - p_2) \|_{Q}
\\
&\leq
|1 - \tau \nu|
\| q_1 - q_2\|_{Q}
+
\tau \bar C \|q_1 - q_2\|_Q.
\end{aligned}
\end{equation*}
For $\nu > \bar C$ and $0 < \tau < 2/(\bar C + \nu)$,
one has 
$|1 - \tau \nu| + \tau \bar C \in [0,1)$ and
$T_\tau$ is a contraction.
In view of the Banach fixed-point theorem,
this implies that, for all such $\nu$ and $\tau$, 
there exists precisely 
one $\bar q \in \Qad$
satisfying 
$\bar q = T_\tau(\bar q)$
and that the sequence generated by 
$q_{k+1} := T_\tau(q_k)$
converges linearly to $\bar q$ for all $q_0 \in \Qad$.
As the equation $\bar q = T_\tau(\bar q)$
is equivalent to the Bouligand stationarity condition
\eqref{eq:Bouligand} by \cref{lem:grad_opt_reform},
since \eqref{eq:P} possesses a global solution
which is necessarily Bouligand stationary 
by \cref{th:solvabilityP,th:first_order_necessary_condition},
and since the iterates 
$\{q_k\}$ of \cref{alg:gradproj} satisfy 
$q_{k+1} =  \proj_{\Qad}(q_k - \tau \nabla_Q \hat J(q_k)) = T_\tau(q_k)$,
the assertions now follow 
immediately. 
\end{proof}

\subsection{A semismooth Newton method}
\label{subsec:semiNewton}

As already mentioned, a main drawback of the 
gradient-projection method discussed in \cref{subsec:gradproj}
is that it typically requires dozens 
of evaluations 
of the projection $\proj_{\Qad}$ 
to identify a stationary point of \eqref{eq:P} up to 
a reasonable tolerance; see the numerical 
experiments in \cref{sec:num_exp}. 
Owing to the requirements on $Q$
in \cref{ass:standing:abstract}
(that are necessary to ensure the existence 
of a minimiser in view of the nonlinear nature of the term 
$\left (A(q)v, w \right )_{L^2(\Omega)^m}$
in \eqref{eq:var_id}), 
each of these evaluations amounts 
to the solution of a non-trivial quadratic optimisation 
problem with constraints and, thus, 
entails significant computational cost. 
To alleviate this problem, in what follows,  
we develop a semismooth Newton method
for the numerical solution of 
\eqref{eq:P}. Since Newton algorithms 
typically require far less iterations 
than gradient-descent methods to achieve a 
certain tolerance and since 
both a Newton- and a gradient-projection step 
require precisely 
one projection operation (see \cref{alg:gradproj} and \cref{th:convergence_Newton} below),
a Newton method has the potential to 
reduce the numerical cost induced by the projection 
evaluations significantly when considering problems of the type \eqref{eq:P}.
(The smaller number of projection operations, of course, then comes 
at the price of the additional overhead 
introduced by the Newton system.) 
For the sake of simplicity, in what follows, 
we develop our semismooth Newton approach only for tracking-type objectives 
of the form $J(v,q)=J_1(v)+(\nu/2) \|q\|_Q^2$
as discussed in \cref{rem:highreg}
and \cref{cor:global_convergence_tracking_type}.
More precisely, throughout this subsection,
we consider the following setting: 

\begin{assumption}[{Standing assumptions for \cref{subsec:semiNewton}}]
\label{ass:standing:Newton}%
In addition to the conditions in \cref{ass:standing:abstract},
we assume that:
\begin{enumerate}[label=\roman*)]
\item
$Q$ is a Hilbert space with Riesz isomorphism 
$\iota_R\colon Q^* \to Q$ 
as in
\cref{ass:standing:appendix_A}\ref{ass:standing:appendix_A:ii};

\item
\label{ass:standing:Newton:ii}
$a,b \in C([0,L])$ are functions 
satisfying the conditions in \cref{ass:standing:appendix_A}\ref{ass:standing:appendix_A:iii};

\item\label{ass:standing:Newton:iii} $A$ is 
twice continuously differentiable with 
second derivative $A''(s) \in \R^{m \times m}$ for all $s \in \R$;
\item\label{ass:standing:Newton:iv}
$J$ has the form 
$J(v,q)=J_1(v)+ \frac\nu{2} \|q\|_Q^2$
with a twice continuously differentiable
$J_1\colon V \to \R$.
\end{enumerate}
\end{assumption}

It is easy to check that the 
conditions on $Q$ in \cref{ass:standing:Newton} 
(respectively, \cref{ass:standing:appendix_A})
are, for example, 
satisfied by the spaces 
$H^s(0,L)$ with $s \in (1/2, 1]$. 
The point of departure for the development of our semismooth
Newton method is the observation that,
for the situation in 
\cref{ass:standing:Newton}
and $\tau = 1/\nu$, 
the necessary optimality condition 
in \cref{lem:grad_opt_reform} reads 
\begin{equation}
\label{eq:random_eq_2636_necop}
0 =
\bar q - \proj_{\Qad}
\left (
\bar q - \tau \nabla_Q \hat J(\bar q)
\right )
=
\bar q - \proj_{\Qad}
 \left (
 - \frac1\nu \iota_R\left ( \bar p \right )
\right  ).
\end{equation}
Here, the function
$\bar p \in L^1(0,L) \subset Q^*$
is defined as in \eqref{eq:KKT}, that is,
$\bar p := S'(\bar q)^* J_1'(S(\bar q))$,
where $S'(\bar q)^* \in \LLLL(V^*,Q^*)$
denotes the adjoint of $S'(\bar q) \in \LLLL(Q,V)$.
If we define $\Sigma :=  \proj_{\Qad} \circ \iota_R \colon Q^* \to Q$, i.e.\ write $\Sigma$ for the solution mapping 
$Q^* \ni g \mapsto \hat q \in Q$ of the variational 
inequality 
\begin{equation}
\label{eq:VI_again_ss}
\text{Find }
\hat q \in \Qad \text{ such that }
(\hat q, q-\hat q)_{Q}
- \left \langle g, q - \hat q\right \rangle_Q\geq 0
\quad\, \forall q \in \Qad,
\end{equation}
and denote by $P$ the map 
$P\colon O \to Q^*$, 
$ P(q) := S'(q)^* J_1'(S(q))$,
defined on the set $O$ from \cref{lem:Sfrechet}
and 
given in explicit form by 
the formulas \eqref{eq:pdef}
and \eqref{eq:adjointequation}
with $v := S(q)$, then 
\eqref{eq:random_eq_2636_necop} can 
be recast as
\begin{equation}\label{eq:newton_Equation}
\bar q-\Sigma\left ( -\frac{1}{\nu } P(\bar q)\right)=0.
\end{equation}
In what follows, our goal will be to show that 
the operator equation \eqref{eq:newton_Equation}
is amenable to a Newton-type algorithm.
As the mapping $\Sigma$ is 
typically nondifferentiable, 
for this, 
we have to show that the function
$Q \supset O \ni q \mapsto  q-\Sigma\left ( -\frac{1}{\nu } P(q)\right ) \in Q$
is \emph{Newton differentiable};
see \eqref{eq:NewtonDiffDef}. 
(Note that, alternatively, one could also work 
with the system 
\eqref{eq:KKT}
in 
$v$, $z$, $q$, and $p$ here instead of the 
reduced equation \eqref{eq:newton_Equation}.
We focus on the reduced approach in what follows.)
As a first step towards 
establishing the Newton differentiability 
of the operator 
$O \ni q \mapsto  q-\Sigma\left ( -\frac{1}{\nu } P(q)\right ) \in Q$,
we prove two lemmas 
which prepare the
 analysis of 
the smoothness properties 
of the function 
$P$.

\begin{lemma}[Twice continuous differentiability of the 
operator in the governing equation]
\label{lem:aux_dif2_op}
The operator $\FF$ in \eqref{eq:FFdef} is twice continuously
differentiable with the following second partial derivatives
$\partial_{vv} \FF\colon 
V \times Q \to \BBBB(V,V;V^*)$,
$\partial_{vq} \FF\colon 
V \times Q \to \BBBB(V,Q;V^*)$, and
$\partial_{qq} \FF\colon 
V \times Q \to \BBBB(Q,Q;V^*)$:
\begin{align*}
&
\left\langle\partial_{vv} \FF(v,q)[\delta v_1,\delta v_2] , w\right \rangle_{V}
=0
&&\forall w \in V\quad \forall\delta v_1,\delta v_2 \in V,
\\
&
\left\langle\partial_{vq} \FF(v,q)[\delta v,\delta q] , w\right \rangle_{V}
=\left (A'(q) \delta q\, \delta v , w \right )_{L^2(\Omega)^m}
&&\forall w \in V\quad \forall\delta v\in V,\delta q \in Q,\\
&
\left \langle  
\partial_{qq} \FF(v,q)[\delta q_1,\delta q_2] , w\right \rangle_{V}
=\left (A''(q) \delta q_1\delta q_2\, v , w \right )_{L^2(\Omega)^m}
&&\forall w \in V\quad \forall \delta q_1,\delta q_2 \in Q.
\end{align*}
\end{lemma}
\begin{proof}
We extend the proof of Lemma \ref{lem:aux_dif_op} and use the notation there.
Only $\FF_2$ needs attention. Its derivative 
is given by
\eqref{eq:F2deriv}. We have to show that
$\FF_4\colon L^\infty(\Omega) \to L^\infty(\Omega)^{m\times m}$,
$\phi\mapsto A'(\phi)$,
is continuously differentiable. By assumption,
$A'\colon \R\to\R^{m\times m}$ is $C^1$,
which implies that 
$A''$ is 
uniformly continuous on bounded sets. We thus can apply
\cite[Theorem 7 and Remark 4]{GKT92} to obtain that
$\FF_4$ is continuously differentiable with
\[
(\FF_4'(\phi)\delta \phi)(x)=A''(\phi(x))\delta \phi(x)\quad
\text{ for a.a. }x\in \Omega
\quad \forall \phi,\delta \phi\in L^\infty(\Omega).
\]
This shows that
$\FF_2$ is twice continuously differentiable with
\[
(\FF_2''(\phi)[\delta \phi_1,\delta \phi_2])(x)=(\FF_4'(\phi)\delta \phi_1)(x)\delta \phi_2(x)
=A''(\phi(x))\delta \phi_1(x)\delta \phi_2(x),\quad 
\forall \phi,
\delta \phi_1,\delta \phi_2\in L^\infty(\Omega),
\]
for a.a.\ $x\in \Omega$.
The formulas for the derivatives now follow easily with the chain rule.
\end{proof}

Next, we use the implicit function theorem to show that the solution
operator $S$ is twice continuously differentiable.
\begin{lemma}[Twice continuous differentiability of $S$]
\label{lem:Sfrechet2}
There exists 
an  open neighbourhood $O \subset Q$ of $\Qad$
such that 
the solution mapping 
$S\colon \Qad \to V$, $q \mapsto v$,
of \eqref{eq:var_id}
can be extended  to a twice continuously 
%Fréchet
differentiable function $S\colon O \to V$ such that  the assertions of \cref{lem:Sfrechet} hold and
such that, 
for every $q \in O$,
the second derivative $\delta_2 v := S''(q)[\delta q_1,\delta q_2] \in V$
of $S$ at $q$
in the directions
$\delta q_1,\delta q_2 \in Q$ is uniquely characterised by the 
problem 
\begin{equation}
\label{eq.DirDif2VE}
\begin{aligned}
&\text{Find }\delta_2 v\in V \text{ such that }\\
&\left (A(q)\delta_2 v, w \right )_{L^2(\Omega)^m} = 
-\left (A'(q)\delta q_1 \delta v_2+
A'(q)\delta q_2 \delta v_1+A''(q)\delta q_1\delta q_2 S(q), w \right )_{L^2(\Omega)^m}
\hspace*{\fill} ~~\forall w \in V,
\end{aligned}
\end{equation}
where $\delta v_i : =S'(q)\delta q_i$ 
is defined by 
\eqref{eq.DirDifVE}
with $\delta q =\delta q_i$ and
$\delta v=\delta v_i$, $i=1,2$.
\end{lemma}
\begin{proof}
%(A(q)S(q)-F,w)=0
%(A(q)S'(q)dq1+A'(q)dq1 S(q),w)=0
%(A(q)S''(q)[dq1,dq2]+A'(q)dq2 S'(q)dq1+A'(q)dq1 S'(q)dq2+A''(q)dq1dq2S(q),w)=0
We just have to extend the proof of \cref{lem:Sfrechet}.
Beyond the assumptions and arguments there, twice continuous
differentiability of the implicit function $S$ only requires that
$\FF$ is twice continuously
differentiable; see
\cite[Theorem 2.3]{Wachsmuth_2014_2}. This, however, was established in \cref{lem:aux_dif2_op}.
The formula for the second derivative of the implicit function
is well known. It results from
differentiating the equation $\FF(S(q),q)=0$ twice w.r.t.\ $q$
and solving for $S''(q)[\delta q_1,\delta q_2]$:
\begin{align*}
\partial_{v} \FF(S(q),q)S''(q)[\delta q_1,\delta q_2]+
\partial_{vq} \FF(S(q),q)[S'(q)\delta q_1,\delta q_2]&+
\partial_{qv}\FF(S(q),q)[\delta q_1,S'(q)\delta q_2]\\
&+
\partial_{qq} \FF(S(q),q)[\delta q_1,\delta q_2]=0;
\end{align*}
see \cite[Equation (5)]{Wachsmuth_2014_2}.
\end{proof}

In combination with the smoothness
properties of $J_1$
in \cref{ass:standing:Newton},
\cref{lem:Sfrechet2} yields that 
$Q \supset O \ni q\mapsto J_1(S(q)) \in \R$
is a $C^2$-map. This, in turn, 
implies that the function
$P\colon O \to Q^*$, 
$ P(q) := S'(q)^* J_1'(S(q))$,
is $C^1$ on the open neighbourhood
$O$ of $\Qad$ it is defined on. 
In what follows, we
show that 
$P$ is even $C^1$ when understood as a function 
with values in 
$L^1(0,L)$. 
This will be essential for the 
analysis of the Newton differentiability of 
the map in the operator equation 
\eqref{eq:newton_Equation}. 
To obtain the higher regularity of $P$,
we use the formulas 
in \cref{lem:derredObjective}. 

\begin{lemma}\label{lem_p_lemma_smooth}
There exists an open neighbourhood $O\subset Q$ of $\Qad$ such that the 
function $P$ given by 
$ P(q) := S'(q)^* J_1'(S(q))$
is well defined and continuously 
differentiable as a map 
$P\colon O \to L^1(0,L)$.
For all $q\in O$, $P$ satisfies 
$P(q) = p \in L^1(0,L) \subset Q^*$, 
where $p\in L^1(0,L)$ is characterised by \eqref{eq:pdef}
with $z=z(v,q)\in V$ being the 
solution of the adjoint equation \eqref{eq:adjointequation} with $v=S(q)$.
\end{lemma}

\begin{proof}
Let $O\subset Q$ be the open neighbourhood from \cref{lem:Sfrechet}.
It
follows immediately 
from \cref{lem:Sfrechet,lem:derredObjective} that $P$ is well defined as a function 
from $O$ to $L^1(0,L)$ and can be expressed 
in terms of the equations 
\eqref{eq:pdef} and \eqref{eq:adjointequation}.
It remains to prove that 
$P\colon O \to L^1(0,L)$ is continuously 
differentiable. To this end,
we use the implicit function
theorem to show that the solution map
$V\times O \ni (v,q) \mapsto z\in V$ 
of \eqref{eq:adjointequation}
is continuously differentiable and then apply the chain rule to show that $P\colon O \to L^1(0,L)$ is continuously differentiable:
With arguments very similar to the proof of \cref{lem:aux_dif_op}, one checks that the
map $\GG:V\times O\times V\to V^*$ defined by 
\[
\langle\GG(v,q,z),w\rangle_V
:=(A(q)^\top z,w)_{L^2(\Omega)^m}+\langle\partial_v J(v,q),w\rangle_V
\quad\forall w \in V
\]
is continuously differentiable with partial derivatives
\begin{align*}
&\langle\partial_v \GG(v,q,z)\delta v,w\rangle_V
=\partial_{vv} J(v,q)[w,\delta v],\\
&\langle\partial_q \GG(v,q,z)\delta q,w\rangle_V
=
(A'(q)^\top \delta q \,z,w)_{L^2(\Omega)^m}+
\partial_{vq} J(v,q)[w,\delta q],\\
&\langle\partial_z \GG(v,q,z)\delta z,w\rangle_V
=(A(q)^\top \delta z,w)_{L^2(\Omega)^m}.
\end{align*}
Using parts of the proof of \cref{lem:derredObjective}
and the notation therein,
we obtain that 
$\partial_z \GG(v,q,z)=\partial_v\FF(v,q)^*\in\LLLL(V,V^*)$
is an isomorphism for all $q\in O$. The implicit function theorem
now 
implies the existence of a unique
continuously differentiable 
solution map $Z: V\times O\to V$,
$(v,q) \mapsto z$,
of the adjoint
equation
\[
\langle\GG(v,q,z),w\rangle_V=0
\quad\forall w \in V.
\]
The derivatives of $Z$ are given by
\[
\partial_v Z(v,q)=-\partial_z \GG(v,q,z)^{-1}\partial_v \GG(v,q,z),\quad
\partial_q Z(v,q)=-\partial_z \GG(v,q,z)^{-1}\partial_q \GG(v,q,z).
\]
The trilinear bounded mapping $\FF_5: V\times L^\infty(\Omega)^{m\times m}\times V\to L^1(0,L)$
defined by 
\[
\FF_5(v,B,z)(x_1) : =
\int_D
\big [ B(x) v(x)
\big ]^\top z(x)
\d (x_2,\ldots,x_{n+1})
\]
is infinitely continuously differentiable. By \cref{lem:Sfrechet}, $S$ is continuously differentiable and, as shown in the proof of \cref{lem:aux_dif2_op}, 
the function $\FF_4\colon L^\infty(\Omega) \to L^\infty(\Omega)^{m\times m}$,
$\phi\mapsto A'(\phi)$,
is continuously differentiable too. Thus, we can apply the chain rule to conclude that
\[
P: O\to L^1(0,L)
,
\qquad
P(q)=\FF_5(S(q),\FF_4(q),Z(S(q),q)),
\]
is continuously differentiable,
where we used the
embeddings $Q \compactly L^\infty(0,L)\hookrightarrow L^\infty(\Omega)$.
\end{proof}

Since it is now proved that $P$ is $C^1$ from an open neighbourhood $O$ of $\Qad$ to $L^1(0,L)$, the mapping $\Sigma$ is the only part of
the operator equation \eqref{eq:newton_Equation} that 
remains to be analysed. 
As this operator is precisely the 
nonsmooth part of \eqref{eq:newton_Equation}, 
its treatment is the most delicate issue when setting up a Newton algorithm. Only recently, it has been demonstrated 
in a series of papers \cite{ChristofWachsmuthUniObstacle,ChristofWachsmuthBiObstacle,Rauls2020}
how the Newton differentiability of the solution operator 
$\Sigma\colon Q^* \to Q$, $g \mapsto \hat q$, 
of an obstacle-type variational inequality 
of the form \eqref{eq:VI_again_ss} can be established. 
The following \cref{th_main_semismooth_main_body} 
summarises the result that 
is obtained by proceeding along the lines of \cite{ChristofWachsmuthUniObstacle,ChristofWachsmuthBiObstacle,Rauls2020} in the situation of \cref{ass:standing:Newton}.
As its proof requires some modifications of the arguments in \cite{ChristofWachsmuthUniObstacle,ChristofWachsmuthBiObstacle,Rauls2020},
to spare the reader from having to go through the analysis of
the original sources in detail, and to facilitate its future use, 
we prove \cref{th_main_semismooth_main_body}  in a largely self-contained 
form in \cref{sec:appendix}.

\begin{theorem}[Newton Differentiability of $\Sigma$]
\label{th_main_semismooth_main_body}
The function $\Sigma :=  \proj_{\Qad} \circ \iota_R \colon Q^* \to Q$
is Newton differentiable 
as a map $\Sigma\colon \MM([0,L]) \to Q$
when endowed with the Newton 
derivative 
$ \D \Sigma \colon \MM([0,L]) \to \LLLL(\MM([0,L]), Q)$
defined by 
\begin{equation}
\label{eq:DSigma_def_2}
	\delta = D\Sigma(g)h
	\qquad
	\iff
	\qquad
	\delta \in \QQ(g),\quad 
	\left (\delta, q \right )_Q = \left \langle h, q\right \rangle_{Q}
	\quad 
	\forall q \in \QQ(g)
\end{equation}
with $\QQ(g) := \{q \in Q \mid q(x) = 0 
\text{ for all } x \in [0,L] \text{ with }  \Sigma(g)(x) = a(x) \text{ or } \Sigma(g)(x) = b(x)\}$.
\end{theorem}

\begin{proof}
See \cref{th_main_appendix} and the preceding results in the appendix.
\end{proof}

With \cref{th_main_semismooth_main_body} at hand,
we now arrive at:

\begin{corollary}
\label{cor:chain}
Let $O$, $P$, $\Sigma$ and 
$D\Sigma$ be as in \cref{lem_p_lemma_smooth} and \cref{th_main_semismooth_main_body}. Define 
\begin{equation}
\label{eq:residue_def}
    \RR\colon O \to Q,\qquad 
    \RR(q) := q-\Sigma\left (- \frac1\nu P(q) \right ).
\end{equation}
Then, $\RR$ is Newton differentiable with Newton derivative
\begin{equation}
\label{eq:Dresidue_def}
    D \RR\colon O \to \LLLL(Q,Q),
    \qquad
    D \RR(q) := \Id +\frac{1}{\nu}D\Sigma\left( - \frac1\nu P(q)\right)P'(q).
\end{equation}
\end{corollary}

\begin{proof}
Since $P\colon O \to L^1(0,L)$ is continuously differentiable, 
$\Sigma\colon \MM([0,L])  \to Q$ is Newton differentiable,
$L^1(0,L) \hookrightarrow \MM([0,L])$,
and the Newton derivatives $D\Sigma(q)$ are bounded, 
the assertion follows immediately from the chain rule for 
Newton derivatives \cite[Lemma A.1]{Alphonse2026}.
\end{proof}

Note that it is essential in the proof of \cref{cor:chain} 
to have the Newton differentiability of $\Sigma$ on $L^1(0,L)$
as $P$ cannot be expected to map into $L^s(0,L)$ for some
$s>1$ in general; cf.\ \eqref{eq:pdef}. Precisely the $L^1$-case 
is excluded in \cite{ChristofWachsmuthUniObstacle,ChristofWachsmuthBiObstacle,Rauls2020}
(where also only the space $H_0^1(\Omega)$ is considered). This is the main reason why we provide a detailed proof of \cref{th_main_semismooth_main_body}
in \cref{sec:appendix}. The standard convergence 
analysis of the semismooth Newton method now 
provides us with:

\begin{theorem}[Superlinear convergence of the semismooth Newton method]
\label{th:convergence_Newton}
Let $O$, $P$, $\Sigma$ and 
$D\Sigma$ be as in \cref{lem_p_lemma_smooth} and \cref{th_main_semismooth_main_body}, and let 
$\RR\colon O \to Q$ and $D \RR\colon O \to \LLLL(Q,Q)$
be defined as in \eqref {eq:residue_def} and 
\eqref{eq:Dresidue_def}, respectively. 
Suppose that $\bar q \in \Qad$ is 
a stationary point of \eqref{eq:P} in the sense 
of \cref{th:first_order_necessary_condition} 
and that there exist $r, C >0$ such that 
the Newton derivative $D \RR (q)$ is continuously 
invertible 
with $\|D \RR (q)^{-1}\|_{\LLLL(Q,Q)} < C$
for all $q \in O$ satisfying $\|q - \bar q\|_Q < r$. 
Then, there exists 
$\rho > 0$ such that, for all 
$q_0 \in Q$ satisfying $\|q_0 - \bar q\|_Q < \rho$,
the semismooth Newton steps
\begin{equation}
\label{eq:Newton_step_def}
    q_{k+1} := q_k - D \RR (q_k)^{-1} \RR(q_k) 
\end{equation}
can be executed for all $k \in \N_0$ and 
produce a sequence of iterates
$\{q_k\}$ which converges $q$-superlinearly 
to $\bar q$ in $Q$ for $k \to \infty$.
\end{theorem}

\begin{proof}
In view of the Newton differentiability of $\RR$ established in 
\cref{cor:chain}, the assertion of this theorem follows 
immediately from standard results; see, e.g., \cite[Theorem 2.1]{Alphonse2026}. 
\end{proof}

As the Newton derivative 
\begin{equation*}  
D \RR\colon O \to \LLLL(Q,Q),
    \qquad
    D \RR(q) := \Id +\frac{1}{\nu}D\Sigma\left( - \frac1\nu P(q)\right)P'(q),
\end{equation*}
is a perturbation of the identity with the 
perturbation term scaled by $1/\nu$, 
one can expect that the 
uniform continuous invertibility of the operators 
$D \RR(q) \in \LLLL(Q,Q)$ is ensured for all $q \in O$
if the Tikhonov parameter $\nu > 0$ is large enough. 
This is indeed the case as the following result shows, 
which also establishes a direct link to the analysis of the 
gradient-projection method in \cref{subsec:gradproj}.

\begin{corollary}
\label{cor:invert_Newton_large_Tikhonov}
Suppose that
the quantities in \eqref{eq:P} satisfy not only 
the conditions in 
\cref{ass:standing:Newton} but also 
those
in  \cref{ass:standing_grad_proj} with constants 
$C_i > 0$, $i=1,\ldots, 6$. Let further
$C_7 > 0$ be defined as in \eqref{eq:C7_def}
and assume that the Tikhonov parameter $\nu$ satisfies 
\[
\nu > 
\bar C :=
\Bigg (
 \frac{C_1 C_5 C_7^2 \|F\|_{V^*}}{\alpha^2} 
+ 
\frac{C_2 C_4 C_6 C_7^2\|F\|_{V^*}^2}{\alpha^4} 
+ \frac{C_2 C_6 C_7 \|F\|_{V^*} }{\alpha^2}
+
2\frac{C_1C_4 C_6 C_7^2  \|F\|_{V^*} }{\alpha^3}
\Bigg ).
\]
Then, \eqref{eq:P}
possesses precisely one solution
$\bar q \in \Qad$, this solution is the only 
 Bouligand stationary point of \eqref{eq:P}
 and  the Newton derivatives $D\RR(q)$ in 
\cref{th:convergence_Newton} are invertible for all 
$q \in \Qad$ and satisfy
\begin{equation}
\label{eq:inverse_bound_newton}
\|D\RR(q)^{-1}\|_{\LLLL(Q,Q)} \leq
\frac{\nu}{\nu - \bar C}
\qquad \forall q \in \Qad.
\end{equation}
\end{corollary}

\begin{proof}
It
follows immediately from 
\cref{cor:global_convergence_tracking_type} that \eqref{eq:P} possesses precisely 
one solution, respectively, Bouligand stationary point.
It remains to prove the invertibility of the 
Newton derivatives $D\RR(q)$ 
for all $q \in \Qad$ with the bound in 
\eqref{eq:inverse_bound_newton}. To this end,
we note that the estimates in the proof of 
\cref{cor:global_convergence_tracking_type} and 
the definition of the map $P$ in 
\cref{lem_p_lemma_smooth} yield
\begin{equation*}
 C_7\|P(q_1) - P(q_2)\|_{L^1(0,L)}
 \leq  \bar C \|q_1 - q_2\|_Q \qquad \forall q_1, q_2 \in \Qad. 
\end{equation*}
In combination with the $C^1$-regularity of $P$
on $O \supset \Qad$ in \cref{lem_p_lemma_smooth}, 
the Lipschitz continuity property above provides us with 
the estimate 
$C_7\|P'(q)\|_{\LLLL(Q,L^1(0,L))} \leq \bar C$ for all $q\in \Qad$. 
From the definition of $D\Sigma$ in 
\eqref{eq:DSigma_def_2}, we further get that 
$D \Sigma (g) =\proj_{\QQ(g)} \circ  \iota_R$ holds for all $g \in L^1(0,L) \subset Q^*$, where 
$\proj_{\QQ(g)}$ denotes the orthogonal projection 
in $Q$ onto the subspace $\QQ(g)$ in \cref{th_main_semismooth_main_body}. 
As projections are non-expansive and since 
$|\langle p,q\rangle_Q | \leq C_7 \|q\|_Q \|p\|_{L^1(0,L)}$
for all $q \in Q$ and $p \in L^1(0,L)$, this 
allows us to conclude that 
\begin{equation}
\label{eq:randomeq_ujene}
\begin{aligned}
\left \|
\frac{1}{\nu}D\Sigma\left( - \frac1\nu P(q)\right)P'(q)
\right \|_{\LLLL(Q,Q)}
 \leq
\frac{1}{\nu}
\left \|
\iota_R\right \|_{\LLLL(Q^*,Q)}
C_7
\|P'(q)\|_{\LLLL(Q, L^1(0,L))}
 \leq \frac{\bar C}{\nu}\quad \forall q \in \Qad.
\end{aligned}
\end{equation}
The assertion of the corollary now follows immediately 
by means of the Neumann series. 
\end{proof}
 
\begin{remark}~
\begin{itemize}
\item Owing to the $C^1$-regularity of $P$ in \cref{lem_p_lemma_smooth} and the fact 
that $\|D\Sigma(g)\|_{\LLLL(Q^*, Q)} \leq 1$
holds for all $g \in \MM([0,L])$, the
norm on the left-hand side of 
\eqref{eq:randomeq_ujene} is also bounded 
by a constant in $[0,1)$ for all $q$ in 
a suitably chosen open
neighbourhood of $\Qad$. This means that the uniform 
invertibility estimate \eqref{eq:inverse_bound_newton} 
(with a slightly worse constant on the right-hand side)
carries over to an open neighbourhood of $\Qad$ too, 
as required for the application of \cref{th:convergence_Newton}.
\item Note that,
in the situation of 
\cref{cor:invert_Newton_large_Tikhonov},
we have $1/\nu \leq 2/(\bar C + \nu)$.
This implies that the global convergence 
result for \cref{alg:gradproj} in 
\cref{cor:global_convergence_tracking_type}
applies to \eqref{eq:P} 
if the 
step size is chosen as $\tau = 1/\nu$. 
From the proof of \cref{cor:global_convergence_tracking_type}, 
we further obtain that, for this choice of $\tau$,
the gradient-projection steps can be written as
\[
q_{k+1} 
=  \proj_{\Qad}\left (q_k - \tau \nabla_Q \hat J(q_k) \right) 
=  \Sigma\left ( \iota_R^{-1} (q_k - \tau \nabla_Q \hat J(q_k) \right)  
=  \Sigma\left (  \frac{1}{\nu} J_1'(q_k)
\right ) =: T(q_k),
\]
with the right-hand side of this equation 
defining a contractive operator 
$T\colon \Qad \to \Qad$.
In total, this shows that, 
in the situation 
of \cref{cor:invert_Newton_large_Tikhonov}, 
\eqref{eq:newton_Equation} 
corresponds to a 
contractive fixed-point equation
and the gradient-projection steps 
in \cref{alg:gradproj} with $\tau = 1/\nu$
coincide with the steps of the 
fixed-point iteration associated with \eqref{eq:newton_Equation}.
Interpreting the gradient-projection steps
with $\tau = 1/\nu$
along these lines makes it possible 
to globalise the 
semismooth Newton iteration in \eqref{eq:Newton_step_def} 
by means of 
the methods
developed in \cite[Section 2]{Alphonse2026}
and, in doing so, to construct a globally 
q-superlinear convergent algorithm 
for the solution of \eqref{eq:P}
that is additionally also robust w.r.t.\ inexactness;
see \cite[Theorem 2.1]{Alphonse2026}. 
As numerical experiments show that 
vanilla semismooth Newton steps are often more cost-effective 
than globalised steps (see \cite[Section 4]{Alphonse2026}), we do not discuss
globalisation strategies in detail here.

\item 
We expect that the locally uniform 
continuous invertibility of the 
Newton derivatives required in 
\cref{th:convergence_Newton} can be guaranteed 
by means of a no-gap second-order 
sufficient optimality condition and 
a strict complementarity assumption;
cf.\ the analysis in $L^2$-spaces in \cite{Casas2024}.
This topic will be discussed in detail in future work. 
\end{itemize}
\end{remark}

\section{Numerical experiments}
\label{sec:num_exp}
We conclude this paper with numerical experiments 
which compare the performance and convergence behaviour of 
the projected gradient-descent method in 
\cref{alg:gradproj} and the 
semismooth Newton iteration in 
\cref{th:convergence_Newton}. 

As a test configuration, we consider the following instance of 
the model problem \eqref{eq:model_problem_elliptic} which we used to motivate 
our analysis in \cref{sec:intro}:
\begin{equation}
\label{eq:model_problem_elliptic_num_ex}
\left.
	\begin{aligned}
		\text{Minimise} 
		\quad &  
        \frac12 \|u - u_\mathrm{obs}\|_{L^2(\Omega)}^2 + 
        \frac{\nu}{2} \|q\|_{H^1(0,L)}^2 \\
        \text{w.r.t.}
        \quad &u \in H_0^1(\Omega), \quad  q \in H^1(0,L),\\
		\text{s.t.}
		\quad & 
        -\nabla\cdot (B(q)\nabla u)=f \text{ in } \Omega,
\quad
u=0 \text{ on } \partial \Omega,
  \\
		\text{and} \quad & a \leq q \leq b \text{ a.e.\ in } (0,L).
\end{aligned}
\quad \right \}
\end{equation}
The problem data is chosen as
\begin{gather*}
\Omega = (0,1)^2,
\qquad
L=1,
\qquad
u_\mathrm{obs}(x_1,x_2) 
= 
S(q_{\mathrm{exact}})(x_1, x_2) +
0.005\,( \cos(\pi x_1) + \cos(\pi x_2)),
\\
\nu >0,
\qquad 
f(x_1, x_2) 
=
10 \sin(2 \pi x_1) \sin(\pi x_2),
\qquad
a \equiv - 0.15,
\qquad
b \equiv 0.15,
\end{gather*}
with
$q_{\mathrm{exact}}(x_1)
:= 
\max(-0.15, \min(0.15, x_1- 0.3))$
and 
$S$ as the solution map of 
the governing PDE. 
The observed state $u_\mathrm{obs}$ is thus 
equal to the state 
associated with $q_{\mathrm{exact}}$ plus a
(for the sake of reproducibility deterministic)
perturbation simulating measurement errors. 
As the field $B\colon [a,b] \to \R^{2\times 2}$, we use
$B(s) := \widehat B(s + 0.25)$, with $\widehat B\colon (0, 0.5) \to \R^{2\times 2}$,
$r \mapsto \widehat B(r)$,
as the effective tensor that arises when
the homogenisation approach in \cref{sec:homogen}
is followed
and the function
\begin{equation*}
\BBB_r(y) 
:=
\big [
c_1 + (c_2 - c_1)
\mathds{1}_{\{\zeta \in Y \colon |\zeta - m_Y| \leq r\}}(y)
\big ] \Id\qquad \forall y \in Y = (0,1)^2\qquad \forall r \in (0, 0.5)
\end{equation*}
in \eqref{eq:coeff_again_42} is considered as a model 
for the micro-structure with $c_1 := 10$ and $c_2 = 10^{-6}$. 
Note that this corresponds to 
a micro-perforation with circular holes of radius $r \in (0, 0.5)$
(or, more precisely, a situation with circular micro-inclusions of radius $r$ with a negligible principal coefficient),
analogous to the scenarios described in \cref{sec:intro,sec:homogen}.
The shift by $0.25$ in the argument of $\smash{\widehat B}$ in the definition 
of $B$ accomplishes that the Tikhonov regularisation in \eqref{eq:model_problem_elliptic_num_ex}
introduces a bias towards the micro-perforation radius $0.25$.
This facilitates producing minimisers $\bar q$ for which both of the constraints 
$a \leq q$ and $q \leq b$ are active.
It is easy to check that the setting above satisfies 
all of the conditions in 
Assumptions
\ref{ass:standing:homo}, \ref{ass:standing:abstract},
\ref{ass:standing_grad_proj},
\ref{ass:standing:Newton},
and \ref{ass:standing:appendix_A}
(cf.\ \cref{ex:scalar}). 

\begin{figure}[ht]
\begin{center}
  \includegraphics[width=6.5cm]{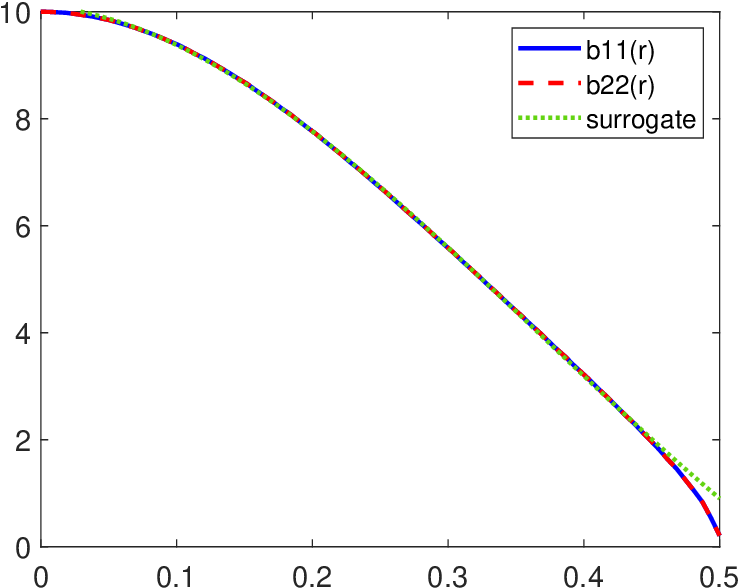}
  \caption{The components $\widehat B_{11}, \widehat B_{22}\colon (0, 0.5) \to \R$ of $\widehat B$ (calculated with a 
  finite element code) and their
  polynomial approximation from \eqref{eq:surrogate_polynomial}. 
  The components $\widehat B_{12}, \widehat B_{21}$ are identically zero.\\[-1.5cm]
  }\label{fig:1}
  \end{center}
\end{figure}

For the numerical solution of \eqref{eq:model_problem_elliptic_num_ex}, 
we use a finite element approach with globally continuous and piecewise affine
ansatz functions subordinate to Friedrichs--Keller triangulations of $\Omega$ and the interval 
$(0,L)$, respectively, 
with partition widths $h = 1/N$, $N \in \N$. 
As a surrogate model for $\widehat B\colon (0, 0.5) \to \R^{2\times 2}$,
we choose the function
%\[
%	\widehat B_\mathrm{sur}\colon (0, 0.5) \to \R^{2\times 2},
%	\quad
%	\widehat B_\mathrm{sur}(r) := 
%    \begin{cases}
%    \left (
%    55.624\, r^3 - 61.322\,r^2 - 1.641\,r + 10.101 
%    \right ) \Id& \text{ if } r \leq 0.45,
%    \\
%    \displaystyle
%    \left (0.001 + 
%    \frac{2.013}{124.128\,r^2 - 100.271\,r + 20.986}
%    \right ) \Id & \text{ if } r > 0.45.
%    \end{cases}
%\]
\begin{equation}
\label{eq:surrogate_polynomial}
	\widehat B_\mathrm{sur}\colon (0, 0.5) \to \R^{2\times 2},
	\quad
	\widehat B_\mathrm{sur}(r) := 
    \left (
    55.624\, r^3 - 61.322\,r^2 - 1.641\,r + 10.101 
    \right ) \Id.
\end{equation}
This model has been identified
in 
an offline phase
as a suitable substitute for $\widehat B$ 
on the interval $[0.1, 0.4]$ relevant for 
\eqref{eq:model_problem_elliptic_num_ex} 
(see \cref{fig:1}) and trivially satisfies the requirements in
Assumptions 
\ref{ass:standing_grad_proj}
and \ref{ass:standing:Newton}
(at least after a
redefinition on $\R \setminus [0.1, 0.4]$---we use 
a rational continuation on this set).

% linewidth 2, helvetica 14

\pagebreak

The integration 
of the term $(B(q)\nabla u) \cdot \nabla v$ in the weak formulation of the state equation of 
\eqref{eq:model_problem_elliptic_num_ex} (and the associated adjoint equation etc.) is realised in our implementation
by lumping the mass of $B(q)$ to the midpoints of the cells of the underlying mesh.
To allow for a fair comparison 
of the gradient-projection algorithm discussed in 
\cref{subsec:gradproj} and the semismooth Newton method developed 
in \cref{subsec:semiNewton}, for both of these algorithms, we use the 
termination criterion
\begin{equation}
\label{eq:randomeq263636}
	\left \| q_k - \proj_{\Qad}
\left (
q_k  - 10^{4} \cdot \nabla_{H^1(0,L)} \hat J(q_k)
\right )
\right \|_{H^1(0,L)} < 10^{-8}
\end{equation}
motivated by the necessary optimality condition in \cref{lem:grad_opt_reform} in our numerical experiments.
Note that the quantity \eqref{eq:randomeq263636} is 
not naturally computed in the two considered 
algorithms and induces an additional overhead (in particular 
the additional evaluation of the projection $\pi_{\Qad}$).
The computation time necessary for checking 
\eqref{eq:randomeq263636} is not included
in the runtimes that we document below. 
It is necessary 
to consider a separate termination criterion for the comparison of the 
gradient-projection method and the semismooth Newton method
as these two algorithms typically use different
quantities to determine whether a given tolerance has been reached
(the gradient-projection method relies on the norm $\|q_{k+1} - q_k\|_Q$ while the Newton method uses the norm of the residue).
As the initial guess, we use $q_0 \equiv 0$ in all experiments.

\begin{figure}[ht]
\begin{center}
  \includegraphics[height=4.8cm]{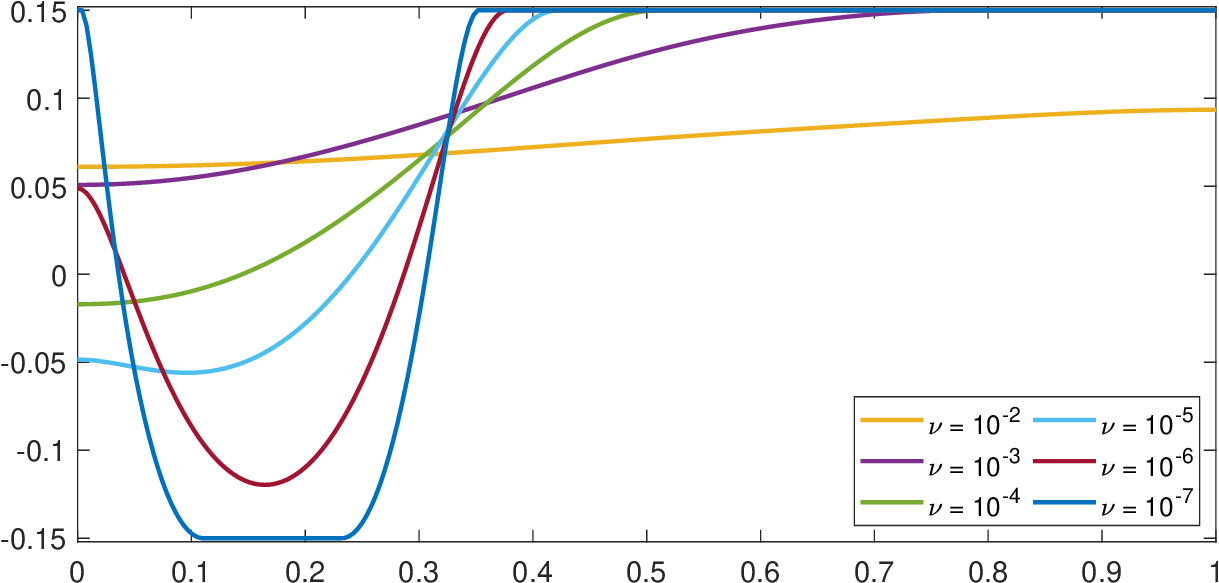}
  \caption{The solutions $\bar q$ of \eqref{eq:model_problem_elliptic_num_ex} determined 
  by means of the semismooth Newton iteration in \eqref{eq:Newton_step_def}
  for $h = 1/256$ and $\nu = 10^{-2}, \ldots, 10^{-7}$. 
  For $\nu = 10^{-5}$, $10^{-6}$, and $10^{-7}$, path-following w.r.t.\ $\nu$ was used to ensure convergence.
  \cref{alg:gradproj} produces the exact same solutions. 
  }\label{fig:2}
  \end{center}
  \vspace{-0.4cm}
\end{figure}

\begin{figure}[ht]
\begin{center}
  \includegraphics[width=6.8cm]{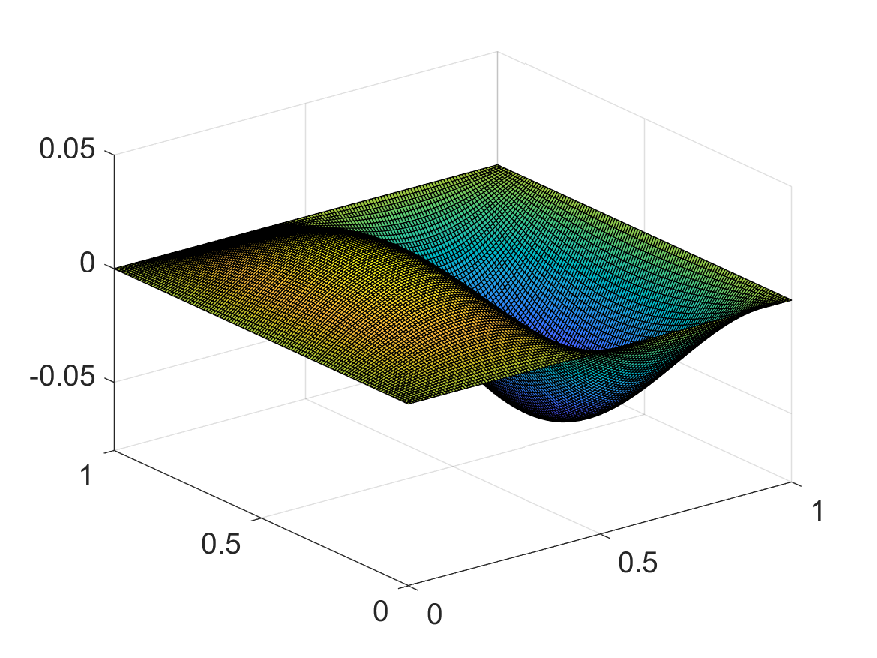}
  ~~~
  \includegraphics[width=6.8cm]{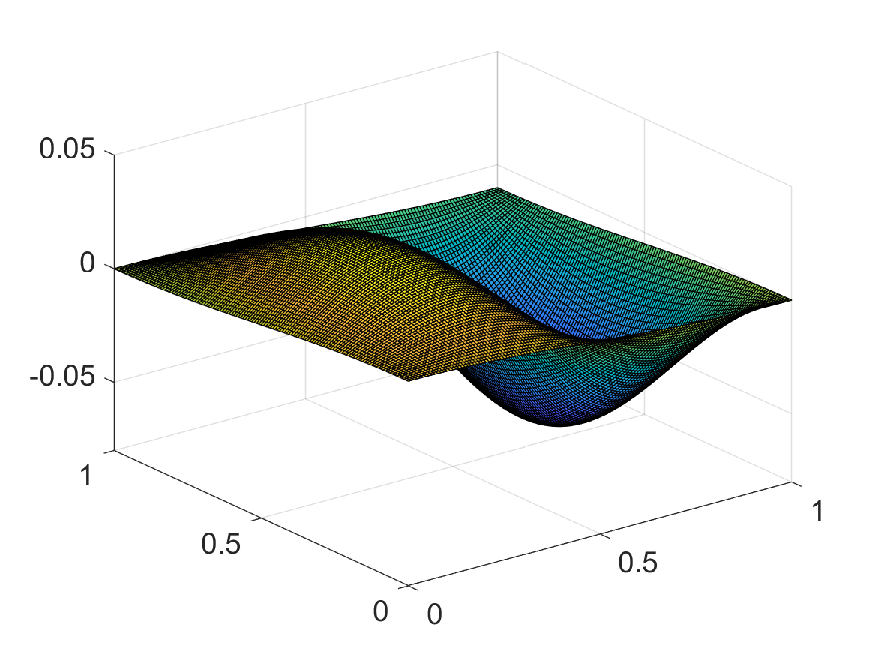}
  \vspace{-0.2cm}
  \caption{The function 
  $S(q_{\mathrm{exact}})$ (left) and 
  the observed data $u_\mathrm{obs}$ (right)
  in the situation 
  of \eqref{eq:model_problem_elliptic_num_ex}.\\[-1cm]
  }\label{fig:3}
  \end{center}
\end{figure}

The results obtained in the situation above by means of 
\cref{alg:gradproj} and the Newton iteration in \eqref{eq:Newton_step_def}
can be seen in Figures \ref{fig:2} and \ref{fig:3} and \cref{tab:1,tab:2}.
Figure~\ref{fig:2} shows that whether and to what extent the bounds $\pm 0.15$ are active depends on the size of the regularization parameter $\nu$. Since $u_\mathrm{obs}$ is only a perturbation of $S(q_{\mathrm{exact}})$, the optimal control $\bar q$ gets closer but does not converge to $q_{\mathrm{exact}}$ as $\nu$ approaches $0$.

\vspace{0.3cm}

% linewidth 2, helvetica 14
\pagebreak

\begin{table}[H]
\centering
\begin{tabular}{c |c | c | c | c | c | c | c}
\noalign{\smallskip}
\diagbox{$\nu$}{$h$}  &  $\frac{1}{32}$ &  $\frac{1}{64}$ & $\frac{1}{128}$ & $\frac{1}{256}$   & $\frac{1}{512}$   & $\frac{1}{1024}$    \\
\noalign{\smallskip}\hline\noalign{\smallskip}
$10^{-1}$ & 
\begin{tabular}{@{}c@{}} $4$ \\   $0.209$s    \end{tabular}  & 
\begin{tabular}{@{}c@{}} $4$ \\   $0.310$s    \end{tabular}  & 
\begin{tabular}{@{}c@{}} $4$ \\   $0.713$s    \end{tabular}  & 
\begin{tabular}{@{}c@{}} $4$ \\   $3.102$s    \end{tabular}  & 
\begin{tabular}{@{}c@{}} $4$ \\   $15.373$s    \end{tabular}  & 
\begin{tabular}{@{}c@{}} $4$ \\   $66.553$s    \end{tabular}  \\ 
\noalign{\smallskip}\hline\noalign{\smallskip}
$10^{-2}$ & 
\begin{tabular}{@{}c@{}} $5$ \\   $0.252$s    \end{tabular}  & 
\begin{tabular}{@{}c@{}} $6$ \\   $0.426$s    \end{tabular}  & 
\begin{tabular}{@{}c@{}} $6$ \\   $0.962$s    \end{tabular}  & 
\begin{tabular}{@{}c@{}} $6$ \\   $4.313$s    \end{tabular}  & 
\begin{tabular}{@{}c@{}} $6$ \\   $23.515$s    \end{tabular}  & 
\begin{tabular}{@{}c@{}} $6$ \\   $97.527$s    \end{tabular}  \\
\noalign{\smallskip}\hline\noalign{\smallskip}
$10^{-3}$ & 
\begin{tabular}{@{}c@{}} $6$ \\   $0.254$s    \end{tabular}  & 
\begin{tabular}{@{}c@{}} $7$ \\   $0.434$s    \end{tabular}  & 
\begin{tabular}{@{}c@{}} $7$ \\   $1.035$s    \end{tabular}  &  
\begin{tabular}{@{}c@{}} $7$ \\   $4.856$s    \end{tabular}  &  
\begin{tabular}{@{}c@{}} $7$ \\   $24.399$s    \end{tabular}  & 
\begin{tabular}{@{}c@{}} $7$ \\   $103.829$s    \end{tabular}  \\
\noalign{\smallskip}\hline\noalign{\smallskip}
$10^{-4}$ & 
\begin{tabular}{@{}c@{}} $6$ \\   $0.264$s    \end{tabular}  &
\begin{tabular}{@{}c@{}} $6$ \\   $0.404$s    \end{tabular}  &
\begin{tabular}{@{}c@{}} $6$ \\   $0.951$s    \end{tabular}  &
\begin{tabular}{@{}c@{}} $7$ \\   $4.935$s    \end{tabular}  &
\begin{tabular}{@{}c@{}} $7$ \\   $24.612$s    \end{tabular}  &
\begin{tabular}{@{}c@{}} $7$ \\   $105.368$s    \end{tabular}  
\end{tabular} 
\caption{Number of Newton iterations and computation times 
(including matrix assembly and FE-overhead but excluding 
the time for checking \eqref{eq:randomeq263636}) 
needed to achieve \eqref{eq:randomeq263636} 
for various mesh widths $h$ and
Tikhonov parameters $\nu$ in the situation of \eqref{eq:model_problem_elliptic_num_ex}.
The computations were carried out in MATLAB/R2024b on a system with an 
IntelCore i7-10510U CPU with 16.0GB RAM.}
\label{tab:1}
\end{table}

\begin{table}[H]
\centering
\begin{tabular}{c |c | c | c | c | c | c | c}
\noalign{\smallskip}
\diagbox{$\tau$}{$h$}  &  $\frac{1}{32}$ &  $\frac{1}{64}$ & $\frac{1}{128}$ & $\frac{1}{256}$  & $\frac{1}{512}$   & $\frac{1}{1024}$    \\
\noalign{\smallskip}\hline\noalign{\smallskip}
$325$ & 
\begin{tabular}{@{}c@{}}   $41$  \\    $0.219$s     \end{tabular} & 
\begin{tabular}{@{}c@{}}   $41$  \\    $0.572$s     \end{tabular} & 
\begin{tabular}{@{}c@{}}   $41$  \\    $2.046$s     \end{tabular} & 
\begin{tabular}{@{}c@{}}   $41$  \\    $9.834$s     \end{tabular} & 
\begin{tabular}{@{}c@{}}   $41$  \\    $49.640$s     \end{tabular} & 
\begin{tabular}{@{}c@{}}   $41$  \\    $217.996$s     \end{tabular}   \\ 
\noalign{\smallskip}\hline\noalign{\smallskip}
$350$ & 
\begin{tabular}{@{}c@{}}   $38$  \\    $0.214$s     \end{tabular} & 
\begin{tabular}{@{}c@{}}   $38$  \\    $0.535$s     \end{tabular} & 
\begin{tabular}{@{}c@{}}   $38$  \\    $1.920$s     \end{tabular} & 
\begin{tabular}{@{}c@{}}   $38$  \\    $9.016$s     \end{tabular} & 
\begin{tabular}{@{}c@{}}   $38$  \\    $46.277$s     \end{tabular} & 
\begin{tabular}{@{}c@{}}   $38$  \\    $204.618$s     \end{tabular}   \\ 
\noalign{\smallskip}\hline\noalign{\smallskip}
$375$ & 
\begin{tabular}{@{}c@{}}   $35$  \\    $0.211$s     \end{tabular} & 
\begin{tabular}{@{}c@{}}   $36$  \\    $0.507$s     \end{tabular} & 
\begin{tabular}{@{}c@{}}   $36$  \\    $1.828$s     \end{tabular} & 
\begin{tabular}{@{}c@{}}   $35$  \\    $8.374$s     \end{tabular} & 
\begin{tabular}{@{}c@{}}   $36$  \\    $43.636$s     \end{tabular} & 
\begin{tabular}{@{}c@{}}   $36$  \\    $187.543$s     \end{tabular}   \\  
\noalign{\smallskip}\hline\noalign{\smallskip}
$400$ & 
\begin{tabular}{@{}c@{}}   $47$  \\    $0.253$s     \end{tabular} & 
\begin{tabular}{@{}c@{}}   $49$  \\    $0.661$s     \end{tabular} & 
\begin{tabular}{@{}c@{}}   $49$  \\    $2.372$s     \end{tabular} & 
\begin{tabular}{@{}c@{}}   $48$  \\    $11.273$s     \end{tabular} & 
\begin{tabular}{@{}c@{}}   $48$  \\    $57.664$s     \end{tabular} & 
\begin{tabular}{@{}c@{}}   $48$  \\    $247.478$s     \end{tabular}   \\ 
\noalign{\smallskip}\hline\noalign{\smallskip}
$425$ & 
\begin{tabular}{@{}c@{}}   $71$  \\    $0.317$s     \end{tabular} & 
\begin{tabular}{@{}c@{}}   $75$  \\    $0.890$s     \end{tabular} & 
\begin{tabular}{@{}c@{}}   $75$  \\    $3.508$s     \end{tabular} & 
\begin{tabular}{@{}c@{}}   $73$  \\    $16.679$s     \end{tabular} & 
\begin{tabular}{@{}c@{}}   $74$  \\    $87.402$s     \end{tabular} & 
\begin{tabular}{@{}c@{}}   $74$  \\    $376.253$s     \end{tabular}   \\ 
\noalign{\smallskip}\hline\noalign{\smallskip}
$450$ & 
\begin{tabular}{@{}c@{}}   $131$  \\    $0.502$s     \end{tabular} & 
\begin{tabular}{@{}c@{}}   $143$  \\    $1.587$s     \end{tabular} & 
\begin{tabular}{@{}c@{}}   $145$  \\    $6.589$s     \end{tabular} & 
\begin{tabular}{@{}c@{}}   $137$  \\    $30.863$s     \end{tabular} & 
\begin{tabular}{@{}c@{}}   $141$  \\    $164.685$s     \end{tabular} & 
\begin{tabular}{@{}c@{}}   $141$  \\    $723.016$s     \end{tabular}   \\ 
\end{tabular} 
\caption{Number of gradient-projection steps 
and computation times 
(including matrix assembly and FE-overhead but excluding 
the time for checking \eqref{eq:randomeq263636}) 
needed by \cref{alg:gradproj} 
to achieve \eqref{eq:randomeq263636} for various mesh widths $h$ and step sizes $\tau$
with $\nu = 10^{-3}$. The optimal step size is (approximately) $\tau = 375$.
For $\tau = 475$ and larger step sizes, 
\cref{alg:gradproj} does not converge. The computing infrastructure 
was the same as in \cref{tab:1}.}
\label{tab:2}
\end{table}

As can be seen in \cref{tab:1,tab:2}, both \cref{alg:gradproj} 
and the semismooth Newton method from \cref{th:convergence_Newton} 
exhibit a mesh-independent convergence behaviour, i.e.\ 
achieve \eqref{eq:randomeq263636} with a number of iterations
which is stable w.r.t.\ mesh refinement. It can further be seen 
that the semismooth Newton method indeed requires far less 
iterations and, thus, evaluations of the projection $\pi_{\Qad}$ than the 
projected gradient-descent method to solve \eqref{eq:model_problem_elliptic_num_ex}
with $\nu = 10^{-3}$. This lower number of projection evaluations also manifests 
itself in lower runtimes (with a sole exception
for $h = 1/32$, where the effort for assembling the Newton matrix 
outweighs the benefit of having a lower number of projection operations);
see \cref{tab:1,tab:2}. 
When compared with \cref{alg:gradproj} with the optimal choice of $\tau$,
the semismooth Newton method achieves a speed-up of $\approx 1.8$
for small $h$ and, thus, outperforms the classical 
gradient-descent method significantly, 
even if the hyperparameters of the latter
have been optimised in an offline phase.\pagebreak

\begin{figure}[ht]
\begin{center}
  \includegraphics[width=6.8cm]{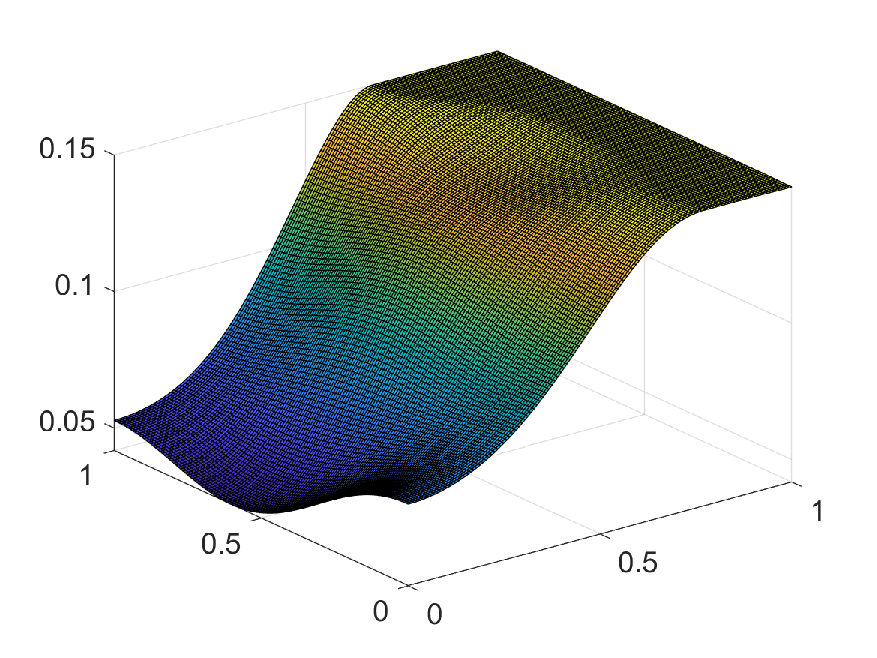}
  ~~~
  \includegraphics[width=6.8cm]{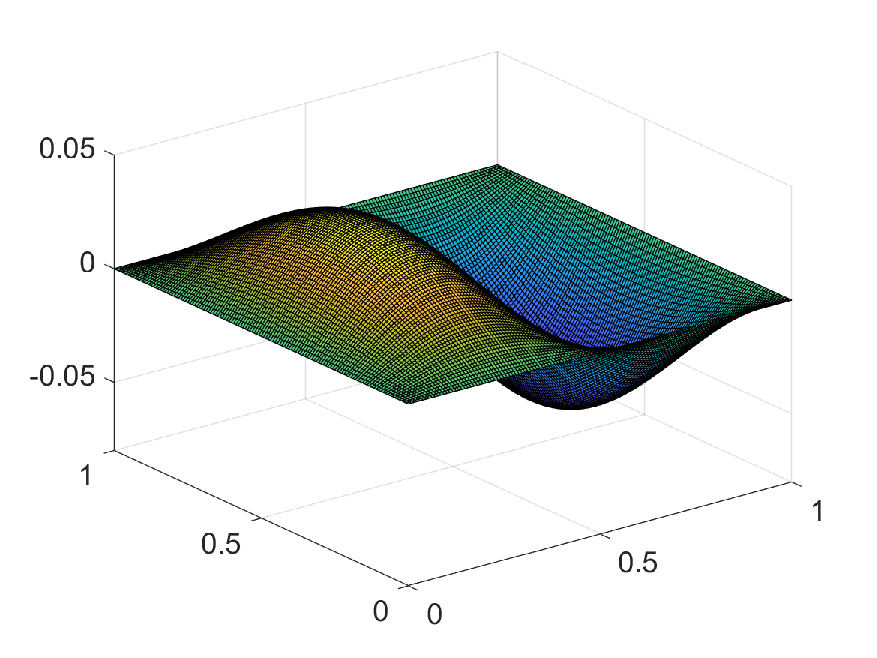}
  \vspace{-0.2cm}
  \caption{The solution $\bar q$ (left) and 
  its state $\bar u$ (right) for $\nu = 10^{-3}$ and $h=1/128$ in the situation that the ansatz space $H^1(0,1)$
for $q$ in problem \eqref{eq:model_problem_elliptic_num_ex}
is replaced by $H^1(\Omega)$; cf.\ \cref{fig:3}.
  }\label{fig:4}
  \end{center}
\end{figure}

Although currently not covered by the theoretical results of \cref{sec:abstractsetting}, we remark that the benefits of our semismooth Newton approach are even more pronounced 
if $q$ is
allowed to vary not only in one 
but in all spatial directions.
If $q$ is sought, for example, in the space 
$H^1(\Omega)$ instead of $H^1(0,1)$ in 
\eqref{eq:model_problem_elliptic_num_ex},
then  
evaluating $\pi_{\Qad}$ becomes 
much more expensive and the cost of the projections 
also scales 
less favourably with $h$. 
This entails that the reduction of the number of
projection evaluations achieved by the semismooth Newton method 
has an even higher impact and that the semismooth 
Newton algorithm outperforms the gradient-projection method 
to an even larger extent; see \cref{tab:3} and \Cref{fig:4}.
\vspace{0.2cm}

\begin{table}[ht]
\centering
\begin{tabular}{c |c | c | c | c | c | c }
\noalign{\smallskip}
$h$  &  $\frac{1}{16}$ &  $\frac{1}{32}$ & $\frac{1}{64}$ &  $\frac{1}{128}$   & $\frac{1}{256}$    \\
\noalign{\smallskip}\hline\noalign{\smallskip}
\begin{tabular}{@{}c@{}}  Projected Gradient Descent  \\   with $\tau =325$    \end{tabular} & 
\begin{tabular}{@{}c@{}}   $45$  \\    $0.325$s    \end{tabular} & 
\begin{tabular}{@{}c@{}}   $45$  \\    $0.939$s    \end{tabular} & 
\begin{tabular}{@{}c@{}}   $45$  \\    $4.669$s    \end{tabular} & 
\begin{tabular}{@{}c@{}}   $45$  \\    $21.773$s    \end{tabular} & 
\begin{tabular}{@{}c@{}}   $45$  \\    $132.760$s    \end{tabular}   \\ 
\noalign{\smallskip}\hline\noalign{\smallskip}
\begin{tabular}{@{}c@{}}  Projected Gradient Descent  \\   with $\tau =375$    \end{tabular} & 
\begin{tabular}{@{}c@{}}   $38$  \\    $0.301$s    \end{tabular} & 
\begin{tabular}{@{}c@{}}   $39$  \\    $0.744$s    \end{tabular} & 
\begin{tabular}{@{}c@{}}   $39$  \\    $4.214$s    \end{tabular} & 
\begin{tabular}{@{}c@{}}   $39$  \\    $19.306$s    \end{tabular} & 
\begin{tabular}{@{}c@{}}   $39$  \\    $111.248$s    \end{tabular}   \\ 
\noalign{\smallskip}\hline\noalign{\smallskip}
\begin{tabular}{@{}c@{}}  Projected Gradient Descent  \\   with $\tau =425$    \end{tabular} & 
\begin{tabular}{@{}c@{}}   $76$  \\    $0.532$s    \end{tabular} & 
\begin{tabular}{@{}c@{}}   $86$  \\    $1.775$s    \end{tabular} & 
\begin{tabular}{@{}c@{}}   $87$  \\    $8.548$s    \end{tabular} & 
\begin{tabular}{@{}c@{}}   $88$  \\    $43.104$s    \end{tabular} & 
\begin{tabular}{@{}c@{}}   $88$  \\    $245.461$s    \end{tabular}   \\ 
\noalign{\smallskip}\hline\noalign{\smallskip}
\noalign{\smallskip}\hline\noalign{\smallskip}
Semismooth Newton & 
\begin{tabular}{@{}c@{}}   $7$  \\    $0.233$s    \end{tabular} & 
\begin{tabular}{@{}c@{}}   $7$  \\    $0.403$s    \end{tabular} & 
\begin{tabular}{@{}c@{}}   $7$  \\    $1.439$s    \end{tabular} & 
\begin{tabular}{@{}c@{}}   $7$  \\    $7.041$s    \end{tabular} & 
\begin{tabular}{@{}c@{}}   $7$  \\    $54.828$s    \end{tabular}   \\ 
\end{tabular} 
\caption{Number of iterations  
and computation times 
(including matrix assembly and FE-overhead but excluding 
the time for checking \eqref{eq:randomeq263636}) 
needed by \cref{alg:gradproj} 
and the semismooth Newton method
\eqref{eq:Newton_step_def} 
to solve \eqref{eq:model_problem_elliptic_num_ex} 
for different mesh widths $h$
in the situation that the ansatz space $H^1(0,1)$
for $q$ 
is replaced by $H^1(\Omega)$.
The problem data, discretisation, and termination criterion
were exactly the same as in the experiments of 
Figures \ref{fig:2} and \ref{fig:3} and \cref{tab:1,tab:2}
and the Tikhonov parameter was chosen as $\nu = 10^{-3}$. It can be seen
that the semismooth Newton method provides an even higher speedup
than before and more than halves the computation times
for small $h$.
}
\label{tab:3}
\end{table}

\vspace{0.2cm}

As already mentioned, we expect 
that for fully distributed design variables $q$, it is still 
possible to perform the convergence analysis in 
\cref{subsec:semiNewton} provided the 
PDE governing the identification problem allows for the derivation of
Gröger-regularity results in suitable spaces. 
This topic will be investigated in future work. 

\appendix

\section{Newton differentiability for one-dimensional variational inequalities with box constraints}
\label{sec:appendix}

In this appendix, we prove the Newton differentiability result
in \cref{th_main_semismooth_main_body} needed 
for the discussion of the semismooth Newton method in \cref{subsec:semiNewton}. 
To prove \cref{th_main_semismooth_main_body}, we employ arguments that have 
been developed for classical 
$H^1(\Omega)$-elliptic
obstacle-problems with vanishing boundary conditions 
on domains of arbitrary dimension 
in \cite{ChristofWachsmuthUniObstacle,ChristofWachsmuthBiObstacle,Rauls2020}. 
We present these arguments in a complete and largely self-contained manner in what follows 
to show that, in the one-dimensional setting, 
it is indeed possible to obtain the Newton-differentiability of 
the map 
$\Sigma$ in 
\cref{subsec:semiNewton}
for $\MM([0,L])$-inputs and to 
avoid the assumption of zero boundary conditions---two 
points that are essential for the  analysis of \cref{subsec:semiNewton} and 
the application background of \eqref{eq:P}  but have not been discussed in 
\cite{ChristofWachsmuthUniObstacle,ChristofWachsmuthBiObstacle,Rauls2020}. 
We also use this opportunity to highlight some shortcuts, generalisations and improvements which are possible
in the one-dimensional setting but not in the multi-dimensional one studied in 
\cite{ChristofWachsmuthUniObstacle,ChristofWachsmuthBiObstacle,Rauls2020}.
We begin by listing our standing assumptions: 

\begin{assumption}[{Standing assumptions for Appendix A}]
\label{ass:standing:appendix_A}~

\begin{enumerate}[label=\roman*)]
\item $L>0$ is given and fixed;
\item
\label{ass:standing:appendix_A:ii}
$Q$ is a real separable Hilbert space with Riesz isomorphism 
$\iota_R\colon Q^* \to Q$ such that:
\begin{enumerate}
\item $Q \hookrightarrow C([0,L])$ compactly and densely;
\item 
$\varphi(q)\in Q$ for all $q \in Q$ and all 
nondecreasing, piecewise affine-linear and 
continuous functions 
$\varphi\colon \R \to \R$
with $\varphi(0) = 0$ 
(acting by superposition), and for every
such $\varphi$ 
there exists a constant $C >0$ 
such that, for every $q \in Q$, we have 
\begin{equation}
\label{eq:maxprop}
\|\varphi(q)\|_Q \leq C \|q\|_Q,
\qquad
\left (
q, \varphi(q)
\right )_Q
\geq
0
\qquad
\text{and}
\qquad
\left (
q, \varphi(q)
\right )_Q = 0
~\Leftrightarrow
~
\varphi(q)=0;
\end{equation}
\item
for every $\psi \in C^\infty([0,L])$, the map $Q \ni q \mapsto \psi q \in Q$
is well-defined and in $\LLLL(Q, Q)$;
\end{enumerate}
\item
\label{ass:standing:appendix_A:iii}
$a,b \in C([0,L])$ are given such that 
there exist $\varepsilon > 0$ and $\tilde q \in Q$
satisfying $a \leq \tilde q \leq b$, $a + \varepsilon \leq b$ pointwise in $[0,L]$,
and $\tilde \mu := \iota_R^{-1}(\tilde q)\in \MM([0,L]) \subset Q^*$.
\end{enumerate}
\end{assumption}

Henceforth, we use the notation 
\begin{equation*}
\begin{aligned}
K_a &:=  \{q \in Q \mid a \leq q \text{ pointwise in } [0,L] \},
\\
K^b&:=  \{q \in Q \mid  q \leq b \text{ pointwise in } [0,L] \},
\\
K_a^b&:=  \{q \in Q \mid  a \leq q \leq b \text{ pointwise in } [0,L] \}.
\end{aligned}
\end{equation*}
Note that all of these sets are convex, closed and nonempty by \cref{ass:standing:appendix_A}. 
We further denote by $\Sigma_a$, $\Sigma^b$ and $\Sigma_a^b$
the solution maps $Q^* \ni g \mapsto \hat q \in Q$ of the variational inequalities 
\begin{equation*}
\text{Find }
\hat q \in K
\text{ such that }
\left ( \hat q, q-\hat q \right )_{Q}
- \left \langle g,q - \hat q \right \rangle_Q
\geq 0
~~\forall q \in K,
\end{equation*}
with $K = K_a, K^b, K_a^b$, respectively (so that 
$\Sigma_a^b$ corresponds 
to the function $\Sigma$ in \cref{subsec:semiNewton}). 
From standard results, we obtain: 

\begin{lemma}[Basic properties of $\Sigma_a$, $\Sigma^b$ and $\Sigma_a^b$]
\label{lem:basic_props_sigma}
The solution maps $\Sigma_a$, $\Sigma^b$ and $\Sigma_a^b$ are well-defined 
and non-expansive as functions from $Q^*$ to $Q$.
Moreover, it holds
\[
	\Sigma_a (\lambda g_1 + (1-\lambda) g_2) 
	\leq
	\lambda \Sigma_a (g_1)  + (1-\lambda) \Sigma_a (g_2) 
	\text{ pointwise in } [0,L] \quad 
	\forall \lambda \in [0,1]
	\quad
	\forall g_1, g_2 \in Q^*,
\]
and $\Sigma_a\colon Q^* \to Q$ is directionally differentiable, with the directional derivative
$\delta := \Sigma_a'(g;h) \in Q$ at $g \in Q^*$ in direction $h \in Q^*$ being characterised by the auxiliary problem 
\begin{equation}
\label{eq:dirdiff_VI}
\text{Find }
	\delta \in \TT(\hat q; K_a) \cap \ker(\eta)
	\text{ such that }
	\left (\delta, p - \delta\right )_{Q}
- \left \langle h, p - \delta \right \rangle_Q
\geq 0
~~\forall p  \in \TT(\hat q; K_a) \cap \ker(\eta).
\end{equation}
Here, $\TT(\hat q; K_a) = \mathrm{cl}_Q(\R^+(K_a - \hat q))$ is the tangent cone of 
$K_a$ in $Q$ at $\hat q := \Sigma_a(g)$ and $ \ker(\eta) \subset Q$ the kernel of the functional 
$\eta := \iota_R^{-1}(\hat q) - g \in Q^*$. 
In particular, $\|\delta\|_Q^2 = \langle h,\delta \rangle_Q$. 
\end{lemma}

\begin{proof}
The well-definedness and non-expansiveness of the maps 
 $\Sigma_a$, $\Sigma^b$, and $\Sigma_a^b$ 
 follow from \cite[Theorem 4:3-3.1]{Rodrigues1987}. 
The pointwise convexity of $\Sigma_a$ is obtained along standard lines:
If $\lambda \in [0,1]$ and $g_1, g_2 \in Q^*$ are given, then 
$\hat q_1 := \Sigma_a(g_1)$, $\hat q_2 := \Sigma_a(g_2)$,  
$\hat q_{12} := \Sigma_a(\lambda g_1 + (1-\lambda)g_2)$ 
and
$w :=  \hat q_{12} - \lambda \hat q_1 - (1 - \lambda)\hat q_2$ satisfy 
$a \leq \hat q_1 + \max(0, w) \in Q$, $a \leq \hat q_2 + \max(0, w) \in Q$, and
\[
Q \ni \hat q_{12} - \max(0, w) = \hat q_{12} - \max(0, \hat q_{12} - \lambda \hat q_1 - (1 - \lambda) \hat q_2) 
= \min(\hat q_{12}, \lambda \hat q_1 + (1 - \lambda)\hat q_2) \geq a
\]
by our assumptions on $Q$, and it follows from the VIs for
$\hat q_1$, $\hat q_2$ and $\hat q_{12}$ that 
\begin{align*}
\lambda\left [ 
\left ( \hat q_1, \max(0, w) \right )_{Q}
- \left \langle g_1, \max(0, w) \right \rangle_Q
\right ]
\geq 0, 
\\
(1-\lambda)\left [ 
\left ( \hat q_2, \max(0, w) \right )_{Q}
- \left \langle g_2, \max(0, w) \right \rangle_Q
\right ]
\geq 0, 
\\
\left ( \hat q_{12}, - \max(0, w) \right )_{Q}
- \left \langle  \lambda g_1 + (1-\lambda) g_2, -\max(0, w) \right \rangle_Q
\geq 0.
\end{align*}
Adding the above inequalities yields
$0 \geq \left ( w , \max(0, w)\right )_Q$ and, 
by \eqref{eq:maxprop}, 
$\max(0,w) = 0$, i.e.\ $w \leq 0$ in $[0,L]$. 
This proves the pointwise convexity as desired. 
To establish the directional differentiability result for $\Sigma_a$, 
it suffices to note that \cref{ass:standing:appendix_A} ensures that 
$Q$ is a vector lattice in the sense of 
\cite[Definition 4.6]{Wachsmuth2016:2}
satisfying 
\cite[Assumption 4.7]{Wachsmuth2016:2}
and that $K_a$ is a set with lower bound in the sense of 
\cite[Definition 4.9]{Wachsmuth2016:2}.
In combination with \cite[Theorem 4.18]{Wachsmuth2016:2},
this implies that $K_a$ is polyhedric in the sense of 
\cite[Definition 3.1]{Wachsmuth2016:2}, which in turn implies the 
asserted directional differentiability properties of 
$\Sigma_a$ by \cite[Theorem 5.2]{Wachsmuth2016:2}.
To 
obtain the identity $\|\delta\|_Q^2 = \langle h,\delta \rangle_Q$,
it suffices to choose 
$p = 0$ and $p= 2 \delta$ in \eqref{eq:dirdiff_VI},
keeping in mind that $\TT(\hat q; K_a) \cap \ker(\eta)$ is a cone. 
\end{proof}

As a generalised differential for $\Sigma_a$, we consider the following object:

\begin{definition}[Strong--weak Bouligand differential]%
\label{def:Bsubdifferential}%
For all $g \in Q^*$, we define the 
strong--weak Bouligand differential 
$\partial_B^{sw}\Sigma_a (g) \subset \LLLL(Q^*,Q)$ by
\begin{equation}
\label{def:Bsubdiff}
\begin{aligned}
\partial_B^{sw}\Sigma_a (g)
:=
\left \{
G \in \LLLL(Q^*,Q)
~~\Bigg |
\exists \{g_n\} \subset \DD \colon g_n \to g \text{ in } Q^*,~ \Sigma_a'(g_n) \WOT G \text{ in } \LLLL(Q^*, Q)
\right \}.
\end{aligned}
\end{equation}
Here, $\DD$ denotes the set of all $g \in Q^*$ with $\Sigma_a'(g) := \Sigma_a'(g; \cdot) \in \LLLL(Q^*, Q)$
(i.e.\ the set of all Gâteaux points of the function $\Sigma_a$) 
and $\WOT$ convergence w.r.t.\ the weak operator topology.
\end{definition}

The following lemma collects some basic properties of the differential $\partial_B^{sw}\Sigma_a$.

\begin{lemma}[Properties of $\partial_B^{sw}\Sigma_a$]\label{lem:props_swBD}~
\begin{enumerate}[label=\roman*)]
\item\label{lem:props_swBD:i} For every $g \in Q^*$, it holds $\partial_B^{sw}\Sigma_a(g) \neq \emptyset$ and 
$\|G\|_{\LLLL(Q^*, Q)} \leq 1$ for all $G \in \partial_B^{sw}\Sigma_a(g)$. 
\item\label{lem:props_swBD:ii} If $\{g_n\} \subset Q^*$ and $\{G_n\} \subset \LLLL(Q^*, Q)$ 
satisfy 
$g_n\to g$ in $Q^*$ for some $g \in Q^*$, 
$G_n \in \partial_B^{sw} \Sigma_a (g_n)$ for all $n \in \mathbb{N}$,
and $G_n \WOT G$  in $\LLLL(Q^*, Q)$ for some 
$G \in \LLLL(Q^*, Q)$, then $G \in \partial_B^{sw} \Sigma_a(g)$. 
\item\label{lem:props_swBD:iii} For all $g \in Q^*$ and all $G \in \partial_B^{sw}\Sigma_a(g)$, we have 
\[
\Sigma_a( g + h) - \Sigma_a(g) \geq G h~~\text{ pointwise in } [0,L] \qquad \forall h \in Q^*.
\]
\item\label{lem:props_swBD:iv} For every $g \in Q^*$ and $G \in \partial_B^{sw}\Sigma_a(g)$, it holds 
\begin{equation}
\label{eq:bouligand_inequality}
	\|G h\|_Q^2 \leq \left \langle h, Gh \right  \rangle_Q\qquad \forall h \in Q^*. 
\end{equation}
\end{enumerate}
\end{lemma} 
\begin{proof}
The assertions of the lemma follow straightforwardly from the 
assumptions on $Q$ and the mapping properties in \cref{lem:basic_props_sigma}. 
For detailed proofs of the statements \ref{lem:props_swBD:i}, \ref{lem:props_swBD:ii} and 
\ref{lem:props_swBD:iii}, we refer the reader to 
\cite[Corollary 2.6, Lemmas 2.7, 2.8]{ChristofWachsmuthUniObstacle}
(with $X = Q^*$, $Y = Q$, and $U = \{0\}$). To obtain point \ref{lem:props_swBD:iv},
one can argue along the lines of the proof of \cite[Lemma 3.8]{ChristofWachsmuthBiObstacle}:
Let $g \in Q^*$, $G \in \partial_B^{sw}\Sigma_a(g)$ and $h \in Q^*$ be given, 
and let $\{g_n\} \subset \DD$ satisfy the conditions in \eqref{def:Bsubdiff}
for $g$ and $G$. Then, it follows from 
\cref{lem:basic_props_sigma}
that $\|\Sigma_a'(g_n) h \|_Q^2 = \left \langle h, \Sigma_a'(g_n) h \right \rangle_Q$
holds for all $n$.
Since the function $\| \cdot \|_Q^2\colon Q \to \R$ is weakly lower semicontinuous
and since $\Sigma_a'(g_n) h \weakly Gh$ holds in $Q$ by the definition of
the weak operator topology, we can pass to the limit in 
$\|\Sigma_a'(g_n) h \|_Q^2 = \left \langle h, \Sigma_a'(g_n) h \right \rangle_Q$
to obtain $\|G h\|_Q^2 \leq \left \langle h, Gh \right  \rangle_Q$ as desired. 
This completes the proof.
\end{proof}

Using the last two lemmas, we can already prove the 
Newton differentiability of $\Sigma_a$ as a function from $\MM([0,L])$ to $C([0,L])$.

\begin{proposition}[{Newton differentiability of $\Sigma_a$ from $\MM([0,L])$ to  $C([0,L])$}]
\label{prop:Newton:1}
The mapping $\Sigma_a\colon Q^* \to Q$ is Newton differentiable 
with Newton derivative $\partial_B^{sw}\Sigma_a$ when considered as a 
function $\Sigma_a\colon \MM([0,L]) \to C([0,L])$, that is,
for all $g \in \MM([0,L])$, it holds 
\begin{equation}
\label{eq:NewtonDiff111}
\sup_{G \in \partial_B^{sw}\Sigma_a (g + h)} 
\frac{\|\Sigma_a (g + h) - \Sigma_a(g) - G h\|_{C([0,L])}}{\|h\|_{\MM([0,L])}} \to 0
\quad\text{ for } 0 <  \|h\|_{\MM([0,L])} \to 0.
\end{equation}
\end{proposition}
\begin{proof}
We argue along the lines of the proof of 
\cite[Theorem 2.12]{ChristofWachsmuthUniObstacle}, using weak-star instead of weak convergence: 
Note that our assumptions imply that every $g \in \MM([0,L])$ can be identified
with a unique element of $Q^*$. Thus, it makes sense to consider $\Sigma_a$ on $\MM([0,L])$.
Let us suppose now that \eqref{eq:NewtonDiff111} is wrong for some $g \in \MM([0,L])$.
Then, there exist $\{h_n\} \subset \MM([0,L]) \setminus \{0\}$, 
$\{G_n\} \subset \LLLL(Q^*, Q)$, and $\varepsilon > 0$
such that $h_n \to 0$ in $\MM([0,L])$, $G_n \in \partial_B^{sw}\Sigma_a (g + h_n)$ and 
\begin{equation}
\label{eq:randomeq273hge8}
	\liminf_{n \to \infty} \frac{\|\Sigma_a (g + h_n) - \Sigma_a(g) - G_n h_n\|_{C([0,L])}}{\|h_n\|_{\MM([0,L])}} \geq \varepsilon.
\end{equation}
By the theorem of Banach--Alaoglu and its WOT-variant (see \cite[Theorem 2.5]{ChristofWachsmuthUniObstacle}),
\cref{lem:props_swBD}\ref{lem:props_swBD:i}
and the  properties of $Q$, 
we may assume w.l.o.g.\ 
that $e_n := h_n/ \| h_n\|_{\MM([0,L])}$ satisfies $e_n \weaklystar e$ in $\MM([0,L])\cong C([0,L])^*$ 
for some $e \in \MM([0,L])$ and that $G_n \WOT G$ holds in $\LLLL(Q^*, Q)$ 
for some $G \in \LLLL(Q^*, Q)$.
From Schauder's theorem,
the compactness of the embedding $Q \hookrightarrow C([0,L])$
and points \ref{lem:props_swBD:ii} and \ref{lem:props_swBD:iii} of 
\cref{lem:props_swBD}, 
we obtain that  we have $e_n \to e$ in $Q^*$,
$G \in \partial_B^{sw}\Sigma_a(g)$ and
	\[
		\Sigma_a(g) - \Sigma_a(g + h_n) \geq - G_n h_n
		\quad\text{as well as}\quad
		\Sigma_a(g + h_n) - \Sigma_a(g) \geq   G   h_n
		\quad\text{ pointwise in }[0,L]. 
	\]
The above yields
	\begin{equation*}
		0
		\geq  \frac{ \Sigma_a(g + h_n) - \Sigma_a(g) - G_n h_n}{\|h_n\|_{\MM([0,L])}} 
		\geq  \frac{ G h_n - G_n h_n}{\| h_n\|_{\MM([0,L])}} 
		=
		(G - G_n) e_n
		\quad \text{pointwise in }[0,L]
	\end{equation*}
and, after taking pointwise maxima, 
 	\begin{equation}
 	\label{eq:randomeq28bd8een}
 	\begin{aligned}
		\frac{ \left \| \Sigma_a(g + h_n) - \Sigma_a(g) - G_n h_n \right \|_{C([0,L])}}{\|h_n\|_{\MM([0,L])}} 
		&\leq
		 \left \|  (G - G_n) e_n \right \|_{C([0,L])}
		 \\
		 &\leq 
		 \left \|  (G - G_n) e \right \|_{C([0,L])} + \left \|  (G - G_n) (e_n - e) \right \|_{C([0,L])}.
		 \end{aligned}
	\end{equation}
Note that, by the definition of WOT-convergence, we have $(G - G_n) e \weakly 0$ in $Q$
and, thus,  $(G - G_n) e \to 0$ in $C([0,L])$ by the compactness of the embedding 
$Q \hookrightarrow C([0,L])$. On the other hand, 
we also obtain from the same embedding, 
\cref{lem:props_swBD}\ref{lem:props_swBD:i} and the convergence
$e_n \to e $ in $Q^*$
that there exists a constant $C>0$ with  
$\left \|  (G - G_n) (e_n - e) \right \|_{C([0,L])} \leq C \|e_n - e\|_{Q^*} \to 0$.
This shows that the left-hand side of \eqref{eq:randomeq28bd8een} 
converges to zero, in contradiction with \eqref{eq:randomeq273hge8},
and establishes the asserted Newton differentiability property. 
\end{proof}

To get the Newton differentiability of $\Sigma_a$ with values in $Q$, we next exploit 
\cref{lem:props_swBD}\ref{lem:props_swBD:iv}.

\begin{proposition}[{Newton differentiability of $\Sigma_a$ from $\MM([0,L])$ to  $Q$}]
\label{prop:energy_Newton_uni}
The mapping $\Sigma_a\colon Q^* \to Q$ is Newton differentiable 
with Newton derivative $\partial_B^{sw}\Sigma_a$ when considered as a 
map from $\MM([0,L])$ to $Q$, that is,
for all $g \in \MM([0,L])$, it holds 
\begin{equation}
\label{eq:NewtonDiff222}
\sup_{G \in \partial_B^{sw}\Sigma_a (g + h)} 
\frac{\|\Sigma_a (g + h) - \Sigma_a(g) - G h\|_{Q}}{\|h\|_{\MM([0,L])}} \to 0
\quad\text{ for } 0 <  \|h\|_{\MM([0,L])} \to 0.
\end{equation}
\end{proposition}

\begin{proof}
We argue along the lines of the proof of \cite[Theorem 4.4]{ChristofWachsmuthBiObstacle}:
Let us again assume that \eqref{eq:NewtonDiff222} does not hold, i.e.\ 
that there exist $g \in \MM([0,L])$,
 $\{h_n\} \subset \MM([0,L]) \setminus \{0\}$, 
$\{G_n\} \subset \LLLL(Q^*, Q)$ and $\varepsilon > 0$
such that 
$h_n \to 0$ in $\MM([0,L])$, $G_n \in \partial_B^{sw}\Sigma_a (g + h_n)$ and 
\begin{equation}
\label{eq:randomeq273hge8-2}
\liminf_{n \to \infty} \frac{\|\Sigma_a (g + h_n) - \Sigma_a(g) - G_n h_n\|_{Q}}{\|h_n\|_{\MM([0,L])}} \geq \varepsilon.
\end{equation}
We define $e_n := h_n / \|h_n\|_{\MM([0,L])}$ and
set $t_n := \|h_n\|_{\MM([0,L])}$ so that $h_n = t_n e_n$ and $0 < t_n \to 0$.
As before, we assume w.l.o.g.\ that 
$e_n \weaklystar e$ holds in $\MM([0,L])$ for some $e \in \MM([0,L])$,
that $e_n \to e$ in $Q^*$ and that $G_n \WOT G$ holds in $\LLLL(Q^*, Q)$
for some $G \in \partial_B^{sw}\Sigma_a(g)$. From 
\cref{prop:Newton:1}, we already know that 
\begin{equation}
\label{eq:randomeq7373h}
	\left \|
	\frac{\Sigma_a (g + t_n e_n) - \Sigma_a(g) - t_n G_n e_n}{t_n}
	 \right \|_{C([0,L])} \to 0
\end{equation}
holds for $n \to \infty$. Owing to the 
nonexpansiveness of $\Sigma_a\colon Q^* \to Q$
and \cref{lem:props_swBD}\ref{lem:props_swBD:i}, we further know that 
the sequence appearing in the $ C([0,L])$-norm in \eqref{eq:randomeq7373h} is bounded in $Q$. 
Combined, this yields that
\[
	\frac{\Sigma_a (g + t_n e_n) - \Sigma_a(g)}{t_n} - G_n e_n \weakly 0
	\quad \text{ in }Q. 
\]
Since $e_n \to e$ holds in $Q^*$, we obtain from the directional differentiability
and Lipschitz continuity properties of $\Sigma_a\colon Q^* \to Q$ in \cref{lem:basic_props_sigma} that
\begin{equation*}
\begin{aligned}
	&\left \| \frac{\Sigma_a (g + t_n e_n) - \Sigma_a(g)}{t_n} - \Sigma_a'(g;e) \right \|_Q
	\\
	&\qquad\leq
	\left \| \frac{\Sigma_a (g + t_n e ) - \Sigma_a(g)}{t_n} - \Sigma_a'(g;e) \right \|_Q
	+
	\left \| \frac{\Sigma_a (g + t_n e_n) - \Sigma_a (g + t_n e )}{t_n} \right \|_Q
	\\
	&\qquad\leq
	\left \| \frac{\Sigma_a (g + t_n e ) - \Sigma_a(g)}{t_n} - \Sigma_a'(g;e) \right \|_Q
	+
	\|e_n - e\|_{Q^*} \to 0
\end{aligned}
\end{equation*}
for $n \to \infty$. It thus follows that  $G_n e_n \weakly \Sigma_a'(g;e)$ for $n \to \infty$.
Let us define $\delta_n := G_n e_n$ and $\delta := \Sigma_a'(g;e)$. Then, we obtain from 
\cref{lem:props_swBD}\ref{lem:props_swBD:iv} that 
$\|\delta_n\|_Q^2 \leq \left \langle e_n, \delta_n \right  \rangle_Q$ holds for all $n$
and from 
\cref{lem:basic_props_sigma}
that  
$\|\delta \|_Q^2 = \left \langle e, \delta  \right  \rangle_Q$. 
In combination with the convergences $e_n \to e$ in $Q^*$
and $\delta_n \weakly \delta$ in $Q$, this yields
\begin{align*}
	\| \delta - \delta_n\|_Q^2
	&=
	\| \delta_n\|_Q^2
	- 
	2\left ( \delta,  \delta_n\right )_{Q}
	+
	\| \delta \|_Q^2
	\le
	\left \langle e_n, \delta_n \right  \rangle_Q
	- 
	2\left ( \delta,  \delta_n\right )_{Q}
    +
	\left \langle e, \delta  \right  \rangle_Q
	\to
	2  \langle e, \delta  \rangle_Q
	-
	2\|  \delta\|_{Q}^2
	=
	0
	.
\end{align*}
Thus, $\delta_n \to \delta$ in $Q$ and, as a consequence,
\begin{equation*}
\begin{aligned}
	&\frac{\left \| \Sigma_a (g + h_n) - \Sigma_a(g)- G_n h_n\right \|_Q}{\|h_n \|_{\MM([0,L])}} 
	\\
	&\qquad=
	\left \| 
	\frac{\Sigma_a (g + t_n e_n) - \Sigma_a(g)}{t_n} - G_n e_n
	\right \|_Q
	\\
	&\qquad\leq
	\left \| 
	\frac{\Sigma_a (g + t_n e_n) - \Sigma_a(g)}{t_n} - \Sigma_a'(g;e)
	\right \|_Q
	+
	\| \Sigma_a'(g;e)  - G_n e_n \|_Q
	\to 0.
\end{aligned}
\end{equation*}
This contradicts \eqref{eq:randomeq273hge8-2} and shows that
$\Sigma_a$ is Newton differentiable as desired. 
\end{proof}

Note that, in order to apply \cref{prop:energy_Newton_uni} in practice,  
one just requires a single element of the differential 
$\partial_B^{sw}\Sigma_a(g)$ for all $g \in \MM([0,L])$. The following lemma 
provides such an element and, in doing so, also prepares the analysis of the 
mapping $\Sigma_a^b$ associated with the bilateral case. 

\begin{lemma}[An element of  $\partial_B^{sw}\Sigma_a$]
\label{lem:elementBouligandDifferential}
Let $g \in Q^*$ be given 
 and let $G_g \in \LLLL(Q^*, Q)$ be  defined 
by
\[
	\delta = G_g(h)
	\qquad
	\iff
	\qquad
	\delta \in \QQ_a(g),\quad 
	\left (\delta, p \right )_Q = \left \langle h, p\right \rangle_{Q}
	\quad 
	\forall p \in \QQ_a(g)
\]
with $\QQ_a(g) := \{p \in Q \mid p(x) = 0 \text{ for all } x \in [0,L] \text{ with }  \Sigma_a(g)(x) = a(x)\}$.
Then,  $G_g \in \partial_B^{sw}\Sigma_a(g)$. 
\end{lemma}
\begin{proof}
Let $g \in Q^*$ be given. We define $A_g := \{x \in [0,L] \mid  \Sigma_a(g)(x) = a(x)\}$. As $A_g$ is closed, 
there exists a nonnegative measure $\mu \in \MM([0,L])$ whose support satisfies 
$\supp(\mu) = A_g$. Indeed, if $A_g = \emptyset$, we can simply choose $\mu = 0$,
and if $A_g \neq \emptyset$, we can choose a 
countable dense subset $\{x_k\}_{k=1}^\infty \subset A_g$ 
and define $\mu := \sum_{k=1}^\infty 2^{-k} \delta_{x_k}$, where $\delta_{x_k}$
denotes the Dirac measure supported at $x_k$.   
Define $g_n := g - \mu/n$, $n \in \N$. 
Then, it holds $g_n \in Q^*$ and $g_n \to g$ in $Q^*$. We claim that 
$g_n$ is an element of the set of Gâteaux points $\DD$ of $\Sigma_a$. To see this, we first note that 
the nonnegativity of $\mu$ and the equality $\supp(\mu) = A_g$ imply  that 
$\hat q := \Sigma_a(g) \in K_a$ satisfies 
\begin{equation*}
\begin{aligned}
\left ( \hat q, q-\hat q \right )_{Q}
- \left \langle g_n, q - \hat q \right \rangle_Q
&=
\left ( \hat q, q-\hat q \right )_{Q}
- \left \langle g, q - \hat q \right \rangle_Q
+
\frac1n 
\int_{[0,L]}  q - \hat q \, \dd \mu 
\\
&\geq
\frac1n 
\int_{[0,L]}  q - \hat q \, \dd \mu 
\\
&=
\frac1n 
\int_{[0,L]}  q - a\, \dd \mu 
\\
&\geq 0\qquad \forall q \in K_a. 
\end{aligned}
\end{equation*}
Thus, $\Sigma_a(g_n) = \Sigma_a(g) = \hat q$ for all $n \in \N$. Next, we show that,
for every $n \in \N$, we have 
$\TT(\hat q; K_a) \cap \ker(\eta_n) = \QQ_a(g)$,
where $\eta_n := \iota_R^{-1}(\hat q)- g_n \in Q^*$. 
To this end, we note that, for every $p \in \TT(\hat q; K_a)$, 
we have 
\[
	\left \langle \eta_n, p\right \rangle_Q
	=
	\left ( \hat q,  p \right )_{Q}
      - \left \langle g, p \right \rangle_Q
      +
      \frac1n 
       \int_{[0,L]}  p \, \dd \mu 
       \geq
       \frac1n 
        \int_{[0,L]}  p \, \dd \mu 
        \geq 0,
\]
where the first inequality follows from the VI for $\hat q =  \Sigma_a(g)$
and the second one from the fact that 
$p \in \TT(\hat q; K_a) = \closure_Q(\R^+(K_a - \hat q))$ implies 
$p \geq 0$ pointwise on $A_g$. The above and the identity 
$\supp(\mu) = A_g$ show that 
$p \in \TT(\hat q; K_a) \cap \ker(\eta_n)$ can only be true 
if $p \in \QQ_a(g)$, i.e.\ we have $ \TT(\hat q; K_a) \cap \ker(\eta_n) \subset \QQ_a(g)$.
Suppose next that $p \in \QQ_a(g)$ is given, i.e.\ that $p \in Q$ with $p = 0$ on $A_g$. 
Then, it is easy to check that our assumptions on $Q$ imply that 
$p_\varepsilon := \max(0, p - \varepsilon) + \min(0, p + \varepsilon)$
is in $Q$ for all $\varepsilon > 0$ and that $p_\varepsilon \to p$ in $C([0,L])$
and $p_\varepsilon \weakly p$ in $Q$ for $\varepsilon \to 0$. 
From the continuity of $p$ and 
$p = 0$ on $A_g$, we further obtain that 
$p_\varepsilon = 0$ holds in an open neighbourhood of 
$A_g$ for all $\varepsilon > 0$. 
As the function $\hat q - a$ 
has a positive minimum outside of 
this neighbourhood due to its continuity,
the definition of $A_g$
and the compactness of $[0,L]$
and since $p_\varepsilon$ is bounded, 
it follows that, for every $\varepsilon > 0$, there exists 
$c_\varepsilon > 0$ with $Q \ni \hat q \pm c_\varepsilon p_\varepsilon \geq a$. 
The latter  implies $\pm p_\varepsilon \in \TT(\hat q; K_a)$ and, with the lemma of Mazur, the convexity and closedness of $\TT(\hat q; K_a)$
and $p_\varepsilon \weakly p$ in $Q$, that $\pm p \in  \TT(\hat q; K_a)$.
Combining the latter with the VI for $\hat q = \Sigma_a(g_n)$ yields
$p \in \ker(\eta_n)$ and, consequently, $p \in \TT(\hat q; K_a) \cap \ker(\eta_n)$.
This shows that $ \QQ_a(g) \subset \TT(\hat q; K_a) \cap \ker(\eta_n)$ holds
and, by what we have shown before, that $\QQ_a(g) = \TT(\hat q; K_a) \cap \ker(\eta_n)$. 
As $\QQ_a(g)$ is a subspace, we may now deduce from 
\cref{lem:basic_props_sigma} that, for all $n \in \N$ and all $h \in Q^*$,
the directional derivative 
$\delta := \Sigma_a'(g_n;h) \in Q$ is uniquely characterised by the auxiliary problem 
\begin{equation*}
\text{Find }
	\delta \in \QQ_a(g)
	\text{ such that }
	\left (\delta, p \right )_{Q}
- \left \langle h, p \right \rangle_Q
= 0
~~\forall p  \in \QQ_a(g),
\end{equation*}
i.e.\ we have $g_n \in \DD$ with $\Sigma_a'(g_n;\cdot) = G_g \in \LLLL(Q^*, Q)$
for all $n$. In view of \cref{def:Bsubdifferential} and the convergence 
$g_n \to g$ in $Q^*$, this yields $G_g \in \partial_B^{sw}\Sigma_a(g)$ as desired
and completes the proof. 
\end{proof}

Next, we prove a lemma on the approximation of $\Sigma_a^b$ by means 
of semilinear problems.

\begin{lemma}
\label{lem:reg_lemma}
Let $g \in \MM([0,L])\subset Q^*$ be given. Then, the problem 
\begin{equation}
\label{eq:regular_aux_prob}
\text{Find }
\hat q_\varepsilon \in Q
\text{ such that }
\left ( \hat q_\varepsilon, q \right )_{Q}
- \left \langle g,q \right \rangle_Q
+
\int_0^L 
\frac1\varepsilon \min(0, \hat q_\varepsilon - a) q
+
\frac1\varepsilon \max(0,  \hat q_\varepsilon - b)q
\d x
= 0 
~\forall q \in Q
\end{equation}
has a unique solution $\hat q_\varepsilon$ for all $\varepsilon > 0$,
it holds $\hat q_\varepsilon \to \hat q := \Sigma_a^b(g)$ in $Q$ for $\varepsilon \to 0$ and there exist unique
measures $\hat \mu_a, \hat \mu^b \in \MM([0,L])$ such that 
\begin{equation}
\begin{gathered}
\label{eq:reg_probs}
\frac1\varepsilon \min(0, \hat q_\varepsilon - a) \weaklystar \hat \mu_a\text{ for }\varepsilon \to 0,
 \qquad 
\frac1\varepsilon \max(0,  \hat q_\varepsilon - b) \weaklystar \hat \mu^b\text{ for }\varepsilon \to 0,
\\
\left ( \hat q , q \right )_{Q}
=
\left \langle g - \hat \mu_a - \hat \mu^b, q \right \rangle_Q
\quad\, \forall q \in Q,
 \\
 \hat \mu_a \leq 0, \qquad 
 \supp(\hat \mu_a) \subset \{x \in [0,L] \mid \hat q(x) = a(x)\},
  \\
 \hat \mu^b \geq 0, \qquad 
 \supp(\hat \mu^b) \subset \{x \in [0,L] \mid \hat q(x) = b(x)\}.
\end{gathered}
\end{equation}
\end{lemma}

\begin{proof}
The existence of a unique solution $\hat q_\varepsilon $ of \eqref{eq:regular_aux_prob} 
for all $\varepsilon > 0$ follows from \cite[Theorem 4:3-3.1]{Rodrigues1987}. 
From the a-priori estimate in this theorem and the existence of the function 
$\tilde q \in Q$ in \cref{ass:standing:appendix_A},
we also obtain that there exists a constant $C_1 > 0$ independent of $\varepsilon > 0$ 
with 
\begin{equation}
\label{eq:apriori_estimate_s2s}
	\|\hat q_\varepsilon\|_Q \leq C_1  + \|g\|_{Q^*}\qquad \forall \varepsilon > 0. 
\end{equation}
Consider now the function
	\begin{equation}
	\label{eq:reg_sign_def}
		\rho_s(x)
		:=
      \max \left ( -1, \min \left (1, \frac{x}{s} \right ) \right )
 		=
		\begin{cases}
			-1 & \text{for } x \le -s, \\
			\frac xs & \text{for } x \in (- s,s), \\
			1 & \text{for } x \ge s, \\
		\end{cases}
		\qquad s > 0.
	\end{equation}
Then, 
$\rho_s$
is nondecreasing, continuous and piecewise affine-linear 
with $\rho_s(0) = 0$
and it holds
$\rho_s(\hat q_\varepsilon - \tilde q) \in Q$
for all $s>0$ by our assumptions. 
Choosing $\rho_s(\hat q_\varepsilon - \tilde q)$ as the test function in \eqref{eq:regular_aux_prob}  
and exploiting \cref{ass:standing:appendix_A} yields
\begin{equation*}
\begin{aligned}
&\int_0^L 
\frac1\varepsilon 
\big (
\min(0, \hat q_\varepsilon - a) 
-
\min(0, \tilde q - a) 
\big )
\rho_s(\hat q_\varepsilon - \tilde q) 
+
\frac1\varepsilon
\big  ( 
\max(0,  \hat q_\varepsilon - b) - \max(0,  \tilde q - b)
\big )
\rho_s(\hat q_\varepsilon - \tilde q)  
\d x
\\
&=
\int_0^L 
\frac1\varepsilon 
\min(0, \hat q_\varepsilon - a) 
\rho_s(\hat q_\varepsilon - \tilde q) 
+
\frac1\varepsilon
\max(0,  \hat q_\varepsilon - b)
\rho_s(\hat q_\varepsilon - \tilde q)  
\d x
\\
&=
- \left ( \hat q_\varepsilon - \tilde q, \rho_s(\hat q_\varepsilon - \tilde q) \right )_{Q}
- \left ( \tilde q, \rho_s(\hat q_\varepsilon - \tilde q) \right )_{Q}
+ \left \langle g, \rho_s(\hat q_\varepsilon - \tilde q) \right \rangle_Q
\\
&\leq
\|g - \tilde \mu\|_{\MM([0,L])} =: C_2 < \infty.
\end{aligned}
\end{equation*}
Passing to the limit $0 < s \to 0$
in the above provides us with the bound
\begin{equation}
\label{eq:randomeq28hei3n}
\int_0^L 
\frac1\varepsilon 
\big |
\min(0, \hat q_\varepsilon - a) 
\big |
+
\frac1\varepsilon
\big  |
\max(0,  \hat q_\varepsilon - b)
\big | \d x 
\leq 
C_2. 
\end{equation}
Suppose now that 
$\{\varepsilon_n\} \subset \R^+$ is a sequence satisfying 
$\varepsilon_n \to 0$. Then \eqref{eq:apriori_estimate_s2s}
implies that, by passing over to a subsequence (still denoted the same),
we can achieve that $\hat q_{\varepsilon_n} \weakly \bar q$ holds 
for some $\bar q \in Q$.
From \eqref{eq:randomeq28hei3n} and $Q \compactly C([0,L])$,
we get that 
$\bar q \in K_a^b$ holds, and from \eqref{eq:regular_aux_prob} and the same 
arguments as above, that
\begin{equation}
\label{eq:randomeq36345}
\left ( \hat q_{\varepsilon_n}, q - \hat q_{\varepsilon_n} \right )_{Q}
- \left \langle g,  q - \hat q_{\varepsilon_n} \right \rangle_Q
\geq 0
\quad\, \forall q \in K_a^b.
\end{equation}
Passing to the limit in \eqref{eq:randomeq36345}, using the weak lower semicontinuity of $\|\cdot\|_Q^2$,
yields $\bar q = \hat q = \Sigma_a^b(g)$. In combination with a standard subsequence--subsequence
argument, this shows that $\hat q_\varepsilon \weakly \hat q$ holds in $Q$ for 
$0 < \varepsilon \to 0$. 
If we define $\hat \mu_a^\varepsilon := \frac1\varepsilon \min(0, \hat q_\varepsilon - a)$
and $\hat \mu_b^\varepsilon :=  \frac1\varepsilon \max(0,  \hat q_\varepsilon - b)$,
then we further obtain from \eqref{eq:randomeq28hei3n}  that 
$\|\hat \mu_a^\varepsilon\|_{L^1(0,L)} + \|\hat \mu_b^\varepsilon\|_{L^1(0,L)} \leq C_2$.
By passing over to a sequence $0 < \varepsilon_n \to 0$
with \smash{$\hat \mu_a^{\varepsilon_n} \weaklystar \hat \mu_a$}
and \smash{$\hat \mu_b^{\varepsilon_n} \weaklystar \hat \mu^b$} 
in $\MM([0,L])$, by using \eqref{eq:regular_aux_prob},
by exploiting the compactness of the embedding $\MM([0,L]) \hookrightarrow Q^*$
and 
by invoking the Riesz representation theorem, it follows that 
$\hat q_{\varepsilon_n} \to \hat q$ in $Q$. Again a subsequence--subsequence argument 
now provides us with $\hat q_\varepsilon \to \hat q$ in $Q$. 
Note that the definitions of $\hat \mu_a^\varepsilon$ and $\hat \mu_b^\varepsilon$,
the equation \eqref{eq:regular_aux_prob}, 
the convergence $\hat q_\varepsilon \to \hat q$ and the embedding $Q \hookrightarrow C([0,L])$
imply that, in the above situation, the weak-star limits 
$ \hat \mu_a,  \hat \mu^b$ satisfy all of the conditions in \eqref{eq:reg_probs}
along the sequence $\{\varepsilon_n\}$. From the variational identity in  \eqref{eq:reg_probs},
the density of $Q$ in $C([0,L])$, the fact that the supports of the measures 
$\hat \mu_a$ and $\hat \mu^b$ have a positive distance from each other
owing to our assumptions on $a$ and $b$ 
and the uniqueness of the Jordan decomposition, we obtain that 
$\hat \mu_a$ and $\hat \mu^b$ are uniquely determined.
Again a subsequence--subsequence argument now provides us with 
the convergences  $\hat \mu_a^{\varepsilon } \weaklystar \hat \mu_a$ 
and $\hat \mu_b^{\varepsilon } \weaklystar \hat \mu^b$ 
in $\MM([0,L])$ for $\varepsilon \to 0$. Collecting everything establishes the assertion of the lemma. 
\end{proof}

\Cref{lem:reg_lemma} holds completely analogously 
for the unilateral cases. We omit stating 
this as a separate result here. For later use, we introduce:

\begin{definition}[Multiplier Maps]
\label{def_mult_maps}
We define 
\begin{equation*}
\begin{aligned}
M_a\colon \MM([0,L]) \to \MM([0,L]),
\qquad g \mapsto g - \iota_R^{-1}(\Sigma_a(g)),
\\
M^b\colon \MM([0,L]) \to \MM([0,L]),
\qquad g \mapsto g -  \iota_R^{-1}(\Sigma^b(g)),
\\
M_a^b\colon \MM([0,L]) \to \MM([0,L]),
\qquad g \mapsto g -  \iota_R^{-1}(\Sigma_a^b(g)).
\end{aligned}
\end{equation*}
For the bilateral case, we further define 
\[
	[M_a^b]_\pm \colon \MM([0,L]) \to \MM([0,L]),
	\qquad 
	g \mapsto [M_a^b(g)]_\pm,
\]
where $[M_a^b(g)]_\pm$ denotes the Jordan decomposition 
of $M_a^b(g)\in \MM([0,L])$; cf.\ \cref{lem:reg_lemma}. 
\end{definition}

Note that, by definition, $[M_a^b(g)]_+$ and $[M_a^b(g)]_-$ are
identical to the measures 
$\hat \mu^b$ and $\hat \mu_a$ in the system \eqref{eq:reg_probs}.
Analogously, $M_a$ and $M^b$
correspond to the multipliers that appear in the complementarity systems 
of the unilateral VIs which involve only the lower 
bound $a$, respectively, the upper bound $b$. 
From \cref{lem:reg_lemma}, we obtain:

\begin{lemma}[Lipschitz continuity of the multiplier maps]
\label{lem:LipschitzMMM}
The maps $M_a$, $M^b$, $M_a^b$, $[M_a^b]_+$ and $[M_a^b]_-$
are nonexpansive
as functions from $\MM([0,L])$ to $\MM([0,L])$.
\end{lemma}
\begin{proof}
We prove the assertions for the bilateral case. 
The unilateral cases are analogous.  
Let $g_1, g_2 \in \MM([0,L])$ with solutions 
$\hat q_1 := \Sigma_a^b(g_1), \hat q_2 := \Sigma_a^b(g_2) \in Q$
be given, and let $\hat q_\varepsilon^1, \hat q_\varepsilon^2$ be defined 
by \eqref{eq:regular_aux_prob}. 
Denote by $\rho_s$ the function in \eqref{eq:reg_sign_def}. Then, 
choosing $\rho_s(\hat q_\varepsilon^1 - \hat q_\varepsilon^2) \in Q$
as the test function in the variational identities for 
$\hat q_\varepsilon^1$ and $\hat q_\varepsilon^2$
and subtracting the resulting equations yields 
\begin{equation*}
\begin{aligned}
&\int_0^L 
\left |\left ( \frac1\varepsilon \min(0, \hat q_\varepsilon^1 - a) 
-
\frac1\varepsilon \min(0, \hat q_\varepsilon^2 - a)\right )\rho_s(\hat q_\varepsilon^1 - \hat q_\varepsilon^2)
\right | 
\\
&\qquad 
+
\left |
\left (
\frac1\varepsilon \max(0,  \hat q_\varepsilon^1 - b)
-
\frac1\varepsilon \max(0,  \hat q_\varepsilon^2 - b)
\right )\rho_s(\hat q_\varepsilon^1 - \hat q_\varepsilon^2)
\right |
\d x
\\
&=
\int_0^L 
\left ( \frac1\varepsilon \min(0, \hat q_\varepsilon^1 - a) 
-
\frac1\varepsilon \min(0, \hat q_\varepsilon^2 - a)\right )\rho_s(\hat q_\varepsilon^1 - \hat q_\varepsilon^2)
\\
&\qquad 
+
\left (
\frac1\varepsilon \max(0,  \hat q_\varepsilon^1 - b)
-
\frac1\varepsilon \max(0,  \hat q_\varepsilon^2 - b)
\right )\rho_s(\hat q_\varepsilon^1 - \hat q_\varepsilon^2)
\d x
\\
& = -
\left ( \hat q_\varepsilon^1 -  \hat q_\varepsilon^2, \rho_s(\hat q_\varepsilon^1 - \hat q_\varepsilon^2) \right )_{Q}
+ \left \langle g_1 - g_2, \rho_s(\hat q_\varepsilon^1 - \hat q_\varepsilon^2) \right \rangle_Q
\\
&\leq
\|g_1 - g_2\|_{\MM([0,L])}.
\end{aligned}
\end{equation*}
Letting $0 < s \to 0$ in the above provides us with the estimate
\begin{equation*}
\begin{aligned}
& 
\left \| \frac1\varepsilon \min(0, \hat q_\varepsilon^1 - a) 
-
\frac1\varepsilon \min(0, \hat q_\varepsilon^2 - a)
\right \|_{L^1(0,L)} 
+
\left \|
\frac1\varepsilon \max(0,  \hat q_\varepsilon^1 - b)
-
\frac1\varepsilon \max(0,  \hat q_\varepsilon^2 - b)
\right \|_{L^1(0,L)} 
\\
&\hspace{8cm}\leq
\|g_1 - g_2\|_{\MM([0,L])}.
\end{aligned}
\end{equation*}
Due to the weak-star lower semicontinuity of the $\MM([0,L])$-norm,
\cref{lem:reg_lemma}
and the definition of $[M_a^b]_\pm$, passing to the limit $0 < \varepsilon \to 0$ in the above gives
\[
\left \|
[M_a^b]_-(g_1) - [M_a^b]_-(g_2)
\right \|_{\MM([0,L])} 
+
\left \|
[M_a^b]_+(g_1) - [M_a^b]_+(g_2)
\right \|_{\MM([0,L])} 
\leq
\|g_1 - g_2\|_{\MM([0,L])}.
\]
This establishes the assertion for $[M_a^b]_-$, $[M_a^b]_+$, and $M_a^b$ and completes the proof. 
\end{proof}

Note that, owing to \cref{prop:energy_Newton_uni} and the chain rule,
we know that the maps $M_a$ and $M^b$ are not only one-Lipschitz but also Newton differentiable:

\begin{lemma}[Newton derivatives for the unilateral multiplier maps]
\label{lem:Newton_Ma_Mb}
The multiplier mappings $M_a, M^b\colon \MM([0,L]) \to \MM([0,L])$ are Newton differentiable as 
functions from $ \MM([0,L])$ to $Q^*$ when endowed with the Newton derivatives 
$DM_a, DM^b \colon \MM([0,L]) \to \LLLL(Q^*, Q^*)$ given by 
\[
		DM_a(g) := \Id - \iota_R^{-1} \circ D\Sigma_a(g),
		\qquad 
		DM^b(g) := \Id - \iota_R^{-1} \circ D\Sigma^b(g)
		\qquad
		\forall g \in \MM([0,L]).
\]
Here and in what follows, 
the functions 
$D\Sigma_a, D\Sigma^b\colon \MM([0,L]) \to \LLLL(Q^*, Q)$ are defined by 
\[
	\delta = D\Sigma_a(g)h
	\qquad
	\iff
	\qquad
	\delta \in \QQ_a(g),\quad 
	\left (\delta, p \right )_Q = \left \langle h, p\right \rangle_{Q}
	\quad 
	\forall p \in \QQ_a(g)
\]
with $\QQ_a(g) := \{p \in Q \mid p(x) = 0 \text{ for all } x \in [0,L] \text{ with }  \Sigma_a(g)(x) = a(x)\}$
and
\[
	\delta = D\Sigma^b(g)h
	\qquad
	\iff
	\qquad
	\delta \in \QQ^b(g),\quad 
	\left (\delta, p \right )_Q = \left \langle h, p\right \rangle_{Q}
	\quad 
	\forall p \in \QQ^b(g)
\]
with $\QQ^b(g) := \{p \in Q \mid p(x) = 0 \text{ for all } x \in [0,L] \text{ with }  \Sigma^b(g)(x) = b(x)\}$,
respectively.
\end{lemma}
\begin{proof}
This follows straightforwardly from 
\cref{prop:energy_Newton_uni} and \cref{lem:elementBouligandDifferential}. 
\end{proof}

To establish the Newton differentiability of $M_a^b$ and, in doing so, 
the one of $\Sigma_a^b$, we argue along the lines of 
\cite[Section 6]{ChristofWachsmuthBiObstacle}, taking some shortcuts that 
are made possible by the embedding $Q \hookrightarrow C([0,L])$
and the Lipschitz continuity properties in \cref{lem:LipschitzMMM}. 
The main idea is to use a reasoning as in the proof of the implicit function theorem. 
To this end, we introduce:

\begin{definition}
\label{def:Phietc}
We define 
(using the notation
from \cref{def_mult_maps} and \cref{lem:Newton_Ma_Mb}):
\begin{equation*}
\begin{aligned}
&\begin{aligned}
&\HH 
:=
Q^*
\times
Q^*
,
&&
\|(g_1, g_2)\|_\HH := \left (\|g_1\|_{Q^*}^2 + \|g_2\|_{Q^*}^2 \right )^{1/2},
\\
&\Phi\colon \MM([0,L])
\to \MM([0,L])^2,\qquad
&&
\Phi(g)
:=
\begin{pmatrix}
\Phi_a(g)
\\
\Phi^b(g)
\end{pmatrix}
:=
\begin{pmatrix}
  [M_a^b(g)]_-
\\
  [M_a^b(g)]_+
\end{pmatrix},
\\
&F\colon \MM([0,L])^3
\to
\HH,
&&
F( \mu_a,\mu^b,g)
:=
\begin{pmatrix}
\mu_a
-
M_a(g - \mu^b)
\\
\mu^b
-
M^b(g - \mu_a)
\end{pmatrix},
\end{aligned}
\\
&\D F\colon \MM([0,L])^3 \to \LLLL(\HH, \HH) \times \LLLL(Q^*,\HH),
\\
&\D F(\mu_a, \mu^b, g) 
	:=
		\left [
		 \begin{pmatrix}
	\Id   & \D M_a(g - \mu^b)
		\\
	\D M^b(g - \mu_a)  	& \Id
	\end{pmatrix},	
	\begin{pmatrix}
		- \D M_a(g - \mu^b)
		\\
		- \D M^b(g - \mu_a) 
	\end{pmatrix}
	\right ].
\end{aligned}
\end{equation*}
\end{definition}

We note:

\begin{lemma}
\label{lem:NewtonDiffF}
The function $F\colon \MM([0,L])^3
\to
\HH$ is Newton differentiable with derivative $DF$.
\end{lemma}

\begin{proof}
This follows immediately from \cref{lem:Newton_Ma_Mb}.
\end{proof}

\begin{lemma}
\label{lem:FP_equation}
For every $g \in \MM([0,L])$, there holds $F(\Phi(g), g) = 0 \in \HH$. 
\end{lemma}

\begin{proof}
This follows immediately from the definition of $\Phi$ and the system \eqref{eq:reg_probs}.
\end{proof}

To argue along the lines of the implicit function theorem, we require that 
the first component of the operator $DF$ is locally uniformly invertible. 
To prove this, we need some lemmas:

\begin{lemma}
\label{lem:projections}
Let $g \in \MM([0,L])$. Set
$G_a := \D M_a(g - \Phi^b(g))$ and $G^b := \D M^b(g - \Phi_a(g))$,
and let  $\QQ_a := \QQ_a(g - \Phi^b(g))$ and $\QQ^b := \QQ^b(g - \Phi_a(g))$ 
be as in \cref{lem:Newton_Ma_Mb}.
Define 
$W_a(g) := \left [ \iota_R^{-1}(\QQ_a )\right ]^\perp$
and 
$W^b(g) := \left [ \iota_R^{-1}(\QQ^b)\right ]^\perp$,
where the orthogonal complement is taken w.r.t.\ the 
canonical inner product $(w_1, w_2)_{Q^*} := ( \iota_R(w_1), \iota_R(w_2))_Q$
on $Q^*$. Then, it holds 
\[
	G_a = \proj_{W_a(g)}\qquad \text{and}\qquad G^b = \proj_{W^b(g)},
\]
where $\proj$ denotes an orthogonal projection in $Q^*$ w.r.t.\ $(\cdot, \cdot)_{Q^*}$. 
\end{lemma}

\begin{proof}
From the definition of $\D\Sigma_a(g - \Phi^b(g))$, 
we obtain that $\D\Sigma_a(g - \Phi^b(g))\iota_R^{-1} = \proj_{\QQ_a}$.
This implies that $\iota_R^{-1}\D\Sigma_a(g - \Phi^b(g))$ is the 
$Q^*$-orthogonal projection onto $\iota_R^{-1}(\QQ_a)$,
which in turn 
entails that $G_a = \D M_a(g - \Phi^b(g)) = \Id - \iota_R^{-1}\D\Sigma_a(g - \Phi^b(g))$
is the orthogonal projection onto $W_a(g) = \left [\iota_R^{-1}(\QQ_a) \right ]^\perp$ in $Q^*$.
The assertion for $G^b$ is proved analogously. 
\end{proof}

\begin{lemma}
\label{lem:stability contact sets}
Let $g \in Q^*$ be given. Then, there exist
a constant
$r > 0$
and nonempty disjoint compact sets $Z_a, Z^b \subset [0,L]$ such that
\begin{align*}
\{\Sigma_a^b(g + h) = a\} \subset Z_a
\qquad
\text{and}
\qquad
\{\Sigma_a^b(g + h) = b\} \subset Z^b
\end{align*}
holds for all $h \in Q^*$ with $\|h\|_{Q^*} \leq r$.
\end{lemma}

\begin{proof}
This lemma can be proved completely analogously to 
\cite[Lemma 6.14]{ChristofWachsmuthBiObstacle},
using the embedding $Q \hookrightarrow C([0,L])$,
our assumption $a + \varepsilon \leq b $,
and the Lipschitz continuity of $\Sigma_a^b\colon Q^* \to Q$.
\end{proof}

\begin{lemma}
\label{lem:sine_prototype}
Let $g \in \MM([0,L])$ be given and let $W_a(\cdot)$, $W^b(\cdot)$ be defined
as in \cref{lem:projections}. Then, there exist constants $r>0$, $C > 1$ such that,
for all $h \in \MM([0,L])$ with $\|h\|_{Q^*} \leq r$, it holds 
\begin{equation}
\label{eq:randomeq36363}
C \|v + w\|_{Q^*} \geq \|v\|_{Q^*}
\qquad
\forall v \in W_a(g + h)
\quad
\forall w \in W^b(g + h).
\end{equation}
\end{lemma}

\begin{proof}
Let $g \in \MM([0,L])$ be given, and let $r>0$ and $Z_a, Z^b$ be chosen 
for $g$ as in \cref{lem:stability contact sets}. Choose a function 
$\psi \in C^\infty([0,L])$ such that $\psi = 1$ holds in $Z_a$
and $\psi = 0$ holds in $Z^b$. Recall that our assumptions 
imply that $Q \ni q \mapsto \psi q \in Q$ is a well-defined element 
of $\LLLL(Q, Q)$. We define 
\begin{equation}
\label{eq:Jdef}
	\JJ\colon Q^* \to Q^*,
	\qquad
	\left \langle \JJ(w), q\right \rangle_Q
	:=
	\left \langle w, \psi q\right \rangle_Q
	\qquad
	\forall w \in Q^*
	\quad
	\forall q \in Q.
\end{equation}
It is easy to check that $\JJ \in \LLLL(Q^*, Q^*)$. Consider now some $h \in \MM([0,L])$
with $\|h\|_{Q^*} \leq r$. Then,  
$\{\Sigma_a^b(g + h) = a\} \subset Z_a$ and $\{\Sigma_a^b(g + h) = b\} \subset Z^b$ hold.
In combination with the properties of $\psi$ and the identity 
\begin{align*}
&\QQ_a(g + h - \Phi^b(g + h)) 
\\
&= 
\{p \in Q \mid p(x) = 0 \text{ for all } x \in [0,L] \text{ with }  \Sigma_a(g + h - \Phi^b(g + h))(x) = a(x)\}
\\
&=
\{p \in Q \mid p(x) = 0 \text{ for all } x \in [0,L] \text{ with }  \Sigma_a^b(g + h)(x) = a(x)\}
\end{align*}
obtained from the definition of $\QQ_a(\cdot)$ in \cref{lem:Newton_Ma_Mb} and 
\cref{lem:FP_equation}, this implies in particular that 
$(1 - \psi)q \in \QQ_a(g + h - \Phi^b(g + h)) $ for all $q \in Q$.
Analogously, we get $\psi q \in \QQ^b(g + h - \Phi_a(g + h))$
for all $q \in Q$.
For arbitrary $v \in W_a(g + h)$ and $w \in W^b(g + h)$, it now follows
\begin{equation*}
\begin{aligned}
\left \langle \JJ(v + w), q\right \rangle_{Q}
&=
\left \langle v + w, \psi q  \right \rangle_{Q}
&&\text{by \eqref{eq:Jdef}}
\\
&=
\left ( v + w, \iota_R^{-1} (\psi q )  \right )_{Q^*}
&&\text{by the definitions of $\iota_R$ and $(\cdot, \cdot)_{Q^*}$}
\\
&=
\left ( v, \iota_R^{-1}(\psi q ) \right )_{Q^*}
&&\text{by $\psi q \in \QQ^b(g + h - \Phi_a(g + h))$, the definition of $W^b$}
\\
&=
\left ( v, \iota_R^{-1}(q)  \right )_{Q^*}
&&\text{by $(1-\psi) q \in \QQ_a(g + h - \Phi^b(g + h))$, the definition of $W_a$}
\\
&=
\left \langle v, q \right   \rangle_{Q}
&&\text{by the definition of $\iota_R$},
\end{aligned}
\end{equation*}
i.e.\ we have $\JJ(v + w) = v$. This yields
\[
	\|v\|_{Q^*}
	=
	\|\JJ(v + w)\|_{Q^*}
	\leq
	\|\JJ\|_{\LLLL(Q^*, Q^*)}
	\|v + w\|_{Q^*}
	\qquad
\forall v \in W_a(g + h)
\quad
\forall w \in W^b(g + h),
\]
which implies the assertion of the lemma. 
\end{proof}

\begin{lemma}
\label{lem:uniform_invertibility}
Let $g \in \MM([0,L])$ be given. Then, there exist constants $r, C>0$ such that,
for every $h \in \MM([0,L])$ with $\|h\|_{Q^*} \leq r$, 
the first component of $DF(\Phi_a(g + h), \Phi^b(g + h), g + h)$ is invertible with 
\[
	\left \| \begin{pmatrix}
	\Id   & \D M_a(g + h - \Phi^b(g + h))
		\\
	\D M^b(g + h - \Phi_a(g + h))  	& \Id
	\end{pmatrix}^{-1} \right \|_{\LLLL(\HH, \HH)}
	\leq C.
\]
\end{lemma}
\begin{proof}
Let $g \in \MM([0,L])$ be given and let 
$C$ and $r$ be chosen such that the assertion of \cref{lem:sine_prototype}
holds.
Suppose further that $h \in \MM([0,L])$ satisfies $\|h\|_{Q^*} \leq r$.
 Then, it follows from \eqref{eq:randomeq36363} and 
 \cite[Lemma 5.4]{ChristofWachsmuthBiObstacle}
 that the sine of the minimal angle between the spaces $W_a(g + h)$ and $W^b(g + h)$ in $Q^*$
 is bounded from below by  $1/C > 0$. In view of 
 \cite[Theorem 5.5]{ChristofWachsmuthBiObstacle}, this entails that the operator
 \[
 			R_1 := \begin{pmatrix}
	\Id   & \proj_{W_a(g+h)}
		\\
	\proj_{W^b(g+h)} 	& \Id
	\end{pmatrix}
	\in \LLLL(\HH, \HH)
 \]
 is continuously invertible with $\|R_1^{-1}\|_{\LLLL(\HH, \HH)} \leq 4/(1 - \sqrt{1 - C^{-2}})$. 
 As $R_1$ is precisely the first component of $DF(\Phi_a(g + h), \Phi^b(g + h), g + h)$
 by \cref{lem:projections}, this proves the assertion. 
\end{proof}

Using the uniform invertibility result in \cref{lem:uniform_invertibility}, we can now prove:

\begin{proposition}[Newton differentiability of $\Phi$]
\label{prop:semismoothPhi}%
The function 
$\Phi$
is Newton differentiable 
as a map from $\MM([0,L])$
to $\HH$ with the derivative 
$\D\Phi\colon \MM([0,L]) \to \LLLL(\MM([0,L]), \HH)$ given by 
\begin{equation}
\label{eq:DPhiDef}
\D\Phi(g)
:=
 - 
 \begin{pmatrix}
	\Id   & \D M_a(g - \Phi^b(g))
		\\
	\D M^b(g - \Phi_a(g))  	& \Id
	\end{pmatrix}^{-1}
	\begin{pmatrix}
		- \D M_a(g - \Phi^b(g))
		\\
		- \D M^b(g - \Phi_a(g)) 
	\end{pmatrix}.
\end{equation}
\end{proposition}
\begin{proof}
Let $g \in \MM([0,L])$ be given,
and let $r,C > 0$ be as in the assertion of \cref{lem:uniform_invertibility}.
Suppose further that $h \in \MM([0,L])$ satisfying 
$0 < \|h\|_{Q^*} \leq r$ is given and denote 
the two components of $\GG := DF(\Phi_a(g + h), \Phi^b(g + h), g + h)$
by $\GG_1 \in \LLLL(\HH, \HH)$ and $\GG_2 \in \LLLL(Q^*, \HH)$, respectively. 
Then, it follows from \cref{lem:uniform_invertibility,lem:FP_equation,lem:LipschitzMMM}
and the definition of $\Phi$ that
\begin{equation}
\label{eq:longequation22}
\begin{aligned}
	&
	\frac{\left \| \Phi(g + h) - \Phi(g) +\GG_1^{-1} \GG_2 h \right \|_\HH}{\left \|h\right \|_{\MM([0,L])}}
	\\&
	\leq
	\|\GG_1^{-1}\|_{\LLLL(\HH, \HH)}
	\frac{\left \| \GG_1 \left (\Phi(g + h) - \Phi(g) \right )+ \GG_2 h \right \|_\HH}{\left \|h\right \|_{\MM([0,L])}}
	\\
	&\leq 
	C
	 \frac{	 \left \|
	  (\GG_1~~\GG_2)
	 	\begin{pmatrix}
		\Phi(g + h) - \Phi(g) 
		\\
		h
	\end{pmatrix}\right \|_\HH}{\left \|h\right \|_{\MM([0,L])}}
	\\&
	=
	C
		 \frac{ \left \|
	 F(\Phi(g + h), g+h) - F(\Phi(g), g) - 
	 \GG
	 	\begin{pmatrix}
		\Phi(g+ h) - \Phi(g) 
		\\
		h
	\end{pmatrix}\right \|_\HH}{ \left \|
	 	\begin{pmatrix}
		\Phi(g + h) - \Phi(g) 
		\\
		h
	\end{pmatrix}\right \|_{\MM([0,L])^3}}
	\frac{\left \|
	 	\begin{pmatrix}
		\Phi(g + h) - \Phi(g) 
		\\
		h
	\end{pmatrix}\right \|_{\MM([0,L])^3}}
	{
	\left \|h\right \|_{\MM([0,L])}}
	\\
	&\leq
	3 C 
	 \frac{\left \|
	 F(\Phi(g + h), g+h) - F(\Phi(g), g) - 
	 \GG
	 	\begin{pmatrix}
		\Phi(g + h) - \Phi(g) 
		\\
		h
	\end{pmatrix}\right \|_\HH}{ \left \|
	 	\begin{pmatrix}
		\Phi(g + h) - \Phi(g) 
		\\
		h
	\end{pmatrix}\right \|_{\MM([0,L])^3}}
	.
\end{aligned}
\end{equation}
Owing to \cref{lem:NewtonDiffF} and again the Lipschitz continuity in \cref{lem:LipschitzMMM},
the right-hand side of the above estimate goes to zero for $0 < \|h\|_{\MM([0,L])} \to 0$.
This shows that $\Phi\colon \MM([0,L]) \to \HH$ is Newton differentiable 
with the asserted derivative and completes the proof. 
\end{proof}

\begin{remark}
Because of the embedding $Q \hookrightarrow C([0,L])$ and the 
Lipschitz continuity of $\Phi$ as a function 
from $\MM([0,L])$ to $\MM([0,L])^2$,
we do not require  harmonic analysis results 
in the proof of \cref{prop:semismoothPhi}
in contrast to the proof of \cite[Proposition 6.19]{ChristofWachsmuthBiObstacle}.
This makes it possible to obtain \cref{prop:semismoothPhi} in a general setting 
which is independent of the precise structure of the inner product 
$(\cdot, \cdot)_Q$ and to show that, in the one-dimensional setting, 
we can get Newton differentiability even for $\MM([0,L])$ inputs
(and not only for $L^1(0,L)$-inputs, as remarked in \cite{ChristofWachsmuthBiObstacle}).
\end{remark}

By rewriting \cref{prop:semismoothPhi}, we can now establish the main 
result of this appendix: 

\begin{theorem}[Newton Differentiability for $\Sigma_a^b$]
\label{th_main_appendix}
The function $\Sigma_a^b$ is Newton differentiable 
as a map $\Sigma_a^b\colon \MM([0,L]) \to Q$
when endowed with the 
derivative 
$ \D \Sigma_a^b\colon \MM([0,L]) \to \LLLL(\MM([0,L]), Q)$
defined by 
\begin{equation}
\label{eq:DSigma_def}
	\delta = D\Sigma_a^b(g)h
	\qquad
	\iff
	\qquad
	\delta \in \QQ_a^b(g),\quad 
	\left (\delta, p \right )_Q = \left \langle h, p\right \rangle_{Q}
	\quad 
	\forall p \in \QQ_a^b(g)
\end{equation}
with $\QQ_a^b(g) := \{p \in Q \mid p(x) = 0 
\text{ for all } x \in [0,L] \text{ with }  \Sigma_a^b(g)(x) = a(x) \text{ or } \Sigma_a^b(g)(x) = b(x)\}$.
\end{theorem}
\begin{proof}
From the definitions of $\Phi$ and $M_a^b$, we obtain that
$\Sigma_a^b = \iota_R( \Id - M_a^b) = \iota_R[\Id - (\Phi_a + \Phi^b)]$.
In view of \cref{prop:semismoothPhi}, 
this implies by the chain and sum rule for Newton derivatives that 
$\Sigma_a^b$ is Newton differentiable as a function from $\MM([0,L])$
to $Q$ with derivative $\iota_R[\Id -  (\Id~\Id )\D \Phi  ]$.
It remains to check that the last expression gives rise to 
the operator $\D\Sigma_a^b$ in \eqref{eq:DSigma_def}. 
To this end, we suppose that $g \in \MM([0,L])$ is given. 
From \cref{prop:semismoothPhi} and \cref{lem:projections}, we obtain that 
\begin{equation}
\label{eq:randomeq373h}
(\Id~\Id )\D \Phi
=
-
(\Id~\Id )
 \begin{pmatrix}
	\Id   & \proj_{W_a(g)}
		\\
	\proj_{W^b(g)} 	& \Id
	\end{pmatrix}^{-1}
	\begin{pmatrix}
		- \proj_{W_a(g)}
		\\
		- \proj_{W^b(g)}
	\end{pmatrix},
\end{equation}
where $W_a(g)$ and $W^b(g)$ are defined as in \cref{lem:projections}.
From \cite[Theorem 5.5]{ChristofWachsmuthBiObstacle} and again 
\cref{lem:sine_prototype}, we obtain
that $W_a(g) + W^b(g)$ is a closed subspace of $Q^*$
and that 
the right-hand side of \eqref{eq:randomeq373h}
is identical to the orthogonal projection $ \proj_{W_a(g) + W^b(g)}$.
If we denote by $\QQ_a$ and $\QQ^b$ the subspaces  of $Q$ from \cref{lem:projections},
then it follows further that 
\begin{equation*}
\begin{aligned}
W_a(g) + W^b(g) &= \left [ \iota_R^{-1}(\QQ_a)\right ]^\perp + \left [ \iota_R^{-1}(\QQ^b)\right ]^\perp
=  \left [ \iota_R^{-1}(\QQ_a \cap \QQ^b)\right ]^\perp
=
\left [ \iota_R^{-1}(\QQ_a^b(g))\right ]^\perp,
\end{aligned}
\end{equation*}
where the last equation follows from the definitions of the subspaces $\QQ_a$ and $\QQ^b$
and the identities 
$\Sigma_a (g - \Phi^b(g)) = \Sigma_a^b(g)$
and 
$\Sigma^b (g - \Phi_a(g)) = \Sigma_a^b(g)$
obtained from the system \eqref{eq:reg_probs}. 
Combined, the above shows that 
\[
	 \iota_R\left [\Id -  (\Id~\Id )\D \Phi \right ] 
	=
	 \iota_R\left [\Id -  \proj_{W_a(g) + W^b(g)} \right ] 
	 =
	 \iota_R\left [\proj_{\iota_R^{-1}(\QQ_a^b(g))}  \right ]
  =
	 \left [\proj_{\QQ_a^b(g)}  \right ]\iota_R.
\]
The last expression
reproduces precisely the operator $\D \Sigma_a^b$ when written down variationally. 
Thus, $\D\Sigma_a^b = \iota_R[\Id -  (\Id~\Id )\D \Phi  ]$ 
and the assertion of the theorem is proved. 
\end{proof}

\bibliographystyle{abbrv} %%%% .BST file

\bibliography{references} %%%%% .Bib file

\end{document}